\documentclass[12pt]{article}
\usepackage{amsmath,amssymb,amsthm,mathrsfs,stmaryrd}
\usepackage{url,enumerate,marginnote,booktabs,graphicx}
\usepackage{xcolor} 
\usepackage{float}
\usepackage{indentfirst}
\usepackage[numbers,sort&compress,square]{natbib}
\usepackage[textwidth=16.5cm, textheight=22.5cm, top=2.54cm,left=2.24cm]{geometry}
\usepackage[colorlinks,linkcolor=blue,citecolor=blue]{hyperref}
\allowdisplaybreaks[4]
\numberwithin{equation}{section}
\theoremstyle{plain} 
\newtheorem{theorem}{Theorem}[section]
\newtheorem{lemma}[theorem]{Lemma}
\newtheorem{corollary}[theorem]{Corollary}

\theoremstyle{definition}

\newtheorem{assumption}[theorem]{Assumption}
\newtheorem{remark}[theorem]{Remark}
\begin{document}

\title{Bismut formula and applications for conditional McKean--Vlasov SDEs driven by fractional Brownian motion}

\author{
    \small Wenze Chen, Zhenzhong Zhang\\
    \footnotesize{School of Mathematics and Statistics, Donghua University, Shanghai, 201620, China}\\
}

\date{}
\maketitle

\begin{abstract}
We investigate the stability of Lions derivatives under transition-rate perturbations for conditional McKean--Vlasov SDEs driven by fractional Brownian motion with $H>1/2$. For bounded measurable terminal functions, we establish Lions differentiability, a Bismut formula, and a gradient estimate under bounded first-order drift derivatives with suitable Hölder regularity, without imposing additional higher-order regularity assumptions on the drift's Lions derivative. Fractional Volterra estimates accommodate switching-induced discontinuities, while the Bismut weight has a common representation in terms of the switching process for all admissible transition-rate matrices. Combining this representation with likelihood-ratio methods, we prove quantitative stability, Fr\'echet differentiability, and a quadratic remainder estimate for Lions derivatives under perturbations preserving the positive off-diagonal support. In a nonlinear two-state model, the quadratic remainder estimate holds uniformly over terminal threshold events. We also demonstrate nonzero transition-rate sensitivity of Lions derivatives under fractional noise, whereas the corresponding Brownian model has a rate-independent terminal law.
\end{abstract}



\section{Introduction}

The stability and regularity of McKean--Vlasov dynamics with respect to their initial distributions have long been important topics in mean-field stochastic analysis. Since the evolving distribution enters the coefficients, perturbations of the initial law may influence the entire subsequent dynamics. Understanding this dependence is fundamental to the analysis of mean-field stochastic flows and their Lions derivatives \cite{MR3630288,Cardaliaguet2013}, with further applications to distributionally robust optimization \cite{MR5014209}. When the system evolves in a common Markovian environment, the transition-rate matrix introduces an additional source of uncertainty. This naturally raises the question of whether the Lions derivatives of terminal expectations remain stable under perturbations of the transition rates. Motivated by this problem, we investigate conditional McKean--Vlasov SDEs driven by fractional Brownian motion with $H>1/2$ and common Markovian switching.

An important aspect of this problem is the regularity of terminal expectations when the test functions are merely bounded and measurable. For smooth test functions, differentiability with respect to the initial distribution can be investigated through variational equations for McKean--Vlasov dynamics \cite{MR3630288}. However, this approach does not directly apply to discontinuous observables, such as indicators of terminal events. Bismut--Elworthy--Li formulas offer an alternative by representing derivatives of expectations through stochastic weights, without requiring differentiability of the test functions \cite{MR755001,MR1297021}. This approach has been extended to Lions derivatives of distribution-dependent SDEs, including equations driven by fractional Brownian motion \cite{MR3983053,MR4441502}. These developments provide a natural starting point for studying Lions regularity in conditional switching systems and, subsequently, its stability under transition-rate perturbations.

McKean--Vlasov SDEs, originating from the work of Kac \cite{MR84985} and McKean \cite{MR221595,MR233437}, describe the limiting dynamics of weakly interacting particle systems. Their probabilistic foundations were systematically developed by Sznitman \cite{MR1108185}. When particles share a common random environment, the limiting distribution is generally random and can be characterized through a conditional law. Early developments in this direction include the stochastic mean-field limits studied by Dawson and Vaillancourt \cite{MR1328577} and the particle representations of nonlinear stochastic partial differential equations developed by Kurtz and Xiong \cite{MR1705602}. Subsequent studies established important results on conditional propagation of chaos, well-posedness, and related probabilistic representations in common-noise systems \cite{MR3513598,MR4255126,MR4583663,MR4612111}. When the common environment is governed by a Markov chain, conditional mean-field limits and well-posedness have been further investigated in \cite{MR4035030,MR4518526,MR4690564}. Related regularity results for classical SDEs with Markovian switching, including smoothness of densities and a Bismut-type formula, were obtained by Hu et al.\ \cite{MR3986248}.

Fractional Brownian motion introduces additional analytical difficulties due to its temporal dependence, which distinguishes it from standard Brownian motion. Following the classical construction of Mandelbrot and Van Ness \cite{MR242239}, stochastic calculus and well-posedness theory for fractional differential equations were developed in \cite{MR1677455,MR1893308}. Within the mean-field framework, Buckdahn and Jing \cite{MR3648964} investigated fractional mean-field SDEs and associated stochastic control problems, while Galeati et al.\ \cite{MR4528970} established well-posedness results for distribution-dependent equations with irregular drifts. More recently, Shen and Wang \cite{ShenWang2025} studied well-posedness, conditional propagation of chaos, and averaging for conditional McKean--Vlasov SDEs driven by fractional Brownian motion with Markovian switching. Related conditional fractional equations with historical path dependence were further investigated in \cite{WangShenZhang2026}. These developments provide a foundation for studying conditional mean-field dynamics under fractional noise and motivate further investigation of their regularity properties with respect to initial distributions.

The development of Bismut formulas for distribution-dependent SDEs has provided important tools for such regularity problems. Ba\~nos \cite{MR3765887} investigated derivative representations for mean-field SDEs, while Crisan and McMurray \cite{MR3800831} established Wiener-space integration-by-parts formulas involving both state and Lions derivatives. 
Ren and Wang \cite{MR3983053} subsequently developed
Bismut formulas for Lions derivatives, and Bao et al.\
\cite{MR4219322} extended this approach to
distribution-path-dependent equations.
Using a different approach, Huang and Wang
\cite{MR4225923} obtained derivative estimates
for the laws of McKean--Vlasov SDEs through
heat-kernel parameter expansions.
In the fractional setting, Fan and Ren
\cite{MR3646611} established Bismut formulas
for SDEs driven by fractional Brownian motion,
while Fan et al.\ \cite{MR4441502} obtained
corresponding formulas for Lions derivatives
of distribution-dependent equations.
Further developments allowing certain distribution-dependent fractional noise coefficients were obtained in \cite{MR4890845}. For McKean--Vlasov SDEs with common Brownian noise, Bao et al.\ \cite{MR4975059} established Malliavin differentiability and an integration-by-parts formula. Related probabilistic representations and derivative estimates
for stochastic equations driven by stable processes were
developed by Xie et al.\ \cite{MR3264434}.

Alongside these regularity results, the stability of stochastic systems under perturbations of transition-rate matrices has attracted considerable attention. For continuous-time Markov chains, Glasserman \cite{MR1162948} investigated parameter-derivative estimators, including likelihood-ratio methods, while Nakayama et al.\ \cite{Nakayama1994} developed likelihood-ratio sensitivity analysis for Markovian reliability models. For stochastic differential equations with Markovian switching, Fan and Feng \cite{MR4485992} established stability estimates under perturbations of transition-rate matrices in the fractional Brownian setting. Shao \cite{MR4740950} further studied comparison principles and stability for regime-switching processes.

Existing Bismut formulas and likelihood-ratio methods address different types of perturbations, but their combination in the present setting is not straightforward. For bounded measurable terminal functions, direct differentiation of the test function is unavailable. Moreover, the regime-dependent diffusion coefficient changes discontinuously at switching times, making it necessary to establish suitable estimates for the inverse Volterra operator in constructing an admissible Bismut weight. Perturbations of the transition-rate matrix introduce a further difficulty by changing the probability law of the common switching process. We overcome these difficulties by establishing a Bismut representation under bounded first-order spatial and Lions derivatives of the drift with suitable Hölder regularity, without imposing the additional differentiability conditions on \(D^L b\) required in \cite[Theorem 4.7(I)]{MR4441502} for \(H>1/2\). We then construct the solution and its associated Bismut weight through a common representation in terms of the switching process, valid for all admissible transition-rate matrices. This yields quantitative stability, Fr\'echet differentiability, and a quadratic remainder estimate for Lions derivatives under perturbations preserving the positive off-diagonal support. In a nonlinear two-state model, we further show that fractional Brownian noise can produce nonzero transition-rate dependence of Lions derivatives, whereas the terminal law of the corresponding Brownian model is independent of the switching rate.

The remainder of this paper is organized as follows.
Section~\ref{sec:MF} introduces the model and the necessary preliminaries.
Section~\ref{sec:wellposedness} establishes the well-posedness of the conditional McKean--Vlasov equation.
Section~\ref{sec:bismut} develops the Lions differentiability of terminal expectations and the corresponding Bismut formula.
Section~\ref{sec:applications} investigates transition-rate stability and differentiability, followed by the analysis of a nonlinear two-state model and numerical illustrations.
Auxiliary estimates and technical proofs are collected in the appendices.

\section{Model Formulation}\label{sec:MF}

Throughout this paper, let
$H\in(1/2,1)$ and 
$\mathbb{S}=\{1,2,\ldots,m_0\}$ be a finite state space describing the
regimes of the common random environment. For $x,y\in\mathbb{R}^d$,
$\langle x,y\rangle$ and $|x|$ denote the Euclidean inner product and
norm, respectively. For a matrix $A$, $\|A\|$ denotes its Frobenius
norm. For $p\ge2$, we denote by $\mathcal{P}_p(\mathbb{R}^d)$ the space of all
probability measures on $\mathbb{R}^d$ with finite $p$th moment,
equipped with the $L^p$-Wasserstein distance
\begin{equation*}
    W_p(\mu,\nu)
    :=\left(
        \inf_{\pi\in\Pi(\mu,\nu)}
        \int_{\mathbb{R}^d\times\mathbb{R}^d}
        |x-y|^p\pi(\mathrm{d}x,\mathrm{d}y)
    \right)^{1/p},
\end{equation*}
where $\Pi(\mu,\nu)$ denotes the set of all couplings of $\mu$ and
$\nu$. When $p=2$, we simply write $W_2$ for the quadratic
Wasserstein distance.

We work on a product probability space to separate the common
switching environment from the fractional noise. Let
$(\Omega^0,\mathscr{F}^0,\{\mathscr{F}_t^0\}_{t\ge0},\mathbb{P}^0)$ and
$(\Omega^1,\mathscr{F}^1,\{\mathscr{F}_t^1\}_{t\ge0},\mathbb{P}^1)$
be two complete filtered probability spaces satisfying the usual
conditions. Set $\Omega=\Omega^0\times\Omega^1$, $\mathscr{F}=\mathscr{F}^0\otimes\mathscr{F}^1$, $\mathbb{P}=\mathbb{P}^0\otimes\mathbb{P}^1$,
and let $\{\mathscr{F}_t\}_{t\ge0}$ be the usual augmentation of
$\{\mathscr{F}_t^0\otimes\mathscr{F}_t^1\}_{t\ge0}$. We write
$\mathbb{E}$ for the full expectation under $\mathbb{P}$, and
$\mathbb{E}^0$ and $\mathbb{E}^1$ for expectations with respect to
$\mathbb{P}^0$ and $\mathbb{P}^1$, respectively. Thus, for any
integrable random variable $Z$ on $\Omega$,
\begin{equation*}
    \mathbb{E}[Z]
    =\mathbb{E}^0\big[\mathbb{E}^1[Z]\big]
    =\mathbb{E}^1\big[\mathbb{E}^0[Z]\big].
\end{equation*}

For $p\ge2$, let
$L^p(\Omega^1,\mathscr{F}_0^1,\mathbb{P}^1;\mathbb{R}^d)$
denote the space of all $\mathscr{F}_0^1$-measurable
$\mathbb{R}^d$-valued random variables $\xi$ such that
$\mathbb{E}^1|\xi|^p<\infty$. For $T>0$ and a Banach space $E$,
$L^p([0,T];E)$ denotes the Bochner space of strongly measurable
functions $f:[0,T]\to E$ satisfying
\begin{equation*}
    \|f\|_{L^p([0,T];E)}^p
    :=\int_0^T\|f_s\|_E^p \mathrm{d}s<\infty.
\end{equation*}
In particular, $L^p([0,T];\mathbb{R}^{d\times d})$ is understood
with respect to the Frobenius norm on $\mathbb{R}^{d\times d}$.
Let $C([0,T];\mathbb{R}^d)$ denote the Banach space of continuous
functions $f:[0,T]\to\mathbb{R}^d$, equipped with the supremum norm
\(
    \|f\|_{C([0,T];\mathbb{R}^d)}
    :=\sup_{t\in[0,T]}|f(t)|
\).
We write $L^p(\Omega;C([0,T];\mathbb{R}^d))$ for the Bochner space
of $C([0,T];\mathbb{R}^d)$-valued random variables $X$ satisfying
\begin{equation*}
    \|X\|_{L^p(\Omega;C([0,T];\mathbb{R}^d))}^p
    :=\mathbb{E}\sup_{t\in[0,T]}|X_t|^p<\infty.
\end{equation*}
Similarly, we have
\(
    \|X\|_{L^p(\Omega^1;C([0,T];\mathbb{R}^d))}^p
    :=\mathbb{E}^1\sup_{t\in[0,T]}|X_t|^p
\).
For processes depending on the common environment, we use
$L^p_{\mathscr{F}^0}(\Omega^0\times[0,T];\mathbb{R}^{d\times d})$
to denote the space of all $\{\mathscr{F}_t^0\}_{t\ge0}$-progressively
measurable processes
$u:[0,T]\times\Omega^0\to\mathbb{R}^{d\times d}$ such that
\begin{equation*}
    \|u\|_{L^p_{\mathscr{F}^0}(\Omega^0\times[0,T];\mathbb{R}^{d\times d})}^p
    :=\mathbb{E}^0\int_0^T\|u_s\|^p \mathrm{d}s<\infty.
\end{equation*}

On $(\Omega^0,\mathscr{F}^0,\mathbb{P}^0)$, let $(r_t)_{t\ge0}$
be a continuous-time Markov chain taking values in $\mathbb{S}$,
with conservative transition rate matrix
$Q=(q_{ij})_{m_0\times m_0}$. We take the environmental filtration
to be the usual augmentation of the natural filtration generated by
$r$, namely
\begin{equation*}
    \mathscr{F}_t^0
    :=\sigma(r_s:0\le s\le t)\vee\mathcal{N}^0,
    \quad\text{and}\quad
    \mathscr{F}^0
    :=\sigma(r_s:s\ge0)\vee\mathcal{N}^0,
\end{equation*}
where $\mathcal{N}^0$ denotes the collection of all
$\mathbb{P}^0$-null sets.

On $(\Omega^1,\mathscr{F}^1,\mathbb{P}^1)$, let
$B^H=(B^{H,1},\ldots,B^{H,d})$ be a $d$-dimensional fractional
Brownian motion with Hurst parameter $H\in(1/2,1)$, independent
of $r$. Its covariance function is given by
\begin{equation*}
    \mathbb{E}^1\big[B_t^{H,k}B_s^{H,\ell}\big]
    =\frac12\big(t^{2H}+s^{2H}-|t-s|^{2H}\big)\delta_{k\ell},
    \quad s,t\ge0,\ 1\le k,\ell\le d,
\end{equation*}
where $\delta_{k\ell}$ denotes the Kronecker delta.
We assume that the initial sigma-field $\mathscr F_0^1$ is
independent of $B^H$ and sufficiently rich in the sense that,
for every $p\ge2$ and $\mu\in\mathcal P_p(\mathbb R^d)$,
there exists
$\xi\in L^p(\Omega^1,\mathscr F_0^1,\mathbb P^1;\mathbb R^d)$
with $\mathcal L(\xi)=\mu$. By the product structure, $\xi$ is independent of
both $B^H$ and $r$. We take the spatial filtration to be the usual
augmentation of
\begin{equation*}
    \mathscr{F}_t^1
    :=\mathscr{F}_0^1\vee\sigma(B_s^H:0\le s\le t),
    \quad\text{and}\quad
    \mathscr{F}^1
    :=\mathscr{F}_0^1\vee\sigma(B_s^H:s\ge0)\vee\mathcal{N}^1,
\end{equation*}
where $\mathcal{N}^1$ denotes the collection of all
$\mathbb{P}^1$-null sets.

For later use, we recall some notation from the Malliavin calculus
for fractional Brownian motion; see Section~2.2 of
\cite{MR4441502}. Let $\mathcal{H}$ be the reproducing kernel
Hilbert space associated with $B^H$. We write
$L^2(\Omega;\mathcal{H})$ for the space of $\mathcal{H}$-valued
square-integrable random variables, equipped with the norm
\(
    \|h\|_{L^2(\Omega;\mathcal{H})}^2
    :=\mathbb{E}\|h\|_{\mathcal{H}}^2
\).
The covariance function admits the representation
\begin{equation*}
    R_H(t,s)
    =\int_0^{t\wedge s}K_H(t,r)K_H(s,r) \mathrm{d}r,
    \quad s,t\in[0,T],
\end{equation*}
where $K_H$ is the corresponding Volterra kernel. There exists a
$d$-dimensional Wiener process $W$ such that
\begin{equation*}
    B_t^H=\int_0^t K_H(t,s) \mathrm{d}W_s,
    \quad t\in[0,T].
\end{equation*}
Let $K_H^*:\mathcal{H}\to L^2([0,T];\mathbb{R}^d)$ be the canonical
isometry associated with $K_H$, and define
\(  R_H:=K_H\circ K_H^* \).
We denote by $K_H^{-1}$ the inverse of $K_H$ on its range.

We use Malliavin calculus with respect to $B^H$, keeping the
common environment and the initial $\sigma$-field fixed.
Let $\mathbb D$ and $\delta$ denote the Malliavin derivative
and its adjoint divergence operator, with domain
$\operatorname{Dom}\delta$, respectively.
In particular, $\mathbb D(ZF)=Z\mathbb DF$ whenever $Z$ is
bounded and $\mathscr F^0\otimes\mathscr F_0^1$-measurable,
and $F$ is a smooth cylindrical functional of $B^H$.
For vector-valued random variables, $\mathbb D$ acts
componentwise.
The corresponding operators for the underlying Wiener process
$W$ are denoted by $\mathbb D^W$ and $\delta^W$.
We use the transfer relations $K_H^*\mathbb DF=\mathbb D^WF$ and 
$\delta(u)=\delta^W(K_H^*u)$,
whenever the terms are well defined.
The Malliavin Sobolev space
$\mathbb D^{1,2}(\mathbb R^d)$ is equipped with the norm
\[
    \|F\|_{1,2}^2
    :=
    \mathbb E|F|^2
    +
    \mathbb E\|\mathbb DF\|_{\mathcal H}^2;
\]
see \cite{MR2200233}.
Throughout, integration with respect to $B^H$ is understood
in the divergence sense and denoted by
$\int_0^t\cdot \delta B_s^H$.

For an adapted process $X$ on the product probability space, its
conditional law with respect to the common switching environment
is denoted by $\mu_t:=\mathcal{L}(X_t\mid\mathscr{F}_t^0)$. By the
product structure and adaptedness, for $\mathbb{P}^0$-a.e.\ fixed
$\omega^0$, this conditional law is realized as
\(
    \mu_t(\omega^0)(\mathrm{d}x)
    =\mathbb{P}^1\left(X_t(\omega^0,\cdot)\in\mathrm{d}x\right)
\).
Thus $\mu_t$ is an $\mathscr{F}_t^0$-measurable random probability
measure. 

In this paper, we consider the following conditional 
McKean--Vlasov SDE driven
by fractional Brownian motion:
\begin{equation}\label{eq:main_MKVSDE}
    \begin{cases}
        \mathrm{d}X_t
        =b(t,X_t,\mu_t,r_t) \mathrm{d}t
        +\sigma(t,\mu_t,r_t) \delta B_t^H,\\[0.3em]
        X_0=\xi,\quad r_0=i_0\in\mathbb{S},
    \end{cases}
\end{equation}
where $\mu_t=\mathcal{L}(X_t\mid\mathscr{F}_t^0)$ and
$\xi\in L^p(\Omega^1,\mathscr{F}_0^1,\mathbb{P}^1;\mathbb{R}^d)$ with $\mathcal L(\xi)=\mu$.
The coefficients are measurable maps, namely
\begin{equation*}
    b:\mathbb{R}_+\times\mathbb{R}^d
    \times\mathcal{P}_p(\mathbb{R}^d)\times\mathbb{S}
    \to\mathbb{R}^d,
    \quad\text{and}\quad
    \sigma:\mathbb{R}_+\times\mathcal{P}_p(\mathbb{R}^d)
    \times\mathbb{S}\to\mathbb{R}^{d\times d}.
\end{equation*}

For $T>0$, let $\mathbb H^{T,p}_{\mathscr F^0}$
be the space of all $\{\mathscr F_t^0\}_{t\ge0}$-adapted,
$W_p$-continuous random measure flows
$\nu=(\nu_t)_{t\in[0,T]}$ taking values in
$\mathcal P_p(\mathbb R^d)$, identified up to
indistinguishability, such that
\[
    \|\nu\|_{\mathbb H^{T,p}_{\mathscr F^0}}^p
    :=
    \mathbb E^0\sup_{t\in[0,T]}
    \int_{\mathbb R^d}|x|^p\nu_t(dx)
    <\infty.
\]
For $\mu,\nu\in\mathbb H^{T,p}_{\mathscr F^0}$, set
\(
    \mathcal D_{T,p}(\mu,\nu)
    :=
    \left(
        \mathbb E^0\sup_{t\in[0,T]}
        W_p^p(\mu_t,\nu_t)
    \right)^{1/p}
\).
Since $(\mathcal P_p(\mathbb R^d),W_p)$ is complete,
$(\mathbb H^{T,p}_{\mathscr F^0},\mathcal D_{T,p})$
is a complete metric space
(cf.~\cite[(2.2)--(2.3)]{MR4583663}).

\section{Well-posedness} \label{sec:wellposedness}
In this section, we fix $p\ge 2$. 
\begin{assumption}\label{ass:wellposedness}
The coefficients $b$ and $\sigma$ are measurable in all variables and satisfy the following conditions.
\begin{enumerate}[(i)]
    \item For each $i\in\mathbb{S}$ and each $T>0$, there exists a constant $L_T>0$ such that for all $t\in[0,T]$, $x,y\in\mathbb{R}^d$, and $\mu,\nu\in\mathcal{P}_p(\mathbb{R}^d)$,
    \[
        |b(t,x,\mu,i)-b(t,y,\nu,i)|
         \le L_T\big(|x-y|+W_p(\mu,\nu)\big),
    \]
    \[
    \|\sigma(t,\mu,i)-\sigma(t,\nu,i)\|
         \le L_TW_p(\mu,\nu).
    \]
    \item For each $T>0$, there exists a constant $K_T>0$ such that for all $t\in[0,T]$ and $i\in\mathbb{S}$,
    \begin{equation*}
        |b(t,0,\delta_0,i)|+\|\sigma(t,\delta_0,i)\|
        \le K_T.
    \end{equation*}
    \item The initial value $\xi\in L^p(\Omega^1,\mathscr{F}^1_0,\mathbb{P}^1;\mathbb{R}^d)$  is independent of $B^H$ and $r$.
\end{enumerate}
\end{assumption}

\begin{lemma}\label{lem:deterministic_skorohod_estimate}
Let $H\in(1/2,1)$ and $p>1/H$. Let $u\in L^p_{\mathscr F^0}(\Omega^0\times[0,T];\mathbb R^{d\times d})$.
For $t\in[0,T]$, define
\begin{equation*}
    I_t(u):=\int_0^t u_s \delta B_s^H.
\end{equation*}
Then $I_\cdot(u)$ admits a continuous version and there exists a constant $C_{p,H}>0$ such that
\begin{equation*}
    \mathbb E\sup_{t\in[0,T]}|I_t(u)|^p
    \le
    C_{p,H}T^{pH-1}
    \mathbb E^0\int_0^T\|u_s\|^p \mathrm ds.
\end{equation*}
\end{lemma}

\begin{proof}
By Fubini's theorem, $u(\cdot,\omega^0)$ belongs to
$L^p([0,T];\mathbb R^{d\times d})$
for $\mathbb P^0$-a.e.\ $\omega^0\in\Omega^0$.
Fix such an $\omega^0$. Since $u(\cdot,\omega^0)$ is deterministic with respect to the fractional Brownian motion on the spatial probability space, the divergence integral coincides with the fractional Wiener integral
\begin{equation*}
    \int_0^t u_s(\omega^0) \delta B_s^H
    =
    \int_0^t u_s(\omega^0) \mathrm dB_s^H.
\end{equation*}
We use the deterministic fractional Wiener integral estimate established in the proof of Theorem 3.1 in \cite{MR4441502}: for every $v\in L^p([0,T];\mathbb R^{d\times d})$ with $p>1/H$,
\begin{equation*}
    \mathbb E^1\sup_{t\in[0,T]}
    \left|
        \int_0^t v_s \mathrm dB_s^H
    \right|^p
    \le
    C_{p,H}T^{pH-1}
    \int_0^T\|v_s\|^p \mathrm ds.
\end{equation*}
Applying this estimate to $v_s=u_s(\omega^0)$, we have
\begin{equation*}
    \mathbb E^1\sup_{t\in[0,T]}
    \left|
        \int_0^t u_s(\omega^0) \delta B_s^H
    \right|^p
    \le
    C_{p,H}T^{pH-1}
    \int_0^T\|u_s(\omega^0)\|^p \mathrm ds.
\end{equation*}
Taking $\mathbb E^0$ on both sides gives
\begin{equation*}
    \mathbb E\sup_{t\in[0,T]}|I_t(u)|^p
    \le
    C_{p,H}T^{pH-1}
    \mathbb E^0\int_0^T\|u_s\|^p \mathrm ds.
\end{equation*}
The continuity of $I_\cdot(u)$ follows from the same argument.
\end{proof}

\subsection{A fixed-point argument}
For a given measure flow
$\nu=(\nu_t)_{t\in[0,T]}\in\mathbb H_{\mathscr F^0}^{T,p}$
and $X_0^{(\nu)}=\xi$, consider the frozen equation
\begin{equation}\label{eq:frozen_equation}
    \mathrm{d}X_t^{(\nu)}
        =b(t,X_t^{(\nu)},\nu_t,r_t) \mathrm{d}t
        +\sigma(t,\nu_t,r_t) \delta B_t^H.
\end{equation}

\begin{lemma}\label{lem:frozen_equation}
Let Assumption \ref{ass:wellposedness} hold. Then, for every $T>0$ and $\nu\in\mathbb{H}_{\mathscr{F}^0}^{T,p}$, the frozen equation \eqref{eq:frozen_equation} admits a unique solution $X^{(\nu)}$ satisfying $\mathbb E\sup_{t\in[0,T]}|X_t^{(\nu)}|^p<\infty$.
Moreover, there exists a constant $C_T>0$, independent of $\nu$, such that
\begin{equation*}
    \mathbb E\sup_{t\in[0,T]} \left|X_t^{(\nu)}\right|^p
    \le C_T\left(1+\mathbb{E}^1|\xi|^p+\|\nu\|_{\mathbb{H}_{\mathscr{F}^0}^{T,p}}^p\right).
\end{equation*}
\end{lemma}

For readability, the proof of Lemma \ref{lem:frozen_equation} is deferred to  Appendix \ref{Appendix_frozen equation}.
For each $\nu\in\mathbb H_{\mathscr F^0}^{T,p}$, define
\begin{equation*}
    \Phi(\nu)_t(\omega^0)
    :=
    \mathbb P^1\circ
    \left(X_t^{(\nu)}(\omega^0,\cdot)\right)^{-1},
    \quad  t\in[0,T].
\end{equation*}
By the product structure of the probability space,
$\Phi(\nu)_t$ is $\mathscr F_t^0$-measurable and is a version of
$\mathcal L(X_t^{(\nu)}\mid\mathscr F_t^0)$. 
Moreover, by Lemma
\ref{lem:frozen_equation}, we know
\[
    \mathbb E^0\sup_{t\in[0,T]}
    \int_{\mathbb R^d}|x|^p\Phi(\nu)_t(dx)
    \le
    \mathbb E\sup_{t\in[0,T]}|X_t^{(\nu)}|^p
    <\infty.
\]
The continuity of $X^{(\nu)}$ and the
conditional coupling inequality imply that
$\Phi(\nu)$ has continuous paths and is
therefore progressively measurable.
Hence $\Phi$ maps $\mathbb H_{\mathscr F^0}^{T,p}$ into itself.
\begin{theorem}\label{thm:fixed_point_wellposedness}
Let Assumption \ref{ass:wellposedness} hold. Then, for every $T>0$,
the map $\Phi$ admits a unique fixed point
$\mu\in\mathbb H_{\mathscr F^0}^{T,p}$. Moreover,
\eqref{eq:main_MKVSDE} admits a unique solution $X$ satisfying
\(
    \mu_t
    =
    \mathcal L(X_t\mid\mathscr F_t^0)
\)
for $t\in[0,T]$, and
\begin{equation*}
    \mathbb E\sup_{t\in[0,T]}|X_t|^p
    \le
    C_T\left(1+\mathbb E^1|\xi|^p\right),
\end{equation*}
where $C_T>0$ is independent of $\xi$.
\end{theorem}

\begin{proof}
Fix any
$\mu^{(0)}\in\mathbb H_{\mathscr F^0}^{T,p}$.
For $n\ge0$, define $X^{(n+1)}:=X^{(\mu^{(n)})}$ and
$\mu^{(n+1)}:=\Phi(\mu^{(n)})$.
Equivalently,
\begin{equation}\label{eq:measure_picard_iteration}
    X_t^{(n+1)}
    =
    \xi
    +
    \int_0^t
        b(s,X_s^{(n+1)},\mu_s^{(n)},r_s) \mathrm ds
    +
    \int_0^t
        \sigma(s,\mu_s^{(n)},r_s) \delta B_s^H,
\end{equation}
with
\(
    \mu_t^{(n+1)}
    =
    \mathcal L(X_t^{(n+1)}\mid\mathscr F_t^0)
\).
For $n\ge1$, subtracting the equations for $X^{(n+1)}$ and $X^{(n)}$, we have
\begin{align*}
    X_t^{(n+1)}-X_t^{(n)}
    &=
    \int_0^t
    \Big[
        b(s,X_s^{(n+1)},\mu_s^{(n)},r_s)
        -
        b(s,X_s^{(n)},\mu_s^{(n-1)},r_s)
    \Big] \mathrm ds \\
    &\quad
    +
    \int_0^t
    \Big[
        \sigma(s,\mu_s^{(n)},r_s)
        -
        \sigma(s,\mu_s^{(n-1)},r_s)
    \Big] \delta B_s^H.
\end{align*}
Arguing as in the proof of Lemma \ref{lem:frozen_equation}, it follows
from the Lipschitz condition, H\"older's inequality, and Lemma
\ref{lem:deterministic_skorohod_estimate} that
\begin{equation}
    \mathbb E\sup_{u\in[0,t]}|X_u^{(n+1)}-X_u^{(n)}|^p
     \le
    C_T\int_0^t
        \mathbb E\sup_{u\in[0,s]}
        |X_u^{(n+1)}-X_u^{(n)}|^p \mathrm ds  +
    C_T\int_0^t
        \mathbb E^0
        W_p^p(\mu_s^{(n)},\mu_s^{(n-1)}) \mathrm ds.
\label{eq:picard_difference_preliminary}
\end{equation}
For $n\ge2$ and $\mathbb P^0$-a.e. $\omega^0$, the pair
\(
    \left(
        X_s^{(n)}(\omega^0,\cdot),
        X_s^{(n-1)}(\omega^0,\cdot)
    \right)
\)
induces a coupling of $\mu_s^{(n)}(\omega^0)$ and
$\mu_s^{(n-1)}(\omega^0)$. Therefore,
\begin{equation*}
    W_p^p(\mu_s^{(n)},\mu_s^{(n-1)})
    \le
    \mathbb E^1
    |X_s^{(n)}-X_s^{(n-1)}|^p,
    \quad 
    \mathbb P^0\text{-a.s.}
\end{equation*}
Set
\(
    a_n(t)
    :=
    \mathbb E\sup_{u\in[0,t]}
    |X_u^{(n)}-X_u^{(n-1)}|^p
\) for $n\ge2$.
Taking $\mathbb E^0$ in the preceding conditional coupling estimate gives
\begin{equation*}
    \mathbb E^0
    W_p^p(\mu_s^{(n)},\mu_s^{(n-1)})
    \le
    \mathbb E|X_s^{(n)}-X_s^{(n-1)}|^p
    \le
    a_n(s).
\end{equation*}
Substituting this into
\eqref{eq:picard_difference_preliminary} and applying Gronwall's
inequality, we obtain
\begin{equation}\label{eq:picard_recursive_estimate}
    a_{n+1}(t)
    \le
    C_T\int_0^t a_n(s) \mathrm ds,
    \quad  n\ge2.
\end{equation}
Iterating \eqref{eq:picard_recursive_estimate} yields
\(
    a_{n+1}(T)
    \le
    \frac{(C_TT)^{n-1}}{(n-1)!} a_2(T)
\) for $n \ge 2$.
Moreover, by Lemma \ref{lem:frozen_equation},
\[
    a_2(T)
    =
    \mathbb E\sup_{t\in[0,T]}|X_t^{(2)}-X_t^{(1)}|^p
    \le
    2^{p-1}
    \left(
        \mathbb E\sup_{t\in[0,T]}|X_t^{(2)}|^p
        +
        \mathbb E\sup_{t\in[0,T]}|X_t^{(1)}|^p
    \right)
    <\infty.
\]
Consequently,
\[
    \sum_{n=2}^{\infty}
    \left(
        \mathbb E\sup_{t\in[0,T]}
        |X_t^{(n+1)}-X_t^{(n)}|^p
    \right)^{1/p}
    =
    \sum_{n=2}^{\infty}a_{n+1}(T)^{1/p}
    \le
    a_2(T)^{1/p}
    \sum_{n=2}^{\infty}
    \left(
        \frac{(C_TT)^{n-1}}{(n-1)!}
    \right)^{1/p}
    <\infty.
\]
Hence $\{X^{(n)}\}_{n\ge1}$ is Cauchy with respect to the norm
\(
    \left(
        \mathbb E\sup_{t\in[0,T]}|X_t|^p
    \right)^{1/p}
\).
Thus there exists a continuous process $X$ such that
\begin{equation}\label{eq:Xn_convergence}
    \mathbb E\sup_{t\in[0,T]}|X_t^{(n)}-X_t|^p
    \longrightarrow0,
    \quad  n\to\infty.
\end{equation}
For $n,m\ge1$, the conditional coupling inequality gives
\[
    \mathcal D_{T,p}(\mu^{(n)},\mu^{(m)})^p
    =
    \mathbb E^0\sup_{t\in[0,T]}
    W_p^p(\mu_t^{(n)},\mu_t^{(m)})
    \le
    \mathbb E\sup_{t\in[0,T]}
    |X_t^{(n)}-X_t^{(m)}|^p
    \longrightarrow0.
\]
Thus $\{\mu^{(n)}\}_{n\ge1}$ is Cauchy in the complete metric space
$(\mathbb H_{\mathscr F^0}^{T,p},\mathcal D_{T,p})$.
Denote its limit by $\mu$.
By the Lipschitz continuity of $b$ and H\"older's inequality,
\begin{align*}
    &\quad
    \mathbb E\sup_{t\in[0,T]}
    \left|
        \int_0^t
        \Big[
            b(s,X_s^{(n+1)},\mu_s^{(n)},r_s)
            -
            b(s,X_s,\mu_s,r_s)
        \Big] \mathrm ds
    \right|^p \\
    &\le
    C_T\int_0^T
        \mathbb E\sup_{u\in[0,s]}
        |X_u^{(n+1)}-X_u|^p \mathrm ds
    +
    C_T\int_0^T
        \mathbb E^0
        W_p^p(\mu_s^{(n)},\mu_s) \mathrm ds,
\end{align*}
which converges to zero as $n\to\infty$. Similarly, by Lemma
\ref{lem:deterministic_skorohod_estimate},
\begin{align*}
    &\quad
    \mathbb E\sup_{t\in[0,T]}
    \left|
        \int_0^t
        \Big[
            \sigma(s,\mu_s^{(n)},r_s)
            -
            \sigma(s,\mu_s,r_s)
        \Big] \delta B_s^H
    \right|^p \\
    &\le
    C_{p,H}T^{pH-1}
    \int_0^T
    \mathbb E^0
    \left\|
        \sigma(s,\mu_s^{(n)},r_s)
        -
        \sigma(s,\mu_s,r_s)
    \right\|^p \mathrm ds \\
    & \le
    C_T
    \int_0^T
        \mathbb E^0
        W_p^p(\mu_s^{(n)},\mu_s) \mathrm ds
    \longrightarrow0,
    \quad  n\to\infty.
\end{align*}
Consequently, letting $n\to\infty$ in
\eqref{eq:measure_picard_iteration}, we obtain
\begin{equation}\label{eq:limit_frozen_equation}
    X_t
    =
    \xi
    +
    \int_0^t
        b(s,X_s,\mu_s,r_s) \mathrm ds
    +
    \int_0^t
        \sigma(s,\mu_s,r_s) \delta B_s^H,
    \quad  t\in[0,T].
\end{equation}

It remains to identify $\mu$ as the conditional law flow of $X$.
Let
\(
    \bar\mu_t
    :=
    \mathcal L(X_t\mid\mathscr F_t^0)
\).
By the continuity and moment bound of $X$,
we have $\bar\mu\in\mathbb H^{T,p}_{\mathscr F^0}$.
By the conditional coupling inequality and \eqref{eq:Xn_convergence}, we get
\[
    \mathcal D_{T,p}(\mu^{(n)},\bar\mu)^p
    \le
    \mathbb E\sup_{t\in[0,T]}|X_t^{(n)}-X_t|^p
    \longrightarrow0.
\]
Consequently, as $n \to \infty$,
\[
    \mathcal D_{T,p}(\mu,\bar\mu)
    \le
    \mathcal D_{T,p}(\mu,\mu^{(n)})
    +
    \mathcal D_{T,p}(\mu^{(n)},\bar\mu)
    \longrightarrow0.
\]
Therefore,
\begin{equation}\label{eq:conditional_law_identification}
    \mu_t
    =
    \bar\mu_t
    =
    \mathcal L(X_t\mid\mathscr F_t^0),
    \quad  t\in[0,T].
\end{equation}
By \eqref{eq:limit_frozen_equation}, $X$ solves the frozen equation
associated with $\mu$. Hence, by the uniqueness in Lemma
\ref{lem:frozen_equation} and
\eqref{eq:conditional_law_identification}, we know
\(
    \Phi(\mu)_t
    =
    \mathcal L(X_t\mid\mathscr F_t^0)
    =
    \mu_t
\) for $ t\in[0,T]$.
Therefore, $\mu$ is a fixed point of $\Phi$, and $X$ is a solution of
\eqref{eq:main_MKVSDE}.

We next prove uniqueness. Let $X$ and $Y$ be two solutions of
\eqref{eq:main_MKVSDE} on the same product probability space, driven by
the same $B^H$ and $r$, with the same initial value $\xi$. Set $\mu_t=\mathcal L(X_t\mid\mathscr F_t^0)$ and $\nu_t=\mathcal L(Y_t\mid\mathscr F_t^0)$.
Using the same difference estimate and conditional coupling argument as
above, we obtain
\begin{equation*}
    \mathbb E\sup_{u\in[0,t]}|X_u-Y_u|^p
    \le
    C_T\int_0^t
        \mathbb E\sup_{u\in[0,s]}|X_u-Y_u|^p \mathrm ds,
\end{equation*}
where we have used
\(
    \mathbb E^0
    W_p^p(\mu_s,\nu_s)
    \le
    \mathbb E|X_s-Y_s|^p
\).
It follows from Gronwall's inequality that
\(
    \mathbb E\sup_{t\in[0,T]}|X_t-Y_t|^p
    =0
\).
Hence $X=Y$ almost surely on $[0,T]$.
Finally, if $\rho\in\mathbb H_{\mathscr F^0}^{T,p}$ is another fixed
point of $\Phi$, then
\(
    \rho_t
    =
    \mathcal L(X_t^{(\rho)}\mid\mathscr F_t^0)
\).
Hence $X^{(\rho)}$ is also a solution of \eqref{eq:main_MKVSDE}. By the
uniqueness proved above, $X^{(\rho)}=X$, and therefore $\rho=\mu$. Thus the
fixed point of $\Phi$ is unique.

It remains to establish the moment estimate. Now that the identity
$\mu_t=\mathcal L(X_t\mid\mathscr F_t^0)$ has been established,
applying the growth and integral estimates used in the proof of
Lemma \ref{lem:frozen_equation} to the nonlinear solution $X$ yields, for
$t\in[0,T]$,
\[
    \mathbb E\sup_{s\in[0,t]}|X_s|^p
    \le
    C_T\left(1+\mathbb E^1|\xi|^p\right)
    +
    C_T\int_0^t
        \mathbb E\sup_{v\in[0,u]}|X_v|^p \mathrm du 
    +
    C_T\int_0^t
        \mathbb E^0
        \left[
            \int_{\mathbb R^d}|y|^p\mu_u(\mathrm dy)
        \right]\mathrm du.
\]
Since
\(
    \mu_u=\mathcal L(X_u\mid\mathscr F_u^0)
\),
we have
\begin{equation*}
    \mathbb E^0
    \left[
        \int_{\mathbb R^d}|y|^p\mu_u(\mathrm dy)
    \right]
    =
    \mathbb E|X_u|^p
    \le
    \mathbb E\sup_{v\in[0,u]}|X_v|^p.
\end{equation*}
Therefore,
\begin{equation*}
    \mathbb E\sup_{s\in[0,t]}|X_s|^p
    \le
    C_T\left(1+\mathbb E^1|\xi|^p\right)
    +
    C_T\int_0^t
        \mathbb E\sup_{v\in[0,u]}|X_v|^p \mathrm du.
\end{equation*}
It follows from Gronwall's inequality that
\begin{equation*}
    \mathbb E\sup_{t\in[0,T]}|X_t|^p
    \le
    C_T\left(1+\mathbb E^1|\xi|^p\right).
\end{equation*}
The proof is complete.
\end{proof}

\begin{remark}
Shen and Wang~\cite[Theorem~3.2]{ShenWang2025}
established well-posedness for conditional
McKean--Vlasov SDEs, where their proof employs
a direct Picard iteration~\cite[(3.1)]{ShenWang2025},
with convergence established through an
Osgood--Bihari argument under non-Lipschitz conditions.
\end{remark}

\section{Bismut Formula}
\label{sec:bismut}

For the Bismut analysis, we take \(p=2\), and $D^L$ denotes the Lions derivative with
respect to the probability measure. 

We recall the related notions of
directional, weak and full $L$-differentiability.
Let $F:\mathcal P_2(\mathbb R^d)\to\mathbb R$ and
$\mu\in\mathcal P_2(\mathbb R^d)$. We write
$L^2(\mathbb R^d \to \mathbb R^d , \mu)$ for the space of all
$\mathbb R^d$-valued measurable functions $\phi$ such that
\[
    \|\phi\|_{L^2_\mu} 
    :=
    \left( \int_{\mathbb R^d}
    |\phi(x)|^2\mu(\mathrm dx) \right) ^{\frac{1}{2}}
    <\infty.
\]
For
$\phi\in  L^2(\mathbb R^d \to \mathbb R^d , \mu)$, the directional
$L$-derivative of $F$ at $\mu$ along $\phi$ is defined, whenever the
limit exists, by
\[
    D_\phi^L F(\mu)
    :=
    \lim_{\varepsilon\downarrow0}
    \frac{
        F\left(
            \mu\circ(\operatorname{Id}+\varepsilon\phi)^{-1}
        \right)
        -
        F(\mu)
    }{\varepsilon}.
\]
If
\(
    \phi\mapsto D_\phi^L F(\mu)
\)
is a bounded linear functional on
$L^2(\mathbb R^d \to \mathbb R^d , \mu)$,
then $F$ is said to be weakly $L$-differentiable
at $\mu$, also called intrinsically differentiable
at $\mu$ (see \cite{MR4219322,MR4225923}).
Moreover, $F$ is said to be $L$-differentiable at $\mu$ if it is
weakly $L$-differentiable at $\mu$ and
\[
    \lim_{\|\phi\|_{L^2_\mu} \to0}
    \frac{
        \left|
            F\left(
                \mu\circ(\operatorname{Id}+\phi)^{-1}
            \right)
            -
            F(\mu)
            -
            D_\phi^L F(\mu)
        \right|
    }{
        \|\phi\|_{L^2_\mu} 
    }
    =0.
\]
When the Lions derivative $D^L F(\mu)(\cdot)$ exists,
\[
    D_\phi^L F(\mu)
    =
    \int_{\mathbb R^d}
    \left\langle
        D^L F(\mu)(x),
        \phi(x)
    \right\rangle
    \mu(\mathrm dx).
\]

\begin{remark}
\label{rem:distribution-free-noise}

For the Bismut formula established below, we restrict the
diffusion coefficient to the form $\sigma(t,r_t)$. In \cite{MR4441502}, the well-posedness result allows
distribution-dependent diffusion coefficients, whereas the Bismut formula
for the Lions derivative is established under distribution-free noise.
The difficulty is that the Lions variation then contains an additional
stochastic term involving $D^L\sigma$, which has no corresponding term
in the Malliavin variation.

More recently, \cite{MR4890845} treated distribution-dependent fractional
noise by introducing an additional independent fractional Brownian motion
with a distribution-free and invertible coefficient. The authors also point
out that their method is not applicable when this invertible coefficient
itself depends on the distribution.
\end{remark}

Hence, we consider the following  conditional 
McKean--Vlasov SDE:
\begin{equation}
\label{eq:regularity-model}
    \mathrm dX_t
    =
    b(t,X_t,\mu_t,r_t)\mathrm dt
    +
    \sigma(t,r_t)\delta B_t^H,
    \quad
    H\in(1/2,1),
\end{equation} 
For each $\mu\in\mathcal P_2(\mathbb R^d)$,
choose
$\xi\in L^2(\Omega^1,\mathscr F_0^1,\mathbb P^1;\mathbb R^d)$
with $\mathcal L(\xi)=\mu$.
We write $X$ for the corresponding solution of
\eqref{eq:regularity-model} with $X_0=\xi$ and
$\mu_t=\mathcal L(X_t\mid\mathscr F_t^0)$ and use the notation
\(
    C^{1,(1,0)}
    \left(
        \mathbb R^d\times\mathcal P_2(\mathbb R^d);
        \mathbb R^d
    \right)
\)
for the class of functions $f(x,\mu)$, which are continuously differentiable
with respect to $x$, Lions differentiable with respect to $\mu$, and such
that $(x,\mu)\mapsto\nabla_x f(x,\mu)$ and
$(x,\mu,z)\mapsto D^L f(x,\mu)(z)$ are continuous. See, e.g., \cite[Definition~2.1]{MR4441502}.

\begin{assumption}
\label{ass:regularity-b}

For every $T>0$, there exists a constant $K_T>0$ such that the following
conditions hold uniformly for $i\in\mathbb S$.

\begin{enumerate}[(i)]

    \item
    For every $t\in[0,T]$, it holds that
    \(
        b(t,\cdot,\cdot,i)
        \in
        C^{1,(1,0)}
        \left(
            \mathbb R^d\times\mathcal P_2(\mathbb R^d);
            \mathbb R^d
        \right)
    \)
    and
    \[
        \sup_{\substack{
            t\in[0,T],\ x,z\in\mathbb R^d,\\
            \mu\in\mathcal P_2(\mathbb R^d),\ i\in\mathbb S
        }}
        \left\{
            \|\nabla_x b(t,x,\mu,i)\|
            +
            \|D^L b(t,x,\mu,i)(z)\|
        \right\}
        \le K_T.
    \]

    \item
    There exist
    $\alpha_0\in\left(H-\frac12,1\right]$ and
    $\beta_0,\gamma_0\in\left(1-\frac{1}{2H},1\right]$
    such that, for all
    $s,t\in[0,T]$, $x,y,z_1,z_2\in\mathbb R^d$,
    $\mu,\nu\in\mathcal P_2(\mathbb R^d)$ and $i\in\mathbb S$,
    \begin{align*}
        &\quad
        \|\nabla_x b(t,x,\mu,i)
        -
        \nabla_x b(s,y,\nu,i)\|
        +
        \|D^L b(t,x,\mu,i)(z_1)
        -
        D^L b(s,y,\nu,i)(z_2)\|
        \nonumber\\
        &\le
        K_T
        \left(
            |t-s|^{\alpha_0}
            +
            |x-y|^{\beta_0}
            +
            |z_1-z_2|^{\gamma_0}
            +
            W_2(\mu,\nu)
        \right).
    \end{align*}

\end{enumerate}

\end{assumption}

\begin{assumption}
\label{ass:regularity-sigma}

For every $T>0$, there exist constants $K_T>0$ and
\(
    \zeta_0 \in\left(H-\frac12,1\right]
\)
such that, for every $i\in\mathbb S$ and $t\in[0,T]$,
$\sigma(t,i)$ is invertible and
\begin{equation*}
    \sup_{t\in[0,T],\ i\in\mathbb S}
    \left(
        \|\sigma(t,i)\|
        +
        \|\sigma^{-1}(t,i)\|
    \right)
    \le K_T.
\end{equation*}
Moreover, for all $s,t\in[0,T]$ and $i\in\mathbb S$,
\begin{equation}
\label{eq:sigma-inverse-holder}
    \|\sigma^{-1}(t,i)-\sigma^{-1}(s,i)\|
    \le
    K_T |t-s|^{\zeta_0}.
\end{equation}

\end{assumption}

For $f\in\mathscr B_b(\mathbb R^d)$, we define
\begin{equation}
\label{eq:PT-functional}
    (P_Tf)(\mu)
    :=
    \mathbb E\left[f(X_T)\right],
    \quad
    \mu\in\mathcal P_2(\mathbb R^d).
\end{equation}
For
$\phi\in  L^2(\mathbb R^d \to \mathbb R^d , \mu)$ and
$\varepsilon\in[0,1]$, set $\xi^{\varepsilon,\phi}
    :=
    \xi+\varepsilon\phi(\xi)$, $\mu^{\varepsilon,\phi}
    :=
    \mathcal L(\xi^{\varepsilon,\phi})$,
and let
$X^{\varepsilon,\phi}$
be the corresponding solution of
\eqref{eq:regularity-model} with initial value $\xi^{\varepsilon, \phi}$. We write
\(
    \mu_t^{\varepsilon,\phi}
    :=
    \mathcal L
    \left(
        X_t^{\varepsilon,\phi}
        \mid
        \mathscr F_t^0
    \right)
\).

Referring to \cite{MR3630288}, let
$(\widetilde\Omega^1,\widetilde{\mathscr F}^1,
\widetilde{\mathbb P}^1)$
be a copy of
$(\Omega^1,\mathscr F^1,\mathbb P^1)$.
We use a tilde to denote the corresponding copy and
$\widetilde{\mathbb E}^1$ for the expectation with respect to
$\widetilde{\mathbb P}^1$ only. In particular, the common environment and
the Markov chain $r$ are kept fixed under
$\widetilde{\mathbb E}^1$.
For later use, set
\begin{equation}
\label{eq:epsilon-variation}
    G_t^{\varepsilon,\phi}
    :=
    \nabla_{\phi(\xi)}
    X_t^{\varepsilon,\phi},
    \quad
    t\in[0,T],
\end{equation}
where the initial directional variation is understood in the sense of
Lemma~\ref{lem:initial-variation}, and
\begin{equation}
\label{eq:epsilon-law-variation}
    F_t^{\varepsilon,\phi}
    :=
    \widetilde{\mathbb E}^{1}
    \left[
        D^L b
        \left(
            t,
            X_t^{\varepsilon,\phi},
            \mu_t^{\varepsilon,\phi},
            r_t
        \right)
        \left(
            \widetilde X_t^{\varepsilon,\phi}
        \right)
        \widetilde G_t^{\varepsilon,\phi}
    \right].
\end{equation}
Then, by Lemma~\ref{lem:initial-variation}, we get
\begin{equation}
\label{eq:epsilon-variation-equation}
    \mathrm dG_t^{\varepsilon,\phi}
    =
    \left[
        \nabla_x b
        \left(
            t,
            X_t^{\varepsilon,\phi},
            \mu_t^{\varepsilon,\phi},
            r_t
        \right)
        G_t^{\varepsilon,\phi}
        +
        F_t^{\varepsilon,\phi}
    \right]
    \mathrm dt,
    \quad
    G_0^{\varepsilon,\phi}
    =
    \phi(\xi).
\end{equation}
For $s_0\in[0,T)$, define
$h_{s_0}^{\varepsilon,\phi}$ through
\begin{equation}
\label{eq:explicit-control-epsilon}
    (R_Hh_{s_0}^{\varepsilon,\phi})(t)
    :=
    \begin{cases}
        0,
        & 0\le t\le s_0,
        \\[1mm]
        \displaystyle
        \int_{s_0}^t
        \sigma^{-1}(s,r_s)
        \left[
            \frac{1}{T-s_0}
            G_s^{\varepsilon,\phi}
            +
            \frac{s-s_0}{T-s_0}
            F_s^{\varepsilon,\phi}
        \right]
        \mathrm ds,
        & s_0<t\le T.
    \end{cases}
\end{equation}
For the admissibility argument below, let
\[
    \mathscr G_t
    :=
    \sigma\left(
        \xi,W_s:0\le s\le t
    \right)
    \vee
    \mathscr F_t^0
    \vee
    \mathcal N,
    \quad
    t\in[0,T],
\]
where $\mathcal N$ denotes the collection of $\mathbb P$-null sets.

\begin{lemma}
\label{lem:switching-cameron-martin}

Suppose that Assumptions~\ref{ass:wellposedness}, \ref{ass:regularity-b} and \ref{ass:regularity-sigma} hold. Then, for every
$\varepsilon\in[0,1]$, $s_0\in[0,T)$ and
$\phi\in  L^2(\mathbb R^d \to \mathbb R^d , \mu)$,
the process $h_{s_0}^{\varepsilon,\phi}$ defined by
\eqref{eq:explicit-control-epsilon} belongs to
$\operatorname{Dom}\delta$.
More precisely,
\(
    K_H^{-1}
    (R_Hh_{s_0}^{\varepsilon,\phi})
    \in
    L^2([0,T]\times\Omega;\mathbb R^d)
\)
and is $\{\mathscr G_t\}_{t\in[0,T]}$-adapted.
Moreover,
\begin{equation}
\label{eq:control-integrability-epsilon}
    \int_0^1
    \left(
        \mathbb E
        \left|
            \delta(h_{s_0}^{\varepsilon,\phi})
        \right|^2
    \right)^{1/2}
    \mathrm d\varepsilon
    <\infty.
\end{equation}
In particular, there exists a constant
$C_{T,s_0,\mu}>0$ such that
\begin{equation}
\label{eq:control-estimate-zero}
    \mathbb E
    \left|
        \delta(h_{s_0}^{0,\phi})
    \right|^2
    \le
    C_{T,s_0,\mu}
    \|\phi\|_{L^2_\mu} ^2.
\end{equation}

\end{lemma}

\begin{proof}

By Assumptions~\ref{ass:regularity-b}(ii) and
\ref{ass:regularity-sigma}, we know
\(
    H-\frac12
    <
    \min
    \left\{
        \frac12,
        \zeta_0,
        \alpha_0,
        H\beta_0,
        H\gamma_0
    \right\}
\).
Fix
\begin{equation}
\label{eq:kappa-switching-control}
    \kappa
    \in
    \left(
        H-\frac12,
        \min
        \left\{
            \frac12,
            \zeta_0,
            \alpha_0,
            H\beta_0,
            H\gamma_0
        \right\}
    \right).
\end{equation}
We divide the proof into three steps.
 
\noindent \textit{Step 1: regularity of the variational terms.}
By the definition \eqref{eq:epsilon-law-variation} and
Assumption~\ref{ass:regularity-b}(i), we have
\begin{equation}
\label{eq:F-by-G-bound}
    \left|
        F_t^{\varepsilon,\phi}
    \right|
    \le
    K_T
    \widetilde{\mathbb E}^{1}
    \left|
        \widetilde G_t^{\varepsilon,\phi}
    \right|\le
    K_T
    \left(
        \widetilde{\mathbb E}^{1}
        \left|
            \widetilde G_t^{\varepsilon,\phi}
        \right|^2
    \right)^{1/2}.
\end{equation}
Since
$(\widetilde\Omega^1,\widetilde{\mathscr F}^1,
\widetilde{\mathbb P}^1)$
is a copy of
$(\Omega^1,\mathscr F^1,\mathbb P^1)$, we know
\[
    \widetilde{\mathbb E}^{1}
    \left|
        \widetilde G_t^{\varepsilon,\phi}
    \right|^2
    =
    \mathbb E^{1}
    \left|
        G_t^{\varepsilon,\phi}
    \right|^2,
    \quad
    \mathbb P^0\text{-a.s.}
\]
Hence, it follows from Minkowski's inequality, Assumption~\ref{ass:regularity-b}(i), and \eqref{eq:F-by-G-bound} that
\begin{align*}
    \left(
        \mathbb E^1
        \left|
            G_t^{\varepsilon,\phi}
        \right|^2
    \right)^{1/2}
    &\le
    \|\phi\|_{L^2_\mu} 
    +
    K_T
    \int_0^t
    \left(
        \mathbb E^1
        \left|
            G_s^{\varepsilon,\phi}
        \right|^2
    \right)^{1/2}
    \mathrm ds
    +
    \int_0^t
    \left(
        \mathbb E^1
        \left|
            F_s^{\varepsilon,\phi}
        \right|^2
    \right)^{1/2}
    \mathrm ds
    \\
    &\le
    \|\phi\|_{L^2_\mu} 
    +
    C_T
    \int_0^t
    \left(
        \mathbb E^1
        \left|
            G_s^{\varepsilon,\phi}
        \right|^2
    \right)^{1/2}
    \mathrm ds.
\end{align*}
Therefore, by Gronwall's inequality we obtain
\begin{equation}
\label{eq:G-epsilon-bound}
    \sup_{\varepsilon\in[0,1]}
    \sup_{t\in[0,T]}
    \mathbb E^1
    \left|
        G_t^{\varepsilon,\phi}
    \right|^2
    \le
    C_T
    \|\phi\|_{L^2_\mu} ^2,
    \quad
    \mathbb P^0\text{-a.s.}
\end{equation}
Consequently,
\begin{equation*}
    \sup_{\varepsilon\in[0,1]}
    \sup_{t\in[0,T]}
    \mathbb E
    \left|
        F_t^{\varepsilon,\phi}
    \right|^2
    \le
    C_T
    \|\phi\|_{L^2_\mu} ^2.
\end{equation*}
It also follows from
\eqref{eq:epsilon-variation-equation} that
\[
    G_t^{\varepsilon,\phi}
    -
    G_s^{\varepsilon,\phi}
    =
    \int_s^t
    \nabla_x b
    \left(
        u,
        X_u^{\varepsilon,\phi},
        \mu_u^{\varepsilon,\phi},
        r_u
    \right)
    G_u^{\varepsilon,\phi}
    \mathrm du
    +
    \int_s^t
    F_u^{\varepsilon,\phi}
    \mathrm du .
\]
Therefore, by the boundedness of $\nabla_xb$,
\eqref{eq:F-by-G-bound}, and
\eqref{eq:G-epsilon-bound}, we have
\begin{equation}
\label{eq:G-time-increment}
    \mathbb E^1
    \left|
        G_t^{\varepsilon,\phi}
        -
        G_s^{\varepsilon,\phi}
    \right|^2
    \le
    C_T
    \|\phi\|_{L^2_\mu} ^2
    |t-s|^2,
    \quad
    0\le s<t\le T,
\end{equation}
for $\mathbb P^0$-a.e.\ common-environment path.

We next estimate the time regularity of
$F^{\varepsilon,\phi}$. By its definition, we know
\begin{equation}\label{F-F}
    F_t^{\varepsilon,\phi}
    -
    F_s^{\varepsilon,\phi}
    =
    \widetilde{\mathbb E}^{1}
    \Big[
        A_{s,t}^{\varepsilon,\phi}
        \widetilde G_t^{\varepsilon,\phi}
    \Big]
    +
    \widetilde{\mathbb E}^{1}
    \Bigg[
        D^Lb
        \left(
            s,
            X_s^{\varepsilon,\phi},
            \mu_s^{\varepsilon,\phi},
            r_s
        \right)
        \left(
            \widetilde X_s^{\varepsilon,\phi}
        \right)
        \left(
            \widetilde G_t^{\varepsilon,\phi}
            -
            \widetilde G_s^{\varepsilon,\phi}
        \right)
    \Bigg],
\end{equation}
where
\[
    A_{s,t}^{\varepsilon,\phi}
    :=
    D^Lb
    \left(
        t,
        X_t^{\varepsilon,\phi},
        \mu_t^{\varepsilon,\phi},
        r_t
    \right)
    \left(
        \widetilde X_t^{\varepsilon,\phi}
    \right)
    -
    D^Lb
    \left(
        s,
        X_s^{\varepsilon,\phi},
        \mu_s^{\varepsilon,\phi},
        r_s
    \right)
    \left(
        \widetilde X_s^{\varepsilon,\phi}
    \right).
\]
Let
\(
    \lambda_Q
    :=
    \max_{i\in\mathbb S}(-q_{ii})
\).
Since $r$ is a finite-state continuous-time Markov chain,
\begin{equation}
\label{eq:switching-probability}
    \mathbb P^0
    \left(
        r_t\ne r_s
        \mid
        \mathscr F_s^0
    \right)
    \le
    1-e^{-\lambda_Q(t-s)}
    \le
    \lambda_Q(t-s).
\end{equation}

On $\{r_t=r_s\}$, by
Assumption~\ref{ass:regularity-b}(ii) we have
\[
    \left|
        A_{s,t}^{\varepsilon,\phi}
    \right|
    \le
    K_T
    \Big(
        |t-s|^{\alpha_0}
        +
        \left|
            X_t^{\varepsilon,\phi}
            -
            X_s^{\varepsilon,\phi}
        \right|^{\beta_0}
        +
        \left|
            \widetilde X_t^{\varepsilon,\phi}
            -
            \widetilde X_s^{\varepsilon,\phi}
        \right|^{\gamma_0}
        +
        W_2
        \left(
            \mu_t^{\varepsilon,\phi},
            \mu_s^{\varepsilon,\phi}
        \right)
    \Big),
\]
whereas on $\{r_t\ne r_s\}$, the boundedness of $D^Lb$
gives
\(
    \left|
        A_{s,t}^{\varepsilon,\phi}
    \right|
    \le
    2K_T
\).
For every $R>0$, by the moment estimate for the solution we know
\begin{equation}\label{eq:Xepsilon-moment-local}
    \sup_{
        \varepsilon\in[0,1],\
        \|\phi\|_{L^2_\mu} \le R
    }
    \mathbb E
    \sup_{t\in[0,T]}
    \left|
        X_t^{\varepsilon,\phi}
    \right|^2
    \le
    C_{T,\mu,R}.
\end{equation}
Indeed,
\[
    \mathbb E^1
    \left|
        \xi^{\varepsilon,\phi}
    \right|^2
    \le
    2\mathbb E^1|\xi|^2
    +
    2\varepsilon^2
    \mathbb E^1|\phi(\xi)|^2
    \le
    2\mathbb E^1|\xi|^2+2R^2,
\]
and the estimate follows from Theorem \ref{thm:fixed_point_wellposedness}.
Moreover, for $0\le s<t\le T$,
\[
    X_t^{\varepsilon,\phi}
    -
    X_s^{\varepsilon,\phi}
    =
    \int_s^t
    b
    \left(
        u,
        X_u^{\varepsilon,\phi},
        \mu_u^{\varepsilon,\phi},
        r_u
    \right)
    \mathrm du
    +
    \int_s^t
    \sigma(u,r_u)
    \delta B_u^H .
\]
By the Cauchy--Schwarz inequality, the linear growth
condition on $b$, and \eqref{eq:Xepsilon-moment-local},
we obtain
\[
    \mathbb E
    \left|
        \int_s^t
        b(u,X_u^{\varepsilon,\phi},
        \mu_u^{\varepsilon,\phi},r_u)
        \mathrm du
    \right|^2
    \le
    |t-s|
    \int_s^t
    \mathbb E
    \left|
        b(u,X_u^{\varepsilon,\phi},
        \mu_u^{\varepsilon,\phi},r_u)
    \right|^2
    \mathrm du\le
    C_{T,\mu,R}|t-s|^2.
\]
Moreover, conditioning on the common environment
and using the covariance formula for fractional
Wiener integrals, we obtain
\[
    \mathbb E
    \left|
        \int_s^t
        \sigma(u,r_u) \delta B_u^H
    \right|^2\le
    C_T H(2H-1)
    \int_s^t\int_s^t
    |u-v|^{2H-2}
    \mathrm du \mathrm dv=
    C_T|t-s|^{2H},
\]
where we used the independence of $r$ and $B^H$
and the uniform boundedness of $\sigma$.
Since $H<1$, the drift contribution can be absorbed into the same
bound on $[0,T]$. Thus
\begin{equation}
\label{eq:Xepsilon-time-increment}
    \sup_{
        \varepsilon\in[0,1],\
        \|\phi\|_{L^2_\mu} \le R
    }
    \mathbb E
    \left|
        X_t^{\varepsilon,\phi}
        -
        X_s^{\varepsilon,\phi}
    \right|^2
    \le
    C_{T,\mu,R}
    |t-s|^{2H}.
\end{equation}
For $\mathbb P^0$-a.e.\ $\omega^0$, the pair
\(
    \left(
        X_t^{\varepsilon,\phi}(\omega^0,\cdot),
        X_s^{\varepsilon,\phi}(\omega^0,\cdot)
    \right)
\)
is a coupling of
$\mu_t^{\varepsilon,\phi}(\omega^0)$ and
$\mu_s^{\varepsilon,\phi}(\omega^0)$. Hence
\[
    W_2^2
    \left(
        \mu_t^{\varepsilon,\phi},
        \mu_s^{\varepsilon,\phi}
    \right)
    \le
    \mathbb E^1
    \left|
        X_t^{\varepsilon,\phi}
        -
        X_s^{\varepsilon,\phi}
    \right|^2,
    \quad
    \mathbb P^0\text{-a.s.},
\]
and therefore
\begin{equation}
\label{eq:conditional-W2-time}
    \mathbb E^0
    W_2^2
    \left(
        \mu_t^{\varepsilon,\phi},
        \mu_s^{\varepsilon,\phi}
    \right)
    \le
    \mathbb E
    \left|
        X_t^{\varepsilon,\phi}
        -
        X_s^{\varepsilon,\phi}
    \right|^2.
\end{equation}

By the Cauchy--Schwarz inequality on the copied probability space
and \eqref{eq:G-epsilon-bound}, we get
\begin{align*}
    \mathbb E\left[
        \mathbf1_{\{r_t\ne r_s\}}
        \left|\widetilde{\mathbb E}^{1}
            [A_{s,t}^{\varepsilon,\phi}
             \widetilde G_t^{\varepsilon,\phi}]
        \right|^2
    \right] 
    & \le
    4K_T^2\mathbb E^0\left[
        \mathbf1_{\{r_t\ne r_s\}}
        \mathbb E^1|G_t^{\varepsilon,\phi}|^2
    \right]\\
    &\le
    C_T\|\phi\|_{L^2_\mu} ^2
    \mathbb P^0(r_t\ne r_s)
    \le C_T\|\phi\|_{L^2_\mu} ^2|t-s|.
\end{align*}
On $\{r_t=r_s\}$, it follows from \eqref{eq:Xepsilon-time-increment},
\eqref{eq:conditional-W2-time}, and Jensen's inequality that
\[
    \mathbb E\widetilde{\mathbb E}^{1}
    \left[
        \mathbf1_{\{r_t=r_s\}}
        |A_{s,t}^{\varepsilon,\phi}|^2
    \right]\le C_{T,\mu,R}\left(
        |t-s|^{2\alpha_0}
        +|t-s|^{2H\beta_0}
        +|t-s|^{2H\gamma_0}
        +|t-s|^{2H}
    \right).
\]
Thus, using
$\widetilde{\mathbb E}^{1}|\widetilde G_t^{\varepsilon,\phi}|^2
=\mathbb E^1|G_t^{\varepsilon,\phi}|^2
\le C_T\|\phi\|_{L^2_\mu} ^2$
for $\mathbb P^0$-a.e.\ common environment, we get
\[
     \mathbb E\left[
        \mathbf1_{\{r_t=r_s\}}
        \left|\widetilde{\mathbb E}^{1}
        [A_{s,t}^{\varepsilon,\phi}
         \widetilde G_t^{\varepsilon,\phi}]
        \right|^2
    \right]\le
    C_{T,\mu,R}\|\phi\|_{L^2_\mu} ^2
    \left(
        |t-s|^{2\alpha_0}
        +|t-s|^{2H\beta_0}
        +|t-s|^{2H\gamma_0}
        +|t-s|^{2H}
    \right).
\]
According to the
boundedness of $D^Lb$ and
\eqref{eq:G-time-increment}, the second term on the right-hand side of \eqref{F-F} is bounded by
\(
    C_T
    \|\phi\|_{L^2_\mu} ^2
    |t-s|^2
\).
Combining the preceding estimates yields
\[
     \mathbb E
    \left|
        F_t^{\varepsilon,\phi}
        -
        F_s^{\varepsilon,\phi}
    \right|^2
    \le
    C_{T,\mu,R}
    \|\phi\|_{L^2_\mu} ^2
    \Big(
        |t-s|^{2\alpha_0}
        +
        |t-s|^{2H\beta_0}
        +
        |t-s|^{2H\gamma_0}
        +
        |t-s|^{2H}
        +
        |t-s|
    \Big).
\]
By the choice of $\kappa$ in
\eqref{eq:kappa-switching-control}, all the powers on the
right-hand side dominate $2\kappa$. For every $R>0$ and
$\phi\in  L^2(\mathbb R^d \to \mathbb R^d , \mu)$
with $\|\phi\|_{L^2_\mu} \le R$, the preceding estimates yield
\begin{equation*}
    \sup_{\varepsilon\in[0,1]}
    \mathbb E|F_t^{\varepsilon,\phi}-F_s^{\varepsilon,\phi}|^2
    \le C_{T,\mu,R}\|\phi\|_{L^2_\mu} ^2
    |t-s|^{2\kappa}.
\end{equation*}
Consequently,
\begin{equation}
\label{eq:control-bracket-bound}
    \sup_{\varepsilon\in[0,1]}
    \sup_{t\in[s_0,T]}
    \mathbb E\left|
        \frac{G_t^{\varepsilon,\phi}}{T-s_0}
        +\frac{t-s_0}{T-s_0}F_t^{\varepsilon,\phi}
    \right|^2\le C_{T,s_0}\|\phi\|_{L^2_\mu} ^2.
\end{equation}
Moreover, for $s_0<s<t\le T$,
\begin{equation}
\label{eq:control-bracket-increment}
    \sup_{\varepsilon\in[0,1]}
    \mathbb E\Bigg|
        \frac{G_t^{\varepsilon,\phi}}{T-s_0}
        +\frac{t-s_0}{T-s_0}F_t^{\varepsilon,\phi}
        -\frac{G_s^{\varepsilon,\phi}}{T-s_0}
        -\frac{s-s_0}{T-s_0}F_s^{\varepsilon,\phi}
    \Bigg|^2\le C_{T,s_0,\mu,R}
    \|\phi\|_{L^2_\mu} ^2|t-s|^{2\kappa}.
\end{equation}

\noindent \textit{Step 2: fractional regularity under switching.}
By Assumption~\ref{ass:regularity-sigma}, we know
\begin{equation}
\label{eq:sigma-switching-increment}
    \left\|
        \sigma^{-1}(t,r_t)
        -
        \sigma^{-1}(s,r_s)
    \right\|^2
    \le
    C_T
    \left(
        |t-s|^{2\zeta_0}
        +
        \mathbf 1_{\{r_t\ne r_s\}}
    \right).
\end{equation}
Indeed, on $\{r_t=r_s\}$ this follows from
\eqref{eq:sigma-inverse-holder}, whereas on
$\{r_t\ne r_s\}$ it follows from the uniform boundedness of
$\sigma^{-1}$.
Since the bracketed term is $\mathscr G_s$-measurable,
the Markov property and the independence of
$r$ from $(\xi,W)$ give
\[
    \mathbb E\left[
        \mathbf1_{\{r_t\ne r_s\}}\mid\mathscr G_s
    \right]
    =\mathbb P^0(r_t\ne r_s\mid\mathscr F_s^0)\le\lambda_Q(t-s).
\]
Hence \eqref{eq:control-bracket-bound} and
\eqref{eq:sigma-switching-increment}
 imply
\begin{equation}
\label{eq:sigma-control-increment}
\begin{aligned}
    & \quad \mathbb E
    \Bigg[
        \left\|
            \sigma^{-1}(t,r_t)
            -
            \sigma^{-1}(s,r_s)
        \right\|^2
        \left|
            \frac{1}{T-s_0}
            G_s^{\varepsilon,\phi}
            +
            \frac{s-s_0}{T-s_0}
            F_s^{\varepsilon,\phi}
        \right|^2
    \Bigg] \\
    & \le
    C_{T,s_0}
    \|\phi\|_{L^2_\mu} ^2
    \left(
        |t-s|^{2\zeta_0}
        +
        |t-s|
    \right).
\end{aligned}
\end{equation}
By \eqref{eq:explicit-control-epsilon}, for a.e. $t>s_0$, we obtain
\begin{equation}
\label{eq:control-derivative}
    \frac{\mathrm d}{\mathrm dt}
    (R_Hh_{s_0}^{\varepsilon,\phi})(t)
    =
    \sigma^{-1}(t,r_t)
    \left[
        \frac{1}{T-s_0}
        G_t^{\varepsilon,\phi}
        +
        \frac{t-s_0}{T-s_0}
        F_t^{\varepsilon,\phi}
    \right].
\end{equation}
For every $R>0$ and
$\phi\in  L^2(\mathbb R^d \to \mathbb R^d , \mu)$
with $\|\phi\|_{L^2_\mu} \le R$, we obtain
\begin{equation}
\label{eq:control-derivative-bound}
    \sup_{\varepsilon\in[0,1]}
    \sup_{t\in(s_0,T]}
    \mathbb E\left|
        \frac{\mathrm d}{\mathrm dt}
        (R_Hh_{s_0}^{\varepsilon,\phi})(t)
    \right|^2\le C_{T,s_0}\|\phi\|_{L^2_\mu} ^2.
\end{equation}
For $s_0<s<t\le T$,  we have 
\begin{align*}
    & \quad \frac{\mathrm d}{\mathrm dt}
    (R_Hh_{s_0}^{\varepsilon,\phi})(t)
    -
    \frac{\mathrm d}{\mathrm ds}
    (R_Hh_{s_0}^{\varepsilon,\phi})(s)\\
    & =
    \left[
        \sigma^{-1}(t,r_t)
        -
        \sigma^{-1}(s,r_s)
    \right]
    \left[
        \frac{1}{T-s_0}
        G_s^{\varepsilon,\phi}
        +
        \frac{s-s_0}{T-s_0}
        F_s^{\varepsilon,\phi}
    \right]
    \\
    & \quad +
    \sigma^{-1}(t,r_t)
    \Bigg[
        \frac{1}{T-s_0}
        G_t^{\varepsilon,\phi}
        +
        \frac{t-s_0}{T-s_0}
        F_t^{\varepsilon,\phi}
        -
        \frac{1}{T-s_0}
        G_s^{\varepsilon,\phi}
        -
        \frac{s-s_0}{T-s_0}
        F_s^{\varepsilon,\phi}
    \Bigg].
\end{align*}
Combining the preceding decomposition with
\eqref{eq:sigma-control-increment} and
\eqref{eq:control-bracket-increment}, we obtain
\begin{equation}
\label{eq:control-derivative-increment}
    \sup_{\varepsilon\in[0,1]}
    \mathbb E\left|
        \frac{\mathrm d}{\mathrm dt}
        (R_Hh_{s_0}^{\varepsilon,\phi})(t)
        -\frac{\mathrm d}{\mathrm ds}
        (R_Hh_{s_0}^{\varepsilon,\phi})(s)
    \right|^2
    \le C_{T,s_0,\mu,R}
    \|\phi\|_{L^2_\mu} ^2|t-s|^{2\kappa}.
\end{equation}

We next estimate the inverse Volterra transform.
For an absolutely continuous function $u$ with
$u(0)=0$ and $u\in\operatorname{Ran}(K_H)$,
the inverse Volterra transform is given by
\begin{equation}
\label{eq:KH-inverse-AC}
    (K_H^{-1}u)(t)
    =
    t^{H-\frac12}
    D_{0+}^{H-\frac12}
    \left[
        t^{\frac12-H}u'(t)
    \right],
\end{equation}
where $D_{0+}^{H-\frac12}$ is the left-sided
Riemann--Liouville fractional derivative; see, for example,
\cite[Section~2]{MR3646611} and
\cite[Section~2]{MR4441502}.
For $\alpha\in(0,1)$, its representation is
\[
    D_{0+}^{\alpha}f(t)
    :=
    \frac{1}{\Gamma(1-\alpha)}
    \left[
        \frac{f(t)}{t^\alpha}
        +
        \alpha
        \int_0^t
        \frac{
            f(t)-f(s)
        }{
            (t-s)^{\alpha+1}
        }
        \mathrm ds
    \right].
\]
For $\alpha>0$, the left-sided
Riemann--Liouville fractional integral is defined by
\[
    I_{0+}^{\alpha}f(t)
    :=
    \frac{1}{\Gamma(\alpha)}
    \int_0^t
    (t-s)^{\alpha-1}f(s) \mathrm ds.
\]

Set
\(
    u(t)
    :=
    (R_Hh_{s_0}^{\varepsilon,\phi})(t)
\).
Since $u(t)=0$ for $t\le s_0$, we also have
$u'(t)=0$ for a.e. $t<s_0$.
We first interpret the right-hand side of
\eqref{eq:KH-inverse-AC} as a candidate inverse,
whose validity will be verified below.
Setting $\alpha=H-\frac12$ and splitting the
Weyl integral at $s_0$, we obtain, for $t>s_0$,
\begin{align}
    (K_H^{-1}u)(t)
    &=
    \frac{1}{\Gamma(\frac32-H)}
    \Bigg\{
        (t-s_0)^{\frac12-H}u'(t)
     +
        \left(H-\frac12\right)t^{H-\frac12}u'(t)
        \int_{s_0}^t
        \frac{
            t^{\frac12-H}-s^{\frac12-H}
        }{
            (t-s)^{H+\frac12}
        }
        \mathrm ds \nonumber
    \\
    &\quad \quad \quad \quad \quad \quad +
        \left(H-\frac12\right)t^{H-\frac12}
        \int_{s_0}^t
        \frac{
            s^{\frac12-H}
            \left(u'(t)-u'(s)\right)
        }{
            (t-s)^{H+\frac12}
        }
        \mathrm ds
    \Bigg\}.  \label{eq:KH-inverse-switching-decomposition}
\end{align}
For $t\le s_0$, the left-hand side is zero.

We estimate the three terms in the right side of 
\eqref{eq:KH-inverse-switching-decomposition} separately. For the first term, by
\eqref{eq:control-derivative-bound}, we have
\[
    \int_{s_0}^T
    (t-s_0)^{1-2H}
    \mathbb E
    \left|
        \frac{\mathrm d}{\mathrm dt}
        (R_Hh_{s_0}^{\varepsilon,\phi})(t)
    \right|^2
    \mathrm dt
    \le
    C_{T,s_0}
    \|\phi\|_{L^2_\mu} ^2
    \int_{s_0}^T
    (t-s_0)^{1-2H}
    \mathrm dt
    <\infty,
\] 
because of 
\(
    1-2H>-1
\).
For the second term, if $s_0=0$, the change of variables $s=tr$ implies
\[
    t^{H-\frac12}
    \left|
        \int_0^t
        \frac{
            t^{\frac12-H}
            -
            s^{\frac12-H}
        }{
            (t-s)^{H+\frac12}
        }
        \mathrm ds
    \right|
    =
    t^{\frac12-H}
    \left|
        \int_0^1
        \frac{
            1-r^{\frac12-H}
        }{
            (1-r)^{H+\frac12}
        }
        \mathrm dr
    \right|
    \le
    C_Ht^{\frac12-H}.
\]
If $s_0>0$, the mean value theorem gives $\left|
        t^{\frac12-H}
        -
        s^{\frac12-H}
    \right|
    \le
    C_{H,s_0}|t-s|$ for
$s_0\le s<t\le T$.
Hence
\[
    t^{H-\frac12}
    \left|
        \int_{s_0}^t
        \frac{
            t^{\frac12-H}
            -
            s^{\frac12-H}
        }{
            (t-s)^{H+\frac12}
        }
        \mathrm ds
    \right|
    \le
    C_{H,T,s_0}
    \int_{s_0}^t
    (t-s)^{\frac12-H}
    \mathrm ds
    \le
    C_{H,T,s_0}.
\]
Thus the second term is square-integrable. For the third term,
by Minkowski's inequality and
\eqref{eq:control-derivative-increment}, we have
\begin{align*}
    & \quad \left(
        \mathbb E
        \left|
            t^{H-\frac12}
            \int_{s_0}^t
            \frac{
                s^{\frac12-H}
            }{
                (t-s)^{H+\frac12}
            }
            \Bigg[
                \frac{\mathrm d}{\mathrm dt}
                (R_Hh_{s_0}^{\varepsilon,\phi})(t)
                -
                \frac{\mathrm d}{\mathrm ds}
                (R_Hh_{s_0}^{\varepsilon,\phi})(s)
            \Bigg]
            \mathrm ds
        \right|^2
    \right)^{1/2}
    \\
    &\le
    C_{T,s_0,\mu,R}
    \|\phi\|_{L^2_\mu} 
    t^{H-\frac12}
    \int_{s_0}^t
    s^{\frac12-H}
    (t-s)^{\kappa-H-\frac12}
    \mathrm ds.
\end{align*}
If $s_0=0$, the beta-integral identity gives
\[
    t^{H-\frac12}
    \int_0^t s^{\frac12-H}
       (t-s)^{\kappa-H-\frac12} \mathrm ds
    =
    B\left(\frac32-H,\kappa-H+\frac12\right)
    t^{\kappa-H+\frac12}.
\]
If $s_0>0$, then
\[
    t^{H-\frac12}
    \int_{s_0}^t s^{\frac12-H}
        (t-s)^{\kappa-H-\frac12} \mathrm ds \le C_{T,s_0}
        (t-s_0)^{\kappa-H+\frac12}.
\]
Since $\kappa>H-\frac12$, both bounds are square-integrable.
Combining the estimates of the three terms in
\eqref{eq:KH-inverse-switching-decomposition},  for every $R>0$ and
$\phi\in  L^2(\mathbb R^d \to \mathbb R^d , \mu)$
with $\|\phi\|_{L^2_\mu} \le R$, we obtain
\begin{equation}
\label{eq:KH-L2-local}
    \sup_{\varepsilon\in[0,1]}
    \mathbb E
    \int_0^T
    \left|
        K_H^{-1}
        (R_Hh_{s_0}^{\varepsilon,\phi})(t)
    \right|^2
    \mathrm dt
    \le
    C_{T,s_0,\mu,R}
    \|\phi\|_{L^2_\mu} ^2.
\end{equation}

It remains to verify that the candidate inverse
indeed belongs to the domain of the Volterra
inversion. Put $\alpha=H-\frac12$.
By \eqref{eq:control-derivative-bound} and
$\alpha<1/2$, we have
\[
    \mathbb E\int_0^T
        s^{-2\alpha}|u'(s)|^2 \mathrm ds\le
    C_{T,s_0}\|\phi\|_{L^2_\mu} ^2
    \int_0^T s^{-2\alpha} \mathrm ds
    <\infty.
\]
Consequently,
$s^{-\alpha}u'(s)\in L^2(0,T)$ almost surely.
By the Cauchy--Schwarz inequality, we have
\[
    \left|
        I_{0+}^{1-\alpha}
        \bigl(s^{-\alpha}u'(s)\bigr)(t)
    \right| \le
    \frac{t^{\frac12-\alpha}}
    {\Gamma(1-\alpha)\sqrt{1-2\alpha}}
    \left(
        \int_0^T
        s^{-2\alpha}|u'(s)|^2 \mathrm ds
    \right)^{1/2}
    \longrightarrow0
\]
as $t\downarrow0$, almost surely.
Thus the boundary term in the fractional
inversion formula vanishes. Together with
the square-integrability estimate
\eqref{eq:KH-L2-local} for the candidate inverse,
the fractional inversion identity and
the range characterization of $K_H$ yield $u\in\operatorname{Ran}(K_H)$ and $K_H^{-1}u\in L^2([0,T]\times\Omega;\mathbb R^d)$.
This justifies
\eqref{eq:KH-inverse-switching-decomposition}.

\noindent \textit{Step 3: admissibility and divergence estimate.}

By \eqref{eq:control-derivative} and
\eqref{eq:KH-inverse-switching-decomposition},
the process
\(
    K_H^{-1}(R_Hh_{s_0}^{\varepsilon,\phi})
\)
is $\{\mathscr G_t\}_{t\in[0,T]}$-adapted.
Together with \eqref{eq:KH-L2-local} and the fact that
$W$ remains a Brownian motion under this filtration,
this implies
\(
    K_H^{-1}(R_Hh_{s_0}^{\varepsilon,\phi}) \in \operatorname{Dom}\delta^W
\).
Since
\(
    K_H^*h_{s_0}^{\varepsilon,\phi}
    =
    K_H^{-1}(R_Hh_{s_0}^{\varepsilon,\phi})
\),
the divergence transfer relation yields
\(
    h_{s_0}^{\varepsilon,\phi}
    \in\operatorname{Dom}\delta
\)
and
\begin{equation}
\label{eq:control-divergence-Ito}
    \delta(h_{s_0}^{\varepsilon,\phi})
    =
    \int_0^T
    \left\langle
        K_H^{-1}
        (R_Hh_{s_0}^{\varepsilon,\phi})(t),
        \mathrm dW_t
    \right\rangle.
\end{equation}
Therefore, by the It\^o isometry and
\eqref{eq:KH-L2-local}, for every $R>0$ and
$\phi\in  L^2(\mathbb R^d \to \mathbb R^d , \mu)$
with $\|\phi\|_{L^2_\mu} \le R$, we obtain
\begin{equation}
\label{eq:control-local-uniform-estimate}
    \sup_{\varepsilon\in[0,1]}
    \mathbb E
    \left|
        \delta(h_{s_0}^{\varepsilon,\phi})
    \right|^2
    \le
    C_{T,s_0,\mu,R}
    \|\phi\|_{L^2_\mu} ^2.
\end{equation}
For any fixed
$\phi\in  L^2(\mathbb R^d \to \mathbb R^d , \mu)$,
taking $R>\|\phi\|_{L^2_\mu} $ yields
\[
    \int_0^1
    \left(
        \mathbb E
        \left|
            \delta(h_{s_0}^{\varepsilon,\phi})
        \right|^2
    \right)^{1/2}
    \mathrm d\varepsilon
    \le
    C_{T,s_0,\mu,R}^{1/2}
    \|\phi\|_{L^2_\mu} 
    <\infty,
\]
which proves \eqref{eq:control-integrability-epsilon}.
Moreover, for fixed $s_0\in[0,T)$ and
$\phi\in  L^2(\mathbb R^d \to \mathbb R^d , \mu)$,
the Picard constructions of
$X^{\varepsilon,\phi}$ and $G^{\varepsilon,\phi}$,
together with \eqref{eq:epsilon-law-variation}, admit versions
jointly measurable in $(\varepsilon,t,\omega)$.
Consequently, \eqref{eq:explicit-control-epsilon} and
\eqref{eq:KH-inverse-switching-decomposition} give a jointly measurable,
$\{\mathscr G_t\}$-progressively measurable version of
$K_H^{-1}(R_Hh_{s_0}^{\varepsilon,\phi})(t)$.
By \eqref{eq:KH-L2-local}, we have
\[
    \int_0^1\mathbb E\int_0^T
    \left|
        K_H^{-1}(R_Hh_{s_0}^{\varepsilon,\phi})(t)
    \right|^2
    \mathrm dt \mathrm d\varepsilon<\infty.
\]
Approximating this integrand by jointly measurable simple
predictable processes and using the It\^o isometry yields a jointly measurable
version of
$(\varepsilon,\omega)\mapsto
\delta(h_{s_0}^{\varepsilon,\phi})(\omega)$.
This verifies the measurability requirement in
\eqref{eq:abstract-integrability}.

It remains to establish
\eqref{eq:control-estimate-zero}.
When $\varepsilon=0$, the solution $X$
and its conditional law flow are independent of
the direction $\phi$.
Since \eqref{eq:epsilon-variation-equation}
is linear in $G^{0,\phi}$ with initial value
$\phi(\xi)$, the processes
$G^{0,\phi}$, $F^{0,\phi}$, and
$h_{s_0}^{0,\phi}$ depend linearly on $\phi$.
For any nonzero $\phi$, set
\(
    \widehat\phi
    :=
    \frac{\phi}{\|\phi\|_{L^2_\mu} }
\).
Then $\|\widehat\phi\|_{L^2_\mu} =1$, and
the linearity of the control and the divergence
operator gives
\(
    \delta(h_{s_0}^{0,\phi})
    =
    \|\phi\|_{L^2_\mu} 
    \delta(h_{s_0}^{0,\widehat\phi})
\).
Applying
\eqref{eq:control-local-uniform-estimate}
with $R=1$, we obtain
\[
    \mathbb E
    \left|
        \delta(h_{s_0}^{0,\phi})
    \right|^2
    =
    \|\phi\|_{L^2_\mu} ^2
    \mathbb E
    \left|
        \delta(h_{s_0}^{0,\widehat\phi})
    \right|^2
    \le
    C_{T,s_0,\mu,1}
    \|\phi\|_{L^2_\mu} ^2.
\]
The same estimate holds trivially for $\phi=0$.
Thus \eqref{eq:control-estimate-zero}
follows with
$C_{T,s_0,\mu}:=C_{T,s_0,\mu,1}$.
\end{proof}

We next verify that the above Cameron--Martin direction produces the
required terminal variation.

\begin{lemma}
\label{lem:terminal-matching}

Under the assumptions of
Lemma~\ref{lem:switching-cameron-martin}, for every
$\varepsilon\in[0,1]$ one has
\begin{equation}
\label{eq:matching-process-epsilon}
    \mathbb D_{R_Hh_{s_0}^{\varepsilon,\phi}}
    X_t^{\varepsilon,\phi}
    =
    \frac{t-s_0}{T-s_0}
    G_t^{\varepsilon,\phi},
    \quad
    t\in[s_0,T].
\end{equation}
In particular,
\begin{equation}
\label{eq:terminal-matching-epsilon}
    \mathbb D_{R_Hh_{s_0}^{\varepsilon,\phi}}
    X_T^{\varepsilon,\phi}
    =
    \nabla_{\phi(\xi)}
    X_T^{\varepsilon,\phi}.
\end{equation}

\end{lemma}

\begin{proof}

By Lemma~\ref{lem:malliavin-derivative}, for $t\in[s_0,T]$, we know
\[
    \mathbb D_{R_Hh_{s_0}^{\varepsilon,\phi}}
    X_t^{\varepsilon,\phi}
    =
    \int_{s_0}^t
    \nabla_x b
    \left(
        s,
        X_s^{\varepsilon,\phi},
        \mu_s^{\varepsilon,\phi},
        r_s
    \right)
    \mathbb D_{R_Hh_{s_0}^{\varepsilon,\phi}}
    X_s^{\varepsilon,\phi}
    \mathrm ds
    +
    \int_{s_0}^t
    \sigma(s,r_s)
    \mathrm d
    (R_Hh_{s_0}^{\varepsilon,\phi})(s).
\]
According to the definition
\eqref{eq:explicit-control-epsilon}, we get
\[
    \mathbb D_{R_Hh_{s_0}^{\varepsilon,\phi}}
    X_t^{\varepsilon,\phi}
    =
    \int_{s_0}^t
    \nabla_x b
    \left(
        s,
        X_s^{\varepsilon,\phi},
        \mu_s^{\varepsilon,\phi},
        r_s
    \right)
    \mathbb D_{R_Hh_{s_0}^{\varepsilon,\phi}}
    X_s^{\varepsilon,\phi}
    \mathrm ds
    +
    \int_{s_0}^t
    \left[
        \frac{1}{T-s_0}
        G_s^{\varepsilon,\phi}
        +
        \frac{s-s_0}{T-s_0}
        F_s^{\varepsilon,\phi}
    \right]
    \mathrm ds.
\]
On the other hand, by
\eqref{eq:epsilon-variation-equation}, we have
\[
    \frac{t-s_0}{T-s_0}
    G_t^{\varepsilon,\phi}
    =
    \int_{s_0}^t
    \nabla_x b
    \left(
        s,
        X_s^{\varepsilon,\phi},
        \mu_s^{\varepsilon,\phi},
        r_s
    \right)
    \frac{s-s_0}{T-s_0}
    G_s^{\varepsilon,\phi}
    \mathrm ds
    +
    \int_{s_0}^t
    \left[
        \frac{1}{T-s_0}
        G_s^{\varepsilon,\phi}
        +
        \frac{s-s_0}{T-s_0}
        F_s^{\varepsilon,\phi}
    \right]
    \mathrm ds.
\]
Both processes therefore satisfy the same linear integral equation on
$[s_0,T]$ with zero initial value at $s_0$. By uniqueness, we get
\[
    \mathbb D_{R_Hh_{s_0}^{\varepsilon,\phi}}
    X_t^{\varepsilon,\phi}
    =
    \frac{t-s_0}{T-s_0}
    G_t^{\varepsilon,\phi},
    \quad
    t\in[s_0,T].
\]
Taking $t=T$ and using
\eqref{eq:epsilon-variation}, we conclude
\(
    \mathbb D_{R_Hh_{s_0}^{\varepsilon,\phi}}
    X_T^{\varepsilon,\phi}
    =
    G_T^{\varepsilon,\phi}
    =
    \nabla_{\phi(\xi)}
    X_T^{\varepsilon,\phi}
\),
which completes the proof.
\end{proof}

\begin{lemma}
\label{lem:bismut-weight-continuity}
Suppose that Assumptions~\ref{ass:wellposedness}, \ref{ass:regularity-b} and \ref{ass:regularity-sigma} hold. Then, for every
$\phi\in  L^2(\mathbb R^d \to \mathbb R^d , \mu)$,
\begin{equation}
\label{eq:bismut-weight-epsilon-continuity}
    \lim_{\varepsilon\downarrow0}
    \mathbb E\left|
        \delta(h_0^{\varepsilon,\phi})
        -\delta(h_0^{0,\phi})
    \right|=0.
\end{equation}
Moreover,
\begin{equation}
\label{eq:bismut-weight-uniform-continuity}
    \lim_{\|\phi\|_{L^2_\mu} \to0}
    \sup_{\varepsilon\in(0,1]}
    \frac{
        \mathbb E\left|
            \delta(h_0^{\varepsilon,\phi})
            -\delta(h_0^{0,\phi})
        \right|
    }{\|\phi\|_{L^2_\mu} }=0.
\end{equation}
\end{lemma}

\begin{proof}
We first prove \eqref{eq:bismut-weight-uniform-continuity}.
Fix $\phi\ne0$ with $\|\phi\|_{L^2_\mu} \le1$ and
$\varepsilon\in[0,1]$. Throughout this proof, write $\Delta G_t
    :=G_t^{\varepsilon,\phi}-G_t^{0,\phi}$ and $\Delta F_t
    :=F_t^{\varepsilon,\phi}-F_t^{0,\phi}$.
All estimates below are uniform in $\varepsilon\in[0,1]$
and $0<\|\phi\|_{L^2_\mu}  \le1$.

By \eqref{eq:regularity-model}, we know
\begin{align*}
    X_t^{\varepsilon,\phi}-X_t
    &=\varepsilon\phi(\xi)
      +\int_0^t\Big[
          b(s,X_s^{\varepsilon,\phi},
            \mu_s^{\varepsilon,\phi},r_s)
          -b(s,X_s,\mu_s,r_s)
        \Big]  \mathrm ds.
\end{align*}
For $\mathbb P^0$-a.e. fixed common environment,
the pair $(X_s^{\varepsilon,\phi},X_s)$
under $\mathbb P^1$ couples the corresponding conditional
laws. Hence
\[
    W_2^2(\mu_s^{\varepsilon,\phi},\mu_s)
    \le\mathbb E^1
       |X_s^{\varepsilon,\phi}-X_s|^2,
    \quad\mathbb P^0\text{-a.s.}
\]
The Lipschitz continuity of $b$ and the Cauchy--Schwarz
inequality therefore yield
\[
    \mathbb E^1\sup_{u\in[0,t]}
      |X_u^{\varepsilon,\phi}-X_u|^2\le3\varepsilon^2\|\phi\|_{L^2_\mu}  ^2
      +6K_T^2T\int_0^t\mathbb E^1\sup_{u\in[0,s]}
        |X_u^{\varepsilon,\phi}-X_u|^2
       \mathrm ds.
\]
Applying Gronwall's inequality and then taking
$\mathbb E^0$ gives
\begin{equation*}
    \mathbb E\sup_{t\in[0,T]}
    \left|X_t^{\varepsilon,\phi}-X_t\right|^2
    \le C_T\varepsilon^2\|\phi\|_{L^2_\mu}  ^2.
\end{equation*}
In particular,
\[
    \mathbb E^0W_2^2
    (\mu_t^{\varepsilon,\phi},\mu_t)
    \le
    \mathbb E\left|X_t^{\varepsilon,\phi}-X_t\right|^2
    \le C_T\|\phi\|_{L^2_\mu}  ^2.
\]
It follows from Assumption~\ref{ass:regularity-b}(ii), Jensen's inequality,
and the preceding estimates that
\[
     \sup_{t\in[0,T]}\mathbb E\Bigl\|
        \nabla_xb(t,X_t^{\varepsilon,\phi},
                       \mu_t^{\varepsilon,\phi},r_t)
        -\nabla_xb(t,X_t,\mu_t,r_t)
    \Bigr\|^2
    \le C_T(\|\phi\|_{L^2_\mu}  ^{2\beta_0}+\|\phi\|_{L^2_\mu}  ^2),
\]
and
\[
     \sup_{t\in[0,T]}
    \mathbb E\widetilde{\mathbb E}^1
    \Bigl\|
        D^Lb(t,X_t^{\varepsilon,\phi},
               \mu_t^{\varepsilon,\phi},r_t)
             (\widetilde X_t^{\varepsilon,\phi})
        -D^Lb(t,X_t,\mu_t,r_t)
             (\widetilde X_t)
    \Bigr\|^2  \le C_T(\|\phi\|_{L^2_\mu}  ^{2\beta_0}
              +\|\phi\|_{L^2_\mu}  ^{2\gamma_0}+\|\phi\|_{L^2_\mu}  ^2).
\]

To estimate $\Delta F_t$, subtract the two expressions in
\eqref{eq:epsilon-law-variation}.
By the boundedness of $D^Lb$, the Cauchy--Schwarz inequality
on the original and copied probability spaces, and
\eqref{eq:G-epsilon-bound},
\[
    \frac{\mathbb E|\Delta F_t|}{\|\phi\|_{L^2_\mu}  }
    \le
    K_T\frac{\mathbb E|\Delta G_t|}{\|\phi\|_{L^2_\mu}  }
    +C_T(\|\phi\|_{L^2_\mu}  ^{\beta_0}+\|\phi\|_{L^2_\mu}  ^{\gamma_0}+\|\phi\|_{L^2_\mu}  ).
\]
Subtracting the equations \eqref{eq:epsilon-variation-equation}
for $G^{\varepsilon,\phi}$ and $G^{0,\phi}$, and using the
boundedness of $\nabla_xb$, we obtain
\begin{align*}
    \frac{\mathbb E\sup_{u\in[0,t]}|\Delta G_u|}{\|\phi\|_{L^2_\mu}  }
    &\le
    C_T\int_0^t
    \left(
        \frac{\mathbb E\sup_{u\in[0,s]}|\Delta G_u|}{\|\phi\|_{L^2_\mu}  }
        +\frac{\mathbb E|\Delta F_s|}{\|\phi\|_{L^2_\mu}  }
    \right)\mathrm ds
    +C_T\|\phi\|_{L^2_\mu}  ^{\beta_0}.
\end{align*}
Consequently, applying Gronwall's inequality yields
\begin{equation}
\label{eq:weight-normalized-L1-variation}
    \sup_{\varepsilon\in[0,1]}
    \sup_{t\in[0,T]}
    \frac{\mathbb E|\Delta G_t|+\mathbb E|\Delta F_t|}
         {\|\phi\|_{L^2_\mu}  }
    \le
    C_T(\|\phi\|_{L^2_\mu}  ^{\beta_0}+\|\phi\|_{L^2_\mu}  ^{\gamma_0}+\|\phi\|_{L^2_\mu}  ).
\end{equation}

Fix $\phi\in L^2(\mathbb R^d \to \mathbb R^d , \mu)$ with
$0<\|\phi\|_{L^2_\mu}\le1$.
Define
\(
    v_t^{\varepsilon,\phi}
    :=
    \frac{\sigma^{-1}(t,r_t)}
    {T\|\phi\|_{L^2_\mu}}
    \left(
        \Delta G_t+t\Delta F_t
    \right)
\).
By \eqref{eq:explicit-control-epsilon} with $s_0=0$,
we have
\[
    R_Hh_0^{\varepsilon,\phi}(t)
    -R_Hh_0^{0,\phi}(t)=
    \frac1T\int_0^t
    \sigma^{-1}(s,r_s)
    \left(
        \Delta G_s+s\Delta F_s
    \right)
    \mathrm ds=
    \|\phi\|_{L^2_\mu}
    \int_0^t
    v_s^{\varepsilon,\phi}\mathrm ds.
\]
Since $\sigma^{-1}$ is bounded,
\eqref{eq:weight-normalized-L1-variation} implies, as $\|\phi\|_{L^2_\mu}\to0$,
\begin{equation}
\label{eq:weight-v-small}
    \sup_{\varepsilon\in[0,1]}
    \sup_{t\in[0,T]}
    \mathbb E|v_t^{\varepsilon,\phi}| \le
    C_T
    \left(
        \|\phi\|_{L^2_\mu}^{\beta_0}
        +
        \|\phi\|_{L^2_\mu}^{\gamma_0}
        +
        \|\phi\|_{L^2_\mu}
    \right)
    \longrightarrow0.
\end{equation}
Choose $\kappa$ as in
\eqref{eq:kappa-switching-control}.
Applying \eqref{eq:control-derivative-bound} and
\eqref{eq:control-derivative-increment}
separately to the two controls, with $R=1$,
we get
\begin{equation}
\label{eq:weight-v-increment}
    \sup_{\varepsilon\in[0,1]}
    \sup_{t\in[0,T]}
    \mathbb E|v_t^{\varepsilon,\phi}|^2
    \le C_{T,\mu}, \quad \text{and} \quad \mathbb E
    |v_t^{\varepsilon,\phi}
    -v_s^{\varepsilon,\phi}|^2
    \le C_{T,\mu}|t-s|^{2\kappa},
\end{equation}
where the second estimate is uniform in
$\varepsilon\in[0,1]$.

We first establish the convergence
\begin{equation}
\label{eq:weight-v-L2-probability}
    \lim_{\|\phi\|_{L^2_\mu}\to0}
    \sup_{\varepsilon\in[0,1]}
    \mathbb E
    \|v^{\varepsilon,\phi}\|_{L^2(0,T)}
    =0.
\end{equation}
For $m\ge1$, set $t_j=\frac{j}{m} T$ for $0\le j\le m$.
By decomposing $v^{\varepsilon,\phi}$ over
the time grid and applying the triangle inequality
in $L^2(0,T)$, we obtain
\begin{align*}
    \|v^{\varepsilon,\phi}\|_{L^2(0,T)}
    &=
    \left(
        \sum_{j=0}^{m-1}
        \int_{t_j}^{t_{j+1}}
        \left|
            v_t^{\varepsilon,\phi}
            -v_{t_j}^{\varepsilon,\phi}
            +v_{t_j}^{\varepsilon,\phi}
        \right|^2
        \mathrm dt
    \right)^{1/2}\\
    &\le
    \left(
        \sum_{j=0}^{m-1}
        \int_{t_j}^{t_{j+1}}
        |v_t^{\varepsilon,\phi}
        -v_{t_j}^{\varepsilon,\phi}|^2
        \mathrm dt
    \right)^{1/2}+
    \left(
        \frac{T}{m}
        \sum_{j=0}^{m-1}
        |v_{t_j}^{\varepsilon,\phi}|^2
    \right)^{1/2}\\
    &\le
    \left(
        \sum_{j=0}^{m-1}
        \int_{t_j}^{t_{j+1}}
        |v_t^{\varepsilon,\phi}
        -v_{t_j}^{\varepsilon,\phi}|^2
        \mathrm dt
    \right)^{1/2}
    +
    \sqrt T
    \max_{0\le j<m}|v_{t_j}^{\varepsilon,\phi}|.
\end{align*}
Taking expectations and using Jensen's inequality,
we get
\[
    \mathbb E
    \|v^{\varepsilon,\phi}\|_{L^2(0,T)}\le
    \left(
        \sum_{j=0}^{m-1}
        \int_{t_j}^{t_{j+1}}
        \mathbb E
        |v_t^{\varepsilon,\phi}
        -v_{t_j}^{\varepsilon,\phi}|^2
        \mathrm dt
    \right)^{1/2}+
    \sqrt T
    \sum_{j=0}^{m-1}
    \mathbb E|v_{t_j}^{\varepsilon,\phi}|.
\]
For the first term, by
\eqref{eq:weight-v-increment}, we have
\[
    \left(
        \sum_{j=0}^{m-1}
        \int_{t_j}^{t_{j+1}}
        \mathbb E
        |v_t^{\varepsilon,\phi}
        -v_{t_j}^{\varepsilon,\phi}|^2
        \mathrm dt
    \right)^{1/2}\le
    C_{T,\mu}
    \left(
        \sum_{j=0}^{m-1}
        \int_{t_j}^{t_{j+1}}
        |t-t_j|^{2\kappa}
        \mathrm dt
    \right)^{1/2}\le
    C_{T,\mu}
    \left(\frac{T}{m}\right)^\kappa.
\]
For the second term,
by \eqref{eq:weight-v-small} we have
\[
    \sqrt T
    \sum_{j=0}^{m-1}
    \mathbb E|v_{t_j}^{\varepsilon,\phi}|
    \le
    C_Tm
    \left(
        \|\phi\|_{L^2_\mu}^{\beta_0}
        +
        \|\phi\|_{L^2_\mu}^{\gamma_0}
        +
        \|\phi\|_{L^2_\mu}
    \right).
\]
Hence, for every fixed $m\ge1$, letting
$\|\phi\|_{L^2_\mu}\to0$ yields
\[
    \limsup_{\|\phi\|_{L^2_\mu}\to0}
    \sup_{\varepsilon\in[0,1]}
    \mathbb E
    \|v^{\varepsilon,\phi}\|_{L^2(0,T)}\le
    C_{T,\mu}
    \left(\frac{T}{m}\right)^\kappa.
\]
Since $\kappa>0$, letting $m\to\infty$
proves \eqref{eq:weight-v-L2-probability}.

We next prove that
\begin{equation}
\label{eq:weight-KH-probability}
    \lim_{\|\phi\|_{L^2_\mu}\to0}
    \sup_{\varepsilon\in[0,1]}
    \mathbb E
    \left\|
        K_H^{-1}
        \left(
            \int_0^\cdot
            v_s^{\varepsilon,\phi}\mathrm ds
        \right)
    \right\|_{L^2(0,T)}
    =0.
\end{equation}
Put $\alpha=H-\frac12$. Then
$0<\alpha<\kappa<\frac12$.
The Weyl representation used in
\eqref{eq:KH-inverse-switching-decomposition},
with $s_0=0$, implies, for almost every $t>0$,
\[
    \left|
        K_H^{-1}
        \left(
            \int_0^\cdot
            v_s^{\varepsilon,\phi}\mathrm ds
        \right)(t)
    \right|\le
    C_Ht^{-\alpha}|v_t^{\varepsilon,\phi}|+
    C_Ht^\alpha\int_0^t
    \frac{
        s^{-\alpha}
        |v_t^{\varepsilon,\phi}
        -v_s^{\varepsilon,\phi}|
    }{(t-s)^{1+\alpha}}\mathrm ds.
\]
For any $a\in(0,T)$, by
\eqref{eq:weight-v-increment}
and Jensen's inequality, we get
\begin{align*}
    \mathbb E
    \|t^{-\alpha}v_t^{\varepsilon,\phi}\|_{L^2(0,T)}& \le
    \left(
        \mathbb E\int_0^a t^{-2\alpha}|v_t^{\varepsilon,\phi}|^2
        \mathrm dt
    \right)^{1/2}+
    a^{-\alpha}
    \mathbb E
    \|v^{\varepsilon,\phi}\|_{L^2(0,T)} \\
    & \le
    C_{T,\mu}a^{\frac12-\alpha}
    +
    a^{-\alpha}
    \mathbb E\|v^{\varepsilon,\phi}\|_{L^2(0,T)}.
\end{align*}
For the integral term, split the integration domain into
$t-s<a$ and $t-s\ge a$.
By Minkowski's inequality and
\eqref{eq:weight-v-increment}, we have
\[
    \left\|
        t^\alpha\int_{(t-a)_+}^t
        \frac{
            s^{-\alpha}
            (v_t^{\varepsilon,\phi}
            -v_s^{\varepsilon,\phi})
        }{(t-s)^{1+\alpha}}\mathrm ds
    \right\|_{L^2(\Omega)}\le
    C_{T,\mu}t^\alpha
    \int_{(t-a)_+}^t
    s^{-\alpha}(t-s)^{\kappa-1-\alpha}
    \mathrm ds\le
    C_{T,\mu}a^{\kappa-\alpha}.
\]
For the last inequality, if $t\ge2a$, then
$s\ge t/2$ throughout the integration interval; if
$t<2a$, the substitution $s=tu$ reduces the integral
to a finite beta integral.
Thus the $L^2(\Omega\times[0,T])$ norm of the
near-diagonal part is bounded by
$C_{T,\mu}a^{\kappa-\alpha}$.
For fixed $a>0$, we have
\begin{align*}
    t^\alpha\int_0^{(t-a)_+}
    \frac{
        s^{-\alpha}
        |v_t^{\varepsilon,\phi}
        -v_s^{\varepsilon,\phi}|
    }{(t-s)^{1+\alpha}}\mathrm ds & \le
    C_{T,a}
    \left(
        |v_t^{\varepsilon,\phi}|
        +
        \int_0^T
        s^{-\alpha}|v_s^{\varepsilon,\phi}|
        \mathrm ds
    \right)\\
    & \le
    C_{T,a}
    \left(
        |v_t^{\varepsilon,\phi}|
        +
        \|v^{\varepsilon,\phi}\|_{L^2(0,T)}
    \right),
\end{align*}
where $s^{-\alpha}\in L^2(0,T)$ because
$\alpha<1/2$.
Hence the $L^1(\Omega;L^2(0,T))$ norm of the
part away from the diagonal is bounded by
\(
    C_{T,a}
    \mathbb E
    \|v^{\varepsilon,\phi}\|_{L^2(0,T)}
\).
Combining the preceding estimates, we obtain
\[
    \sup_{\varepsilon\in[0,1]}
    \mathbb E
    \left\|
        K_H^{-1}
        \left(
            \int_0^\cdot
            v_s^{\varepsilon,\phi}\mathrm ds
        \right)
    \right\|_{L^2(0,T)}\le
    C_{T,\mu}
    \left(
        a^{\frac12-\alpha}
        +
        a^{\kappa-\alpha}
    \right)+
    C_{T,\mu,a}
    \sup_{\varepsilon\in[0,1]}
    \mathbb E
    \|v^{\varepsilon,\phi}\|_{L^2(0,T)}.
\]
By \eqref{eq:weight-v-L2-probability},
first letting $\|\phi\|_{L^2_\mu}\to0$
and then $a\downarrow0$ proves
\eqref{eq:weight-KH-probability}, since
$1/2-\alpha>0$ and $\kappa-\alpha>0$.

By the transfer relation and the It\^o representation
\eqref{eq:control-divergence-Ito},
\[
    \frac{
        \delta(h_0^{\varepsilon,\phi})
        -
        \delta(h_0^{0,\phi})
    }{\|\phi\|_{L^2_\mu}}=
    \int_0^T
    \left\langle
        K_H^{-1}
        \left(
            \int_0^\cdot
            v_s^{\varepsilon,\phi}\mathrm ds
        \right)(t),
        \mathrm dW_t
    \right\rangle.
\]
The stochastic integral is well defined by
Lemma~\ref{lem:switching-cameron-martin}
and \eqref{eq:control-local-uniform-estimate}
with $R=1$.
Applying the Burkholder--Davis--Gundy inequality
with exponent $1$, we obtain
\[
    \frac{
        \mathbb E
        \left|
            \delta(h_0^{\varepsilon,\phi})
            -
            \delta(h_0^{0,\phi})
        \right|
    }{\|\phi\|_{L^2_\mu}}\le
    C\mathbb E
    \left\|
        K_H^{-1}
        \left(
            \int_0^\cdot
            v_s^{\varepsilon,\phi}\mathrm ds
        \right)
    \right\|_{L^2(0,T)}.
\]
Consequently, by \eqref{eq:weight-KH-probability}, we conclude
\begin{equation}
\label{eq:weight-normalized-probability}
    \lim_{\|\phi\|_{L^2_\mu}\to0}
    \sup_{\varepsilon\in[0,1]}
    \frac{
        \mathbb E
        \left|
            \delta(h_0^{\varepsilon,\phi})
            -
            \delta(h_0^{0,\phi})
        \right|
    }{\|\phi\|_{L^2_\mu}}
    =0.
\end{equation}
This proves
\eqref{eq:bismut-weight-uniform-continuity}.

Finally, fix $\phi\ne0$.
The initial perturbations corresponding to
$(\varepsilon,\phi)$ and $(1,\varepsilon\phi)$
coincide.
Linearity of the variational equation with respect
to its initial direction and of the control
construction gives $h_0^{1,\varepsilon\phi} =\varepsilon h_0^{\varepsilon,\phi}$ and $h_0^{0,\varepsilon\phi} =\varepsilon h_0^{0,\phi}$.
Consequently, by
\eqref{eq:bismut-weight-uniform-continuity},
\[
    \mathbb E
    \left|
        \delta(h_0^{\varepsilon,\phi})
        -
        \delta(h_0^{0,\phi})
    \right|=
    \|\phi\|_{L^2_\mu}
    \frac{
        \mathbb E
        \left|
            \delta(h_0^{1,\varepsilon\phi})
            -
            \delta(h_0^{0,\varepsilon\phi})
        \right|
    }{\|\varepsilon\phi\|_{L^2_\mu}}
    \longrightarrow0
    \quad \text{as }\varepsilon\downarrow0.
\]
For $\phi=0$ the assertion is immediate, proving
\eqref{eq:bismut-weight-epsilon-continuity}.
\end{proof}

\begin{theorem}[Bismut formula]
\label{thm:bismut-formula}

Suppose that Assumptions~\ref{ass:wellposedness}, \ref{ass:regularity-b} and \ref{ass:regularity-sigma} hold. 
Then, for every $T>0$, $\mu\in\mathcal P_2(\mathbb R^d)$,
$\phi\in L^2(\mathbb R^d \to \mathbb R^d , \mu)$, and
$f\in\mathscr B_b(\mathbb R^d)$, the functional
$\mu\mapsto(P_Tf)(\mu)$ is $L$-differentiable
at $\mu$, and
\begin{equation}
\label{eq:bismut-formula}
    D_\phi^L(P_Tf)(\mu)
    =
    \mathbb E
    \left[
        f(X_T)
        \int_0^T
        \left\langle
            K_H^{-1}(R_Hh_0^{0,\phi})(s),
            \mathrm dW_s
        \right\rangle
    \right],
\end{equation}
where $h_0^{0,\phi}$ is determined by
\eqref{eq:explicit-control-epsilon}
with $s_0=0$ and $\varepsilon=0$.
\end{theorem}

\begin{proof}
By Lemma~\ref{lem:switching-cameron-martin},
Lemma~\ref{lem:terminal-matching}, and
Lemma~\ref{lem:bismut-weight-continuity},
all the assumptions of
Theorem~\ref{thm:abstract-bismut}
are satisfied. Therefore, by
\eqref{eq:abstract-bismut-formula} we get
\[
    D_\phi^L(P_Tf)(\mu)
    =
    \mathbb E
    \left[
        f(X_T)\delta(h_0^{0,\phi})
    \right].
\]
By the divergence representation
\eqref{eq:control-divergence-Ito}, we have
\[
    \delta(h_0^{0,\phi})
    =
    \int_0^T
    \left\langle
        K_H^{-1}(R_Hh_0^{0,\phi})(s),
        \mathrm dW_s
    \right\rangle.
\]
Combining these identities proves
\eqref{eq:bismut-formula}.
\end{proof}

As a direct consequence, we obtain the following
gradient estimate (cf.~\cite[(4.66)]{MR4441502}).

\begin{corollary}
\label{cor:lions-gradient-estimate}

Under the assumptions of
Theorem~\ref{thm:bismut-formula}, for every
$T>0$ and $\mu\in\mathcal P_2(\mathbb R^d)$, there exists a constant
$C_{T,\mu}>0$ such that, for every
$f\in\mathscr B_b(\mathbb R^d)$,
\begin{equation}
\label{eq:lions-gradient-estimate}
    \|D^L(P_Tf)(\mu)\|_{L^2_\mu} 
    \le
    C_{T,\mu}
    \left[
        (P_Tf^2)(\mu)
        -
        \left((P_Tf)(\mu)\right)^2
    \right]^{1/2}.
\end{equation}

\end{corollary}

\begin{remark}
    (i) Unlike \cite[Theorem 4.7(I)]{MR4441502},
    our assumptions do not require the additional derivatives
    of $D^L b$ imposed by Assumption (A3) therein.
    In that work, these derivatives are used to expand
    temporal increments of the perturbed Lions-derivative
    term, whereas we control the inverse Volterra
    transform directly through first-order H\"older
    estimates, uniform $L^2$ increment bounds, and
    a time-splitting argument.

    (ii) At $H=1/2$, fractional Brownian motion
    reduces to standard Brownian motion. In the absence of
    Markovian switching, the Brownian counterpart of
    \eqref{eq:bismut-formula} is consistent with the classical
    distribution-dependent SDE setting, for which Ren and Wang
    \cite{MR3983053} derived a Bismut formula for the
    Lions derivative. When Markovian switching is retained
    but distribution dependence is removed, the corresponding
    setting is consistent with the classical switching diffusion
    framework studied by Hu et al.\ \cite{MR3986248}.

    (iii) For $H<1/2$, the inverse Volterra operator involves
    fractional integration, which may simplify the
    construction of Bismut controls. Indeed,
    \cite[Theorem~4.7]{MR4441502} covers this regime,
    but requires a constant invertible diffusion matrix.
    In contrast, our model allows regime-dependent diffusion
    coefficients, while Section~\ref{sec:wellposedness}
    further permits distribution-dependent diffusion coefficients.
    Extending the corresponding well-posedness estimates
    and Bismut formula to $H<1/2$ requires additional
    analysis of stochastic integration under Markovian switching.
\end{remark}

\section{Applications}\label{sec:applications}

We next investigate how perturbations of transition
rates affect the Lions derivatives of terminal
expectation functionals for conditional McKean--Vlasov SDEs.

\subsection{Likelihood-ratio representation}
\label{subsec:likelihood-ratio}

Let $Q=(q_{ij})_{i,j\in\mathbb S}$ be the conservative
transition rate matrix introduced in Section~\ref{sec:MF},
and fix the initial value $r_0=i_0$. Set 
\(
    \mathcal E:=\{(i,j)\in\mathbb S^2:i\ne j,\ q_{ij}>0\}
\).
We consider the class of conservative generators with the same
positive off-diagonal support,
\[
    \mathcal Q_{\mathcal E}:=
    \left\{\widetilde Q=(\widetilde q_{ij}):
    \begin{array}{ll}
       \widetilde q_{ij}>0,&(i,j)\in\mathcal E,\\
       \widetilde q_{ij}=0,&i\ne j,\ (i,j)\notin\mathcal E,\\
       \widetilde q_{ii}=-\displaystyle\sum_{j\ne i}
       \widetilde q_{ij},&i\in\mathbb S
    \end{array}
    \right\}.
\]
The corresponding space of admissible matrix directions is
\[
    \mathcal T_{\mathcal E}:=
    \left\{A=(a_{ij}): A\mathbf1=0,\
      a_{ij}=0\ \text{for }i\ne j,\ (i,j)\notin\mathcal E
    \right\}.
\]
Thus $\mathcal Q_{\mathcal E}$ is open in its finite-dimensional
affine span. The initial distribution and the coefficients
$b$ and $\sigma$ are kept fixed throughout the perturbation.
If $\mathcal E=\varnothing$, this parameter space is a singleton;
henceforth we assume $\mathcal E\ne\varnothing$.

Let $\mathfrak D_T=D([0,T];\mathbb S)$ be the Skorokhod path
space with Borel $\sigma$-field $\mathscr D_T$, and let $r$ be
its coordinate process. Write $r_t(\gamma)=\gamma(t)$ and $\mathscr D_t=\sigma(r_s:0\le s\le t)$.
For $\widetilde Q\in\mathcal Q_{\mathcal E}$, denote by
$\mathbb P_0^{\widetilde Q}$ the law of the Markov chain with
generator $\widetilde Q$ and initial state $i_0$. On
$\mathfrak D_T\times\Omega^1$, set
\(
    \mathbb P^{\widetilde Q}:=\mathbb P_0^{\widetilde Q}\otimes\mathbb P^1
\),
and write $\mathbb E_0^{\widetilde Q}$ and $\mathbb E^{\widetilde Q}$
for the corresponding expectations. At the reference matrix
$Q$, abbreviate $\mathbb P_0:=\mathbb P_0^Q$, $\mathbb P:=\mathbb P^Q$, $\mathbb E_0:=\mathbb E_0^Q$ and $\mathbb E:=\mathbb E^Q$. 
The reference canonical law coincides with the original
switching law after passage to the canonical path space of the switching process.

For $(i,j)\in\mathcal E$, define the transition count and
state occupation time by
\[
    N_{ij}(T)
    :=\sum_{0<s\le T}\mathbf1_{\{r_{s-}=i,\ r_s=j\}},
    \quad \text{and} \quad
    O_i(T):=\int_0^T\mathbf1_{\{r_s=i\}} \mathrm ds.
\]
Write $N_T:=\sum_{(i,j)\in\mathcal E}N_{ij}(T)$.
For the classical likelihood-ratio method for
continuous-time Markov chains, see
Glasserman~\cite[Section~6]{MR1162948} and
Nakayama et al.~\cite[Section~2]{Nakayama1994}.
Fix $\widetilde Q\in\mathcal Q_{\mathcal E}$ and set
$\Delta=\widetilde Q-Q\in\mathcal T_{\mathcal E}$.
Since the two generators have the same positive transition
edges and the same initial state, their laws on
$\mathfrak D_T$ are equivalent. The likelihood ratio is defined by
\begin{equation}
\label{eq:switching-likelihood-ratio}
    L_T^{\widetilde Q}
    :=\frac{\mathrm d\mathbb P_0^{\widetilde Q}}
            {\mathrm d\mathbb P_0^Q}
    =\prod_{(i,j)\in\mathcal E}
       \left(\frac{\widetilde q_{ij}}{q_{ij}}\right)^{N_{ij}(T)}
       \exp\left(-\sum_{(i,j)\in\mathcal E}
                 \Delta_{ij}O_i(T)\right).
\end{equation}
Consequently, for every nonnegative measurable random variable
$Z$ on the canonical product space, we know
\[
    \mathbb E^{\widetilde Q}[Z]
    =\mathbb E\big[ZL_T^{\widetilde Q}\big],
    \quad \text{and} \quad
    \mathbb E_0 L_T^{\widetilde Q}=1.
\]
For $A\in\mathcal T_{\mathcal E}$, the matrices
$Q+\theta A$ belong to $\mathcal Q_{\mathcal E}$ whenever
$|\theta|$ is sufficiently small. By differentiating the
logarithm of \eqref{eq:switching-likelihood-ratio}, we get
\begin{equation}
\label{eq:matrix-likelihood-score}
    S_A(T)
    :={\left.\frac{\mathrm d}{\mathrm d\theta}
        \log L_T^{Q+\theta A}\right|_{\theta=0}}=\sum_{(i,j)\in\mathcal E}
        a_{ij}\left(\frac{N_{ij}(T)}{q_{ij}}-O_i(T)\right).
\end{equation}
In particular, the score is linear in the matrix direction $A$.
For $(i,j)\in\mathcal E$ and $t\in[0,T]$, let
\[
    M_{ij}(t):=N_{ij}(t)-q_{ij}O_i(t),
\]
where $N_{ij}(t)$ and $O_i(t)$ are defined as above with
$T$ replaced by $t$. The compensated transition-counting
processes are
square-integrable martingales. Since distinct directed
transitions cannot occur at the same time, these martingales
are pairwise orthogonal. Their predictable quadratic
variations imply
\begin{equation}
\label{eq:matrix-score-second-moment}
    \mathbb E_0 S_A(T)=0, \quad  \text{and} \quad
    \mathbb E_0[S_A(T)^2]
    =\sum_{(i,j)\in\mathcal E}
      \frac{a_{ij}^2}{q_{ij}}
      \int_0^T\mathbb P_0(r_s=i) \mathrm ds.
\end{equation}
More generally, for $A,B\in\mathcal T_{\mathcal E}$, we have
\[
    \mathbb E_0[S_A(T)S_B(T)]
    =\sum_{(i,j)\in\mathcal E}
      \frac{a_{ij}b_{ij}}{q_{ij}}
      \int_0^T\mathbb P_0(r_s=i) \mathrm ds.
\]

We next establish quantitative estimates for the likelihood
ratio. Define
\[
    d_Q(Q,\widetilde Q)^2
    :=\max_{i\in\mathbb S}
       \sum_{j:(i,j)\in\mathcal E}
       \frac{(\widetilde q_{ij}-q_{ij})^2}{q_{ij}}.
\]
For $t\in[0,T]$, let $L_t^{\widetilde Q}$ be the
likelihood ratio on $\mathscr D_t$ given by
\eqref{eq:switching-likelihood-ratio} with $T$ replaced
by $t$. It follows from the counting-process likelihood-ratio formula that
\[
    \mathrm dL_t^{\widetilde Q}
    =L_{t-}^{\widetilde Q}
      \sum_{(i,j)\in\mathcal E}
      \frac{\Delta_{ij}}{q_{ij}} \mathrm dM_{ij}(t),
    \quad \text{and} \quad L_0^{\widetilde Q}=1.
\]
Since $\mathcal E$ is finite and $N_T$ has exponential
moments of every order, the likelihood-ratio process has
finite second moments. The jump It\^o formula, orthogonality
of the compensated counting martingales, and the compensation
formula give
\[
    \mathbb E_0[(L_t^{\widetilde Q})^2]
    =1+\int_0^t\mathbb E_0\left[
       (L_{s-}^{\widetilde Q})^2
       \sum_{j:(r_{s-},j)\in\mathcal E}
       \frac{\Delta_{r_{s-}j}^2}{q_{r_{s-}j}}
       \right] \mathrm ds \le 1+d_Q(Q,\widetilde Q)^2
       \int_0^t\mathbb E_0[(L_s^{\widetilde Q})^2]
        \mathrm ds.
\]
Gronwall's inequality implies
\(
    \mathbb E_0[(L_T^{\widetilde Q})^2]
    \le \exp\big(Td_Q(Q,\widetilde Q)^2\big)
\).
Using $\mathbb E_0 L_T^{\widetilde Q}=1$, we obtain
\begin{equation}
\label{eq:matrix-likelihood-l2}
    \big\|L_T^{\widetilde Q}-1\big\|_{L^2(\mathbb P_0)}
    \le\left(e^{Td_Q(Q,\widetilde Q)^2}-1\right)^{1/2}.
\end{equation}

Put
\(
    q_*:=\min_{(i,j)\in\mathcal E}q_{ij}>0
\)
and fix $\rho\in(0,q_*/2)$. Let
$\Delta\in\mathcal T_{\mathcal E}$ satisfy
$\|\Delta\|\le\rho$, where $\|\cdot\|$ is the
Frobenius norm used in Section~\ref{sec:MF}. For
$u\in[0,1]$, define
\[
    R_u:=\sum_{(i,j)\in\mathcal E}
        \Delta_{ij}\left(
          \frac{N_{ij}(T)}{q_{ij}+u\Delta_{ij}}-O_i(T)
        \right).
\]
By differentiating the likelihood ratio along the matrix segment, we get
\[
    \frac{\mathrm d}{\mathrm du}L_T^{Q+u\Delta}
    =L_T^{Q+u\Delta}R_u,
    \quad  \text{and} \quad
    \frac{\mathrm d^2}{\mathrm du^2}L_T^{Q+u\Delta}
    =L_T^{Q+u\Delta}\left[R_u^2-
       \sum_{(i,j)\in\mathcal E}
       \frac{\Delta_{ij}^2N_{ij}(T)}
            {(q_{ij}+u\Delta_{ij})^2}\right].
\]
In particular, $L_T^Q=1$ and
$\left.\frac{\mathrm d}{\mathrm du}L_T^{Q+u\Delta}\right|_{u=0}
=S_\Delta(T)$. Since
$q_{ij}+u\Delta_{ij}\ge q_*/2$,
$0\le O_i(T)\le T$, and $\mathcal E$ is finite,
there is a constant $C_Q>0$, independent of
$u$ and $\Delta$, such that
\[
    |R_u|\le C_Q\|\Delta\|(N_T+T),
    \quad  \text{and} \quad
    |L_T^{Q+u\Delta}|\le
      \exp\big(C_Q\rho(N_T+T)\big).
\]
Consequently, we obtain
\[
    \sup_{u\in[0,1]}
    \left|\frac{\mathrm d^2}{\mathrm du^2}L_T^{Q+u\Delta}\right|
    \le C_Q\|\Delta\|^2
       (1+N_T)^2\exp\big(C_Q\rho(N_T+T)\big).
\]
Under $\mathbb P_0$, the total number $N_T$ of jumps is
stochastically dominated by a Poisson random variable with
parameter $T\max_{i\in\mathbb S}(-q_{ii})$.
Thus the right-hand side has a finite $L^2(\mathbb P_0)$
norm, with a bound uniform over $\|\Delta\|\le\rho$. It follows from
Taylor's formula with integral remainder and Minkowski's
inequality that
\begin{equation}
\label{eq:matrix-likelihood-taylor}
    \big\|L_T^{Q+\Delta}-1-S_\Delta(T)\big\|_{L^2(\mathbb P_0)}
    \le C_{T,Q,\rho}\|\Delta\|^2,
    \quad \Delta\in\mathcal T_{\mathcal E},\ 
           \|\Delta\|\le\rho.
\end{equation}
The constant is independent of the direction of $\Delta$.
In particular, for every fixed $A\in\mathcal T_{\mathcal E}$, we have 
\[
    \lim_{\theta\to0} \mathbb{E}_0 \left| \frac{L_T^{Q+\theta A}-1}{\theta}
    - S_A(T) \right|^2 =0.
\]

\begin{lemma}\label{lem:common-switching-representation}
Suppose that Assumptions~\ref{ass:wellposedness},
\ref{ass:regularity-b} and
\ref{ass:regularity-sigma} hold.
Fix $\mu\in\mathcal P_2(\mathbb R^d)$
and an initial random variable $\xi$ with
$\mathcal L(\xi)=\mu$, as in Section~\ref{sec:bismut}.
There exists a jointly measurable map
\(
    \mathcal X:\mathfrak D_T\times\Omega^1
       \rightarrow C([0,T];\mathbb R^d)
\),
independent of the choice of $\widetilde Q\in\mathcal Q_{\mathcal E}$,
such that, under each $\mathbb P^{\widetilde Q}$, the process
$X_t^{\widetilde Q}:=\mathcal X_t(r,\cdot)$ solves
\eqref{eq:regularity-model} and
\[
    \mathcal L_{\mathbb P^{\widetilde Q}}
       (X_t^{\widetilde Q}\mid\mathscr D_t)(\gamma)
      =\mathcal L_{\mathbb P^1}
       (\mathcal X_t(\gamma,\cdot))
\]
for $\mathbb P_0^{\widetilde Q}$-almost every $\gamma$; the right-hand
side is $\mathscr D_t$-measurable.

For each $\phi\in L^2(\mathbb R^d\to\mathbb R^d,\mu)$,
the Bismut weight determined by
\eqref{eq:explicit-control-epsilon} has a jointly measurable
version $\mathcal W^\phi$, independent of $\widetilde Q$, such that
\[
    \mathcal W^\phi(r,\cdot)=\delta(h_0^{0,\phi})
    \quad\mathbb P^{\widetilde Q}\text{-a.s.},
    \ \widetilde Q\in\mathcal Q_{\mathcal E}.
\]
Moreover, $\phi\mapsto\mathcal W^\phi$ is linear into
$L^2(\mathbb P)$, 
\[
    \mathbb E|\mathcal W^\phi|^2
    \le C_{T,\mu,Q}\|\phi\|_{L^2_\mu}^2,
    \quad \text{and} \quad
    \mathbb E^1[\mathcal W^\phi(\gamma,\cdot)]=0
    \quad\text{for }\mathbb P_0\text{-a.e. }\gamma.
\]
\end{lemma}

\begin{proof}
For each $\ell\in\mathbb S$, let
\(
    Y_t^\ell:=\int_0^t\sigma(s,\ell)\delta B_s^H
\)
be the continuous version supplied by
Lemma~\ref{lem:deterministic_skorohod_estimate}.
For $\gamma\in\mathfrak D_T$, write $J(\gamma)$
for its number of jumps in $(0,T)$, with jump times
$0<\tau_1<\cdots<\tau_{J(\gamma)}<T$ when
$J(\gamma)\ge1$. Set $\tau_0=0$ and
$\tau_{J(\gamma)+1}=T$. Define
\[
    Z_t^\gamma
    :=\sum_{k=0}^{J(\gamma)}
       \left(
       Y_{t\wedge\tau_{k+1}}^{\gamma_{\tau_k}}
       -Y_{t\wedge\tau_k}^{\gamma_{\tau_k}}
       \right).
\]
The jump times and states are Borel functionals of
$\gamma$, so $Z^\gamma$ is jointly measurable and
non-anticipative. For every fixed $\gamma$, we have
\[
    Z_t^\gamma
    =\int_0^t\sigma(s,\gamma_s)\delta B_s^H.
\]
The uniform bound on $\sigma$ and the covariance kernel of
$B^H$ imply, uniformly in $\gamma$,
\[
    \mathbb E^1|Z_t^\gamma-Z_s^\gamma|^2
    \le C_T|t-s|^{2H},\quad s,t\in[0,T].
\]
Since $Z^\gamma$ is Gaussian and $H>1/2$,
standard Gaussian moment estimates and the continuity
criterion yield
\[
    \sup_{\gamma\in\mathfrak D_T}
    \mathbb E^1\sup_{t\le T}|Z_t^\gamma|^2
    \le C_T.
\]

We construct Picard iterates on the spatial probability
space by
\[
\begin{aligned}
    X_t^{(\gamma,0)}&=\xi+Z_t^\gamma,\quad
    \nu_t^{(\gamma,n)}=\mathcal L_{\mathbb P^1}
                         (X_t^{(\gamma,n)}),\\
    X_t^{(\gamma,n+1)}
       &=\xi+Z_t^\gamma
         +\int_0^t b(s,X_s^{(\gamma,n)},
                       \nu_s^{(\gamma,n)},\gamma_s) \mathrm ds.
\end{aligned}
\]
By induction, $(t,\gamma,\omega^1)\mapsto
X_t^{(\gamma,n)}$ is jointly measurable and
non-anticipative. Using a countable measure-determining
family of bounded continuous functions and the
second-moment maps, we also see that
$(t,\gamma)\mapsto\nu_t^{(\gamma,n)}$ is measurable
as a $\mathcal P_2(\mathbb R^d)$-valued map.
The fixed-path Picard argument in the proof of
Theorem~\ref{thm:fixed_point_wellposedness}
yields a jointly measurable limit $\mathcal X$ satisfying
\begin{equation}\label{eq:common-path-fixed-sde}
    \mathcal X_t(\gamma,\cdot)
    =\xi+Z_t^\gamma
    +\int_0^t
       b(s,\mathcal X_s(\gamma,\cdot),
         \nu_s^\gamma,\gamma_s) \mathrm ds,
\end{equation}
where
$\nu_t^\gamma=\mathcal L_{\mathbb P^1}
(\mathcal X_t(\gamma,\cdot))$.
The passage to the limit follows from the Lipschitz
condition and the coupling inequality, while uniqueness
follows from Gronwall's inequality. No step depends on a
transition rate or on the choice of $\widetilde Q$.

By non-anticipativity, $\nu_t^\gamma$ depends only on
$\gamma|_{[0,t]}$ and hence is $\mathscr D_t$-measurable.
Continuity of $\mathcal X(\gamma,\cdot)$ and the coupling
inequality imply continuity of $t\mapsto\nu_t^\gamma$.
For every bounded Borel function $g$ and
$B\in\mathscr D_t$, by Fubini's theorem we get
\[
    \mathbb E^{\widetilde Q}\left[
       \mathbf1_B(r)g(\mathcal X_t(r,\cdot))\right]
    =\int_B\mathbb E^1
       \left[g(\mathcal X_t(\gamma,\cdot))\right]
       \mathbb P_0^{\widetilde Q}(\mathrm d\gamma)
    =\mathbb E^{\widetilde Q}\left[
       \mathbf1_B(r)\int_{\mathbb R^d}
       g(y)\nu_t^r(\mathrm dy)\right].
\]
Consequently, we have $\mathcal L_{\mathbb P^{\widetilde Q}} (\mathcal X_t(r,\cdot)\mid\mathscr D_t) =\nu_t^r$ for $\widetilde Q\in\mathcal Q_{\mathcal E}$.
Together with \eqref{eq:common-path-fixed-sde}, this shows
that $X^{\widetilde Q}:=\mathcal X(r,\cdot)$ solves
\eqref{eq:regularity-model} under every $\mathbb P^{\widetilde Q}$.

Finally, fix $\phi\in L^2(\mathbb R^d\to\mathbb R^d,\mu)$.
The common solution and conditional-law representations,
uniqueness in Lemma~\ref{lem:initial-variation}, and
\eqref{eq:explicit-control-epsilon} show that the
variational processes and the control $h_0^{0,\phi}$
can be chosen independently of $\widetilde Q$.
The change of measure affects only the law of the common switching process;
the spatial It\^o representation
\eqref{eq:control-divergence-Ito} therefore remains valid
under every $\mathbb P^{\widetilde Q}$.
Since $\mathbb P_0^{\widetilde Q}\sim\mathbb P_0$,
the same jointly measurable version $\mathcal W^\phi$
represents the Bismut weight for every $\widetilde Q$.
By linearity of the variational equation and the control,
$\phi\mapsto\mathcal W^\phi$ is linear in
$L^2(\mathbb P)$. Moreover, from
\eqref{eq:control-estimate-zero} we know
\[
    \mathbb E|\mathcal W^\phi|^2
    \le C_{T,\mu,Q}\|\phi\|_{L^2_\mu}^2.
\]
For $\mathbb P_0$-almost every fixed switching path,
the spatial It\^o integral is centered, implying
\(
    \mathbb E^1[\mathcal W^\phi(\gamma,\cdot)]=0
\).
The proof is complete.
\end{proof}

\subsection{Stability of Lions derivatives}
\label{subsec:mixed-lions-sensitivity}

For $\widetilde Q\in\mathcal Q_{\mathcal E}$ and
$f\in\mathscr B_b(\mathbb R^d)$, define
\(
    (P_T^{\widetilde Q}f)(\mu)
    :=\mathbb E^{\widetilde Q}[f(X_T^{\widetilde Q})]
\).
By Lemma~\ref{lem:common-switching-representation},
$X^{\widetilde Q}=\mathcal X(r,\cdot)$ on the canonical
product space. At the reference generator $Q$, its law
agrees with that of the original model in
Section~\ref{sec:bismut}, so $P_T^Q=P_T$.

\begin{theorem}
\label{thm:mixed-lions-sensitivity}
Suppose that Assumptions~\ref{ass:wellposedness},
\ref{ass:regularity-b} and~\ref{ass:regularity-sigma}
hold. Fix $T>0$, $\mu\in\mathcal P_2(\mathbb R^d)$,
$f\in\mathscr B_b(\mathbb R^d)$ and
$Q\in\mathcal Q_{\mathcal E}$. Let $C_{T,\mu,Q}$ be
a constant in the Bismut-weight estimate of
Lemma~\ref{lem:common-switching-representation}.
Then the following statements hold.
\begin{enumerate}[(i)]
    \item For every $\widetilde Q\in\mathcal Q_{\mathcal E}$,
    \begin{equation}
    \label{eq:matrix-lions-stability}
        \big\|D^L(P_T^{\widetilde Q}f)(\mu)
                -D^L(P_T^Qf)(\mu)\big\|_{L^2_\mu}
        \le C_{T,\mu,Q}^{1/2}
        \inf_{c\in\mathbb R}\|f-c\|_\infty
        \left(e^{T d_Q(Q,\widetilde Q)^2}-1\right)^{1/2}.
    \end{equation}
    In particular, for each $R>0$, there exists
    $C_{T,\mu,Q,R}>0$ such that
    \[
        \big\|D^L(P_T^{\widetilde Q}f)(\mu)
                -D^L(P_T^Qf)(\mu)\big\|_{L^2_\mu}
        \le C_{T,\mu,Q,R}\|f\|_\infty
        \|\widetilde Q-Q\|
    \]
    whenever $\widetilde Q\in\mathcal Q_{\mathcal E}$ and
    $\|\widetilde Q-Q\|\le R$.

    \item The map
    \(
        Q\mapsto D^L(P_T^Qf)(\mu)
    \)
    is Fr\'echet differentiable on $\mathcal Q_{\mathcal E}$.
    For every $A\in\mathcal T_{\mathcal E}$ and
    $\phi\in  L^2(\mathbb R^d\to\mathbb R^d,\mu)$,
    \begin{equation}
    \label{eq:mixed-lions-bismut}
        \left\langle
        \partial_A^Q D^L(P_T^Qf)(\mu),\phi
        \right\rangle_{L^2_\mu} =
        \mathbb E^Q\left[
            f(X_T^Q)\delta(h_0^{0,\phi})S_A(T)
        \right],
    \end{equation}
    where $\partial_A^Q$ denotes the derivative
    in the matrix direction $A$.

    \item Let $q_*$ and $\rho\in(0,q_*/2)$ be as in
    Section~\ref{subsec:likelihood-ratio}. There exists
    $C_{T,\mu,Q,\rho}>0$ such that
    \begin{equation}
    \label{eq:matrix-lions-taylor}
        \big\|D^L(P_T^{Q+\Delta}f)(\mu)
        -D^L(P_T^Qf)(\mu)
        -\partial_\Delta^Q D^L(P_T^Qf)(\mu)
        \big\|_{L^2_\mu}
        \le C_{T,\mu,Q,\rho}
        \inf_{c\in\mathbb R}\|f-c\|_\infty
        \|\Delta\|^2,
    \end{equation}
    for every $\Delta\in\mathcal T_{\mathcal E}$ with
    $\|\Delta\|\le\rho$.
\end{enumerate}
\end{theorem}

\begin{proof}
By Theorem~\ref{thm:bismut-formula} applied under
$\mathbb P^{\widetilde Q}$, and by
Lemma~\ref{lem:common-switching-representation},
for every $\widetilde Q\in\mathcal Q_{\mathcal E}$ and
$\phi\in   L^2(\mathbb R^d\to\mathbb R^d,\mu)$,
\begin{equation}
\label{eq:matrix-lions-common-representation}
    D_\phi^L(P_T^{\widetilde Q}f)(\mu)
    =\mathbb E^Q\left[
       f(\mathcal X_T)\mathcal W^\phi
       L_T^{\widetilde Q}\right].
\end{equation}
Here $\mathcal X$ and $\mathcal W^\phi$ are the
common-path solution and Bismut weight from
Lemma~\ref{lem:common-switching-representation},
and
\begin{equation}
\label{eq:matrix-lions-weight-bound}
    \|\mathcal W^\phi\|_{L^2(\mathbb P^Q)}
    \le C_{T,\mu,Q}^{1/2}\|\phi\|_{L^2_\mu}.
\end{equation}
Since $\mathbb E^1[\mathcal W^\phi(\gamma,\cdot)]=0$
for $\mathbb P_0^Q$-almost every $\gamma$, and
$L_T^{\widetilde Q}$ and $S_A(T)$ depend only on the
switching path, we know
\[
    \mathbb E^Q[\mathcal W^\phi L_T^{\widetilde Q}]=0,
    \quad \text{and} \quad
    \mathbb E^Q[\mathcal W^\phi S_A(T)]=0.
\]
Consequently, $f(\mathcal X_T)$ may be replaced by
$f(\mathcal X_T)-c$ for any $c\in\mathbb R$ in the
formulas below.

To prove (i), subtract the representations
\eqref{eq:matrix-lions-common-representation} at
$\widetilde Q$ and $Q$. For every $\phi\in   L^2(\mathbb R^d\to\mathbb R^d,\mu)$,
\[
    D_\phi^L(P_T^{\widetilde Q}f)(\mu)
     -D_\phi^L(P_T^Qf)(\mu)
    =\mathbb E^Q\left[
       (f(\mathcal X_T)-c)\mathcal W^\phi
       (L_T^{\widetilde Q}-1)\right].
\]
By the Cauchy--Schwarz inequality,  \eqref{eq:matrix-lions-weight-bound}
and \eqref{eq:matrix-likelihood-l2}, we get
\[
    \big|D_\phi^L(P_T^{\widetilde Q}f)(\mu)
             -D_\phi^L(P_T^Qf)(\mu)\big|
    \le C_{T,\mu,Q}^{1/2}\|f-c\|_\infty
       \|\phi\|_{L^2_\mu}
       \left(e^{Td_Q(Q,\widetilde Q)^2}-1\right)^{1/2}.
\]
Taking the supremum over $\|\phi\|_{L^2_\mu}\le1$
and then the infimum over $c$ gives
\eqref{eq:matrix-lions-stability}. Since
\(
    d_Q(Q,\widetilde Q)^2
    \le q_*^{-1}\|\widetilde Q-Q\|^2
\),
and $e^u-1\le ue^u$ for $u\ge0$, we also have
\[
    \left(e^{Td_Q(Q,\widetilde Q)^2}-1\right)^{1/2}
    \le\sqrt{T/q_*} e^{TR^2/(2q_*)}
       \|\widetilde Q-Q\|
\]
whenever $\|\widetilde Q-Q\|\le R$. This proves (i).

We prove (ii) and (iii) together.
For $A\in\mathcal T_{\mathcal E}$, consider
the linear functional
\(
    \phi\mapsto
       \mathbb E^Q\left[
       f(\mathcal X_T)\mathcal W^\phi S_A(T)\right]
\)
on $L^2(\mathbb R^d\to\mathbb R^d,\mu)$. By \eqref{eq:matrix-lions-weight-bound},
\eqref{eq:matrix-score-second-moment}, and conditional
centering, for every $c\in\mathbb R$, we have
\[
    \big|\mathbb E^Q[
       f(\mathcal X_T)\mathcal W^\phi S_A(T)]\big|
    \le C_{T,\mu,Q}^{1/2}\|f-c\|_\infty
       \|\phi\|_{L^2_\mu}
       \big(\mathbb E_0^Q[S_A(T)^2]\big)^{1/2}.
\]
By the Riesz representation theorem, there exists
a unique $R_Q(A)\in L^2(\mathbb R^d\to\mathbb R^d,\mu)$ such that
\[
    \langle R_Q(A),\phi\rangle_{L^2_\mu}
    =\mathbb E^Q\left[
        f(\mathcal X_T)\mathcal W^\phi S_A(T)
    \right].
\]
The linearity of $A\mapsto S_A(T)$ implies
that $A\mapsto R_Q(A)$ is linear.
The preceding estimate implies
\begin{equation}
    \label{eq:mixed-lions-norm-estimate}
    \left\| R_Q(A) \right\|_{L_\mu^2} \le C_{T,\mu,Q}^{1/2}
    \inf_{c\in\mathbb R}\|f-c\|_\infty
    \left(
    \sum_{(i,j)\in\mathcal E}\frac{a_{ij}^2}{q_{ij}}
    \int_0^T\mathbb P_0^Q(r_s=i) \mathrm ds
    \right)^{1/2}.
\end{equation}
Since $q_*>0$, \eqref{eq:mixed-lions-norm-estimate}
shows that $A\mapsto R_Q(A)$ is bounded.

Fix $\rho\in(0,q_*/2)$ and
$\Delta\in\mathcal T_{\mathcal E}$ with
$\|\Delta\|\le\rho$. The same common-path
representation yields,
\[
    D_\phi^L(P_T^{Q+\Delta}f)(\mu)
     -D_\phi^L(P_T^Qf)(\mu)
     -\big\langle
       R_Q(\Delta) ,\phi
       \big\rangle_{L^2_\mu}
    =\mathbb E^Q\left[
       (f(\mathcal X_T)-c)\mathcal W^\phi
       \big(L_T^{Q+\Delta}-1-S_\Delta(T)\big)
       \right],
\]
for every $\phi\in   L^2(\mathbb R^d\to\mathbb R^d,\mu)$. 
Cauchy--Schwarz inequality, \eqref{eq:matrix-lions-weight-bound}
and the uniform Taylor estimate
\eqref{eq:matrix-likelihood-taylor} imply
\[
    \Big|D_\phi^L(P_T^{Q+\Delta}f)(\mu)
     -D_\phi^L(P_T^Qf)(\mu)
     -\big\langle
        R_Q(\Delta),\phi
       \big\rangle_{L^2_\mu}\Big|
    \le C_{T,\mu,Q}^{1/2}C_{T,Q,\rho}
       \|f-c\|_\infty
       \|\phi\|_{L^2_\mu}\|\Delta\|^2.
\]
Taking the supremum over $\|\phi\|_{L^2_\mu}\le1$
and then the infimum over $c\in\mathbb R$ yields
\[
    \big\|D^L(P_T^{Q+\Delta}f)(\mu)
       -D^L(P_T^Qf)(\mu)-R_Q(\Delta)\big\|_{L^2_\mu}
    \le C_{T,\mu,Q,\rho}
       \inf_{c\in\mathbb R}\|f-c\|_\infty
       \|\Delta\|^2.
\]
Since $A\mapsto R_Q(A)$ is bounded and linear,
this proves the Fr\'echet differentiability of
$Q\mapsto D^L(P_T^Qf)(\mu)$, with
\[
    \partial_A^Q D^L(P_T^Qf)(\mu)=R_Q(A).
\]
The preceding estimate therefore gives
\eqref{eq:matrix-lions-taylor}.
Identifying
$\mathcal W^\phi=\delta(h_0^{0,\phi})$
under $\mathbb P^Q$ gives
\eqref{eq:mixed-lions-bismut}.
This completes the proof of (ii) and (iii).
\end{proof}

\subsection{A nonlinear conditional mean-field model}
\label{subsec:nonlinear-application}

In this subsection, let $d=1$ and
$\mathbb S=\{-1,1\}$.
We use $\mathfrak D_T=D([0,T];\mathbb S)$
for the corresponding two-state Skorokhod
path space, equipped with its canonical
filtration $(\mathscr D_t)_{t\in[0,T]}$.
For $\lambda\ge0$, let $r^\lambda$ be
the continuous-time Markov chain starting from
$1$ with generator
\[
    Q_\lambda=
    \begin{pmatrix}
        -\lambda & \lambda\\
        \lambda & -\lambda
    \end{pmatrix}.
\]
Denote by $\mathbb P_0^\lambda$ the law of
$r^\lambda$ on $\mathfrak D_T$.
On this canonical path space, $r^\lambda$
is identified with the coordinate process
under $\mathbb P_0^\lambda$.
Set
\(
    \mathbb P^\lambda
    :=
    \mathbb P_0^\lambda\otimes\mathbb P^1
\).
Write $\mathbb E_0^\lambda$ and
$\mathbb E^\lambda$ for the corresponding
expectations. Throughout this subsection,
$\mu\in\mathcal P_2(\mathbb R)$ is a  
probability measure, and $\xi$ is an initial
random variable with $\mathcal L(\xi)=\mu$,
independent of the Markov chain and fractional noise.
For $a,\beta\in\mathbb R$ and $H\in(1/2,1)$, consider
\begin{equation}
\label{eq:nonlinear-two-state-model}
    \mathrm dX_t^{\lambda,a,\mu}
    =\left[
        a\tanh(X_t^{\lambda,a,\mu})
        +\beta\mathbb E^\lambda
        \left(
            X_t^{\lambda,a,\mu}
            \mid\mathscr D_t
        \right)
      \right]\mathrm dt
      +r_t^\lambda\delta B_t^H,
    \quad X_0^{\lambda,a,\mu}=\xi.
\end{equation}
Here the drift is independent of the regime, whereas the diffusion
coefficient changes sign with the regime.
Since $|\tanh'|\le1$, $|\tanh|\le1$,
and the Lions derivative of the mean functional is constant,
Assumptions~\ref{ass:wellposedness},
\ref{ass:regularity-b} and~\ref{ass:regularity-sigma}
hold on every finite time interval.

For $\gamma\in\mathfrak D_T$, let
$X^{\gamma,a,\mu}$ denote the fixed-path solution
constructed as in
Lemma~\ref{lem:common-switching-representation}.
Given $\phi\in L^2(\mathbb R\to\mathbb R,\mu)$,
let $U^{\gamma,a,\mu,\phi}$ denote the initial-law
directional variation with respect to
$\xi\mapsto\xi+\varepsilon\phi(\xi)$.
By Lemma~\ref{lem:initial-variation},
$U^{\gamma,a,\mu,\phi}$ satisfies
\begin{equation}
\label{eq:nonlinear-variational}
    U_t^{\gamma,a,\mu,\phi}
    =\phi(\xi)+\int_0^t\bigl[
        a\operatorname{sech}^2(X_s^{\gamma,a,\mu})
           U_s^{\gamma,a,\mu,\phi}
        +\beta\mathbb E^1 U_s^{\gamma,a,\mu,\phi}
      \bigr]\mathrm ds.
\end{equation}
In particular, the mean-field variation term in
\eqref{eq:epsilon-law-variation} reduces, for a fixed
switching path, to
\(
    \left.F_s^{0,\phi}\right|_{r=\gamma}
    =\beta\mathbb E^1 U_s^{\gamma,a,\mu,\phi}
\).
Accordingly, the control in
\eqref{eq:explicit-control-epsilon}, with $s_0=0$,
becomes
\begin{equation}\label{eq:nonlinear-bismut-control}
    \begin{aligned}
        (R_Hh^{\gamma,a,\mu,\phi})(t)
            & :=\frac1T\int_0^t\gamma_s
            \left[
            U_s^{\gamma,a,\mu,\phi}
            +s\beta\mathbb E^1 U_s^{\gamma,a,\mu,\phi}
            \right]\mathrm ds,\\
        \mathcal W^{\gamma,a,\mu,\phi}
            & :=\int_0^T
            \left(K_H^{-1}(R_Hh^{\gamma,a,\mu,\phi})\right)(s)
            \mathrm dW_s.
    \end{aligned}
\end{equation}
Here $\sigma^{-1}(s,\gamma_s)=\gamma_s$.
The inverse transform and the stochastic integral are understood
for almost every switching path under the relevant positive-rate law;
Lemma~\ref{lem:common-switching-representation} supplies their
jointly measurable common versions. The same lemma and
\eqref{eq:control-estimate-zero} imply that
$\phi\mapsto\mathcal W^{r^\lambda,a,\mu,\phi}$
is linear and, for every $\lambda>0$,
\begin{equation}
\label{eq:nonlinear-weight-centering}
    \mathbb E^1\mathcal W^{\gamma,a,\mu,\phi}=0
    \quad\text{for }\mathbb P_0^\lambda\text{-a.e. }\gamma,\quad \text{and} \quad
    \mathbb E^\lambda
       |\mathcal W^{r^\lambda,a,\mu,\phi}|^2
       \le C_{T,H,\lambda,a,\beta,\mu}
           \|\phi\|_{L^2_\mu}^2,
\end{equation}
where the constant is independent of $\phi$.

For $f\in\mathscr B_b(\mathbb R)$, define
\(
    (P_T^{\lambda,a}f)(\mu)
    :=\mathbb E^\lambda
       \left[f(X_T^{\lambda,a,\mu})\right]
\).
For $\lambda>0$, write
\[
    A:=\begin{pmatrix}-1&1\\1&-1\end{pmatrix},
    \quad \text{and} \quad Q_{\lambda+h}=Q_\lambda+hA.
\]
Thus $A\in\mathcal T_{\mathcal E}$ for the positive
transition-edge set
$\mathcal E=\{(1,-1),(-1,1)\}$ of $Q_\lambda$.
Let $N_T$ be the total number of jumps of $r^\lambda$
on $[0,T]$. By \eqref{eq:matrix-likelihood-score}, we have
\[
    S_A(T)
    =\frac{N_{1,-1}(T)+N_{-1,1}(T)}{\lambda}
        -O_1(T)-O_{-1}(T) =\frac{N_T}{\lambda}-T=:S_\lambda(T).
\]
Since $O_1(T)+O_{-1}(T)=T$,
\eqref{eq:matrix-score-second-moment} implies
\(
    \mathbb E_0^\lambda[S_\lambda(T)^2]
       =\frac{T}{\lambda}
\).
The likelihood ratio in
\eqref{eq:switching-likelihood-ratio} reduces to
\begin{equation}
\label{eq:nonlinear-likelihood-ratio}
    L_T^h
    :=\frac{\mathrm d\mathbb P_0^{\lambda+h}}
            {\mathrm d\mathbb P_0^\lambda}
    =\left(1+\frac{h}{\lambda}\right)^{N_T}e^{-hT},
    \quad |h|<\lambda.
\end{equation}

\begin{theorem}
\label{thm:nonlinear-threshold}
Let $H\in(1/2,1)$, $T>0$, $\lambda>0$,
$a,\beta\in\mathbb R$, and
$\mu\in\mathcal P_2(\mathbb R)$.
For $z\in\mathbb R$, set
$f_z=\mathbf1_{(z,\infty)}$.
Then $P_T^{\lambda,a}f_z$ is Lions differentiable at $\mu$,
and the mapping
\(
    \lambda\mapsto
       D^L(P_T^{\lambda,a}f_z)(\mu)
\)
is differentiable in $L^2(\mathbb R\to\mathbb R,\mu)$ for every $\lambda>0$.
For each $\phi\in L^2(\mathbb R\to\mathbb R,\mu)$,
\begin{equation}
\label{eq:nonlinear-threshold-formula}
\begin{aligned}
    D_\phi^L(P_T^{\lambda,a}f_z)(\mu)
      &=\mathbb E^\lambda\left[
          f_z(X_T^{\lambda,a,\mu})
          \mathcal W^{r^\lambda,a,\mu,\phi}\right],\\
    \partial_\lambda D_\phi^L(P_T^{\lambda,a}f_z)(\mu)
      &=\mathbb E^\lambda\left[
          f_z(X_T^{\lambda,a,\mu})
          \mathcal W^{r^\lambda,a,\mu,\phi}
          S_\lambda(T)\right].
\end{aligned}
\end{equation}
For every fixed $T,H,\lambda,a,\beta,\mu$,
there exist $\eta\in(0,\lambda)$ and $C>0$,
independent of $z$, such that
\begin{equation}
\label{eq:nonlinear-uniform-threshold}
    \sup_{z\in\mathbb R}\big\|
       D^L(P_T^{\lambda+h,a}f_z)(\mu)
        -D^L(P_T^{\lambda,a}f_z)(\mu)
       -h\partial_\lambda
                 D^L(P_T^{\lambda,a}f_z)(\mu)
    \big\|_{L^2_\mu}\le Ch^2,
    \quad |h|<\eta.
\end{equation}
Furthermore, with the constant in
\eqref{eq:nonlinear-weight-centering},
\begin{equation}
\label{eq:nonlinear-threshold-bound}
    \sup_{z\in\mathbb R}
    \left\|\partial_\lambda
        D^L(P_T^{\lambda,a}f_z)(\mu)
    \right\|_{L^2_\mu}
    \le\frac12
        C_{T,H,\lambda,a,\beta,\mu}^{1/2}
        \sqrt{\frac{T}{\lambda}}.
\end{equation}
\end{theorem}

\begin{proof}
The first identity in \eqref{eq:nonlinear-threshold-formula}
follows from Theorem~\ref{thm:bismut-formula} and
\eqref{eq:nonlinear-bismut-control}. For the
second identity, apply
Theorem~\ref{thm:mixed-lions-sensitivity}(ii)
at the reference generator $Q_\lambda$ in direction
$A$. Indeed, $Q_{\lambda+h}=Q_\lambda+hA$ and
$S_A(T)=S_\lambda(T)$. Its Fr\'echet differentiability conclusion gives the derivative in $L^2(\mathbb R\to\mathbb R,\mu)$.

Fix $\rho\in(0,\lambda/2)$ and choose
$\eta\in(0,\min\{\lambda,\rho/2\})$.
Since $\|A\|=2$ in Frobenius norm, the perturbation
$\Delta=hA$ satisfies $\|\Delta\|=2|h|<\rho$
whenever $|h|<\eta$.
Applying Theorem~\ref{thm:mixed-lions-sensitivity}(iii)
with $Q=Q_\lambda$ and $f=f_z$ yields
\[
\begin{aligned}
    & \quad \big\|D^L(P_T^{\lambda+h,a}f_z)(\mu)
          -D^L(P_T^{\lambda,a}f_z)(\mu)
          -h\partial_\lambda D^L(P_T^{\lambda,a}f_z)(\mu)
       \big\|_{L^2_\mu}\\
    & \le C_{T,\mu,Q_\lambda,\rho}
       \inf_{c\in\mathbb R}\|f_z-c\|_\infty
       \|hA\|^2
       =2C_{T,\mu,Q_\lambda,\rho}h^2.
\end{aligned}
\]
The right-hand side is independent of $z$, proving
\eqref{eq:nonlinear-uniform-threshold}.

Finally, applying \eqref{eq:mixed-lions-norm-estimate}
with $A$ and $f=f_z$, and using
\eqref{eq:nonlinear-weight-centering} and
$\inf_{c\in\mathbb R}\|f_z-c\|_\infty=1/2$,
we obtain
\[
    \left\|\partial_\lambda
       D^L(P_T^{\lambda,a}f_z)(\mu)\right\|_{L^2_\mu}
    \le\frac12 C_{T,H,\lambda,a,\beta,\mu}^{1/2}
       \big(\mathbb E_0^\lambda S_\lambda(T)^2\big)^{1/2}
    =\frac12 C_{T,H,\lambda,a,\beta,\mu}^{1/2}
       \sqrt{\frac{T}{\lambda}}.
\]
This bound is also uniform in $z$, completing the proof.
\end{proof}

For comparison, we also consider the
corresponding Brownian model with $H=1/2$,
where the stochastic integral is understood
in the It\^o sense.

\begin{theorem}
\label{thm:nonlinear-memory-persistence}
For the model \eqref{eq:nonlinear-two-state-model},
the following statements hold.
\begin{enumerate}[(i)]
    \item Fix $H=1/2$ and $T>0$.
    For every $a,\beta\in\mathbb R$ and
    $\mu\in\mathcal P_2(\mathbb R)$, the law of
    $X_T^{\lambda,a,\mu}$ is independent of
    $\lambda\ge0$. Consequently, for every
    $f\in\mathscr B_b(\mathbb R)$, whenever the
    corresponding Lions derivative exists,
    \[
        \partial_\lambda
        D^L(P_T^{\lambda,a}f)(\mu)=0
        \quad\text{in }L^2(\mathbb R\to\mathbb R,\mu),\ \lambda>0.
    \]

    \item Fix $H\in(1/2,1)$, $T>0$,
    $\beta\in\mathbb R$, $\mu\in\mathcal P_2(\mathbb R)$,
    and a nonzero direction
    $\phi\in L^2(\mathbb R\to\mathbb R,\mu)$.
    There exist $\lambda_0>0$ and $a_*>0$
    such that, for every $|a|<a_*$, one can
    choose $z_a\in\mathbb R$ satisfying
    \begin{equation}
    \label{eq:nonlinear-threshold-nonzero}
        \left.\partial_\lambda
        D_\phi^L(P_T^{\lambda,a}f_{z_a})(\mu)
        \right|_{\lambda=\lambda_0}         =\mathbb E^{\lambda_0}\left[
          f_{z_a}(X_T^{\lambda_0,a,\mu})
          \mathcal W^{r^{\lambda_0},a,\mu,\phi}
          S_{\lambda_0}(T)\right]\ne0,
    \end{equation}
    where $f_z:=\mathbf1_{(z,\infty)}$.
    In particular,
    \[
        \left.\partial_\lambda
        D^L(P_T^{\lambda,a}f_{z_a})(\mu)
        \right|_{\lambda=\lambda_0}
        \ne0\quad\text{in }L^2(\mathbb R\to\mathbb R,\mu).
    \]
\end{enumerate}
\end{theorem}

\begin{proof}
Suppose first that $H=1/2$.
For any fixed switching path $\gamma$, the process
$\widetilde B_t^\gamma:=\int_0^t\gamma_s \mathrm dB_s$
is a standard Brownian motion under $\mathbb P^1$,
since $\gamma_s^2=1$.
It is independent of $\xi$, whose law is $\mu$.
Uniqueness in law for the corresponding
McKean--Vlasov equation therefore shows that
its fixed-path solution law is independent of
$\gamma$. Hence $(P_T^{\lambda,a}f)(\mu)$
is independent of $\lambda$ for every initial
law $\mu$ and every bounded measurable $f$.
The same is true of its Lions derivative
whenever it exists, proving (i).

For (ii), fix $H\in(1/2,1)$, $T>0$,
$\beta\in\mathbb R$, $\mu\in\mathcal P_2(\mathbb R)$,
and a nonzero direction
$\phi\in L^2(\mathbb R\to\mathbb R,\mu)$.
By homogeneity, it suffices to consider
$\|\phi\|_{L^2_\mu}=1$.
By Theorem~\ref{thm:mixed-lions-sensitivity}(ii),
it suffices to establish nonzero sensitivity
for a bounded test function.
Write $m_\mu:=\int_{\mathbb R}y\mu(\mathrm dy)$, $m_\phi:=\int_{\mathbb R}\phi(y)\mu(\mathrm dy)$, $c_T:=e^{\beta T}-1$,
and set $\psi(y):=\phi(y)+c_Tm_\phi$.
If $\psi=0$ in
$L^2(\mathbb R\to\mathbb R,\mu)$, then
\[
    0=\int_{\mathbb R}\psi(y)\mu(\mathrm dy)
      =e^{\beta T}m_\phi,
\]
which implies $m_\phi=0$ and hence $\phi=0$,
a contradiction.
Thus $\psi(y)\mu(\mathrm dy)$ is a nonzero
finite signed measure.
By the uniqueness, continuity, and conjugate
symmetry of Fourier transforms of finite real
signed measures, there exists $k>0$ such that
\(
    \int_{\mathbb R}e^{iky}\psi(y)\mu(\mathrm dy)\ne0
\).
We may therefore choose $\theta\in\mathbb R$ such that
\[
    C_{\mu,\phi}(k,\theta)
    :=\int_{\mathbb R}
      \cos\bigl(k(y+c_Tm_\mu)+\theta\bigr)
      \psi(y)\mu(\mathrm dy)>0.
\]
Set $f_{k,\theta}(y):=\sin(ky+\theta)$.
We use this smooth bounded function to detect
a nonzero response, keeping $k$ and $\theta$
fixed when differentiating in $\mu$.

First take $a=0$. For each fixed switching
path $\gamma$, the conditional mean is
$\mathbb E^1X_t^{\gamma,0,\mu}=e^{\beta t}m_\mu$,
and therefore
\begin{equation}
\label{eq:nonlinear-linear-base-solution}
    X_T^{\gamma,0,\mu}
    =\xi+c_Tm_\mu
      +\int_0^T\gamma_s \delta B_s^H.
\end{equation}
The stochastic integral is centered Gaussian
with variance
\(
    V_H(\gamma):=H(2H-1)
    \int_0^T\!\int_0^T
      \gamma_s\gamma_t|t-s|^{2H-2}
       \mathrm ds \mathrm dt
\).
Set
\(
    M_{H,k}(\lambda):=\mathbb E_0^\lambda
       \left[e^{-k^2V_H(r^\lambda)/2}\right]
\).
According to the Gaussian characteristic function we have
\begin{equation}
\label{eq:nonlinear-linear-general-law}
    (P_T^{\lambda,0}f_{k,\theta})(\mu)
    =M_{H,k}(\lambda)
      \int_{\mathbb R}
      \sin\bigl(k(y+c_Tm_\mu)+\theta\bigr)
       \mu(\mathrm dy).
\end{equation}
For $\mu^\varepsilon
:=\mathcal L(\xi+\varepsilon\phi(\xi))$,
we have $m_{\mu^\varepsilon}=m_\mu+\varepsilon m_\phi$.
Using \eqref{eq:nonlinear-linear-general-law}
with $f_{k,\theta}$ fixed and differentiating
at $\varepsilon=0$, we get
\begin{equation}\label{eq:nonlinear-linear-directional-sensitivity}
    D_\phi^L(P_T^{\lambda,0}f_{k,\theta})(\mu)
    =kM_{H,k}(\lambda)C_{\mu,\phi}(k,\theta).
\end{equation}
A path with no switch has variance $T^{2H}$,
whereas a path with exactly one switch at
$s\in(0,T)$ has variance
\[
    V_H(s)=2\bigl[s^{2H}+(T-s)^{2H}\bigr]
             -T^{2H}<T^{2H},
\]
since $2H>1$.
The no-switch probability is $e^{-\lambda T}$,
the density of exactly one switch at $s$ is
$\lambda e^{-\lambda T} \mathrm ds$, and the
probability of two or more switches is
$O(\lambda^2)$ as $\lambda\downarrow0$.
Since $0\le e^{-k^2V_H(\gamma)/2}\le1$,
\[
    M_{H,k}(\lambda)
    =e^{-\lambda T}e^{-k^2T^{2H}/2}
       +\lambda e^{-\lambda T}
        \int_0^T e^{-k^2V_H(s)/2} \mathrm ds
       +O(\lambda^2)
    =M_{H,k}(0)+\lambda A_{H,T,k}
       +O(\lambda^2),
\]
where
\begin{equation}
\label{eq:nonlinear-one-switch-positive}
    A_{H,T,k}:=\int_0^T
       \left[e^{-k^2V_H(s)/2}
             -e^{-k^2T^{2H}/2}\right]\mathrm ds>0.
\end{equation}
Thus $M_{H,k}'(0+)>0$. Continuity at zero
and differentiability on $(0,\infty)$ by
likelihood-ratio differentiation imply,
by the mean-value theorem, that
$M_{H,k}'(\lambda_0)>0$ for some $\lambda_0>0$.
Consequently,
\begin{equation}
\label{eq:nonlinear-linear-positive}
    \left.\partial_\lambda
      D_\phi^L(P_T^{\lambda,0}f_{k,\theta})(\mu)
    \right|_{\lambda=\lambda_0}
    =kC_{\mu,\phi}(k,\theta)M_{H,k}'(\lambda_0)>0.
\end{equation}

We now turn to the nonlinear model.
Let $|a|\le1$.
Couple $X^{\gamma,a,\mu}$ and
$X^{\gamma,0,\mu}$ with the same initial
variable and fractional noise.
Since $|\tanh|\le1$, it follows from the difference equation
and Gronwall's inequality that
\begin{equation}
\label{eq:nonlinear-state-perturbation}
    \sup_{\gamma\in\mathfrak D_T}\sup_{t\le T}
    \|X_t^{\gamma,a,\mu}-X_t^{\gamma,0,\mu}\|
       _{L^2(\mathbb P^1)}
    \le C_{T,\beta}|a|.
\end{equation}
For $a=0$, equation
\eqref{eq:nonlinear-variational} gives
\(
    U_t^{\gamma,0,\mu,\phi}
    =\phi(\xi)+(e^{\beta t}-1)m_\phi
\).
Since $\|\phi\|_{L^2_\mu}=1$ and
$|\operatorname{sech}^2|\le1$, it follows from
Gronwall's inequality that
\[
    \sup_\gamma\sup_{t\le T}
    \|U_t^{\gamma,a,\mu,\phi}\|_{L^2(\mathbb P^1)}
    \le C_{T,\beta}.
\]
Moreover, subtracting the variational equations
for $a$ and $0$, we obtain
\[
    \|U_t^{\gamma,a,\mu,\phi}
        -U_t^{\gamma,0,\mu,\phi}\|_{L^2(\mathbb P^1)}
    \le
    |a|\int_0^t
        \|U_s^{\gamma,a,\mu,\phi}\|_{L^2(\mathbb P^1)}
        \mathrm ds
    +
    |\beta|\int_0^t
        \|U_s^{\gamma,a,\mu,\phi}
          -U_s^{\gamma,0,\mu,\phi}\|_{L^2(\mathbb P^1)}
        \mathrm ds.
\]
Consequently, another application of Gronwall's
inequality yields
\begin{equation}
\label{eq:nonlinear-variation-perturbation}
    \sup_\gamma\sup_{t\le T}
    \|U_t^{\gamma,a,\mu,\phi}
        -U_t^{\gamma,0,\mu,\phi}\|_{L^2(\mathbb P^1)}
    \le C_{T,\beta}|a|.
\end{equation}
Define
\(
    G_a(\gamma):=\mathbb E^1\left[
      f_{k,\theta}'(X_T^{\gamma,a,\mu})
      U_T^{\gamma,a,\mu,\phi}\right]
\).
Since $\|f_{k,\theta}'\|_\infty\le k$
and $\|f_{k,\theta}''\|_\infty\le k^2$,
\eqref{eq:nonlinear-state-perturbation} and
\eqref{eq:nonlinear-variation-perturbation},
together with the uniform bound on
$U^{\gamma,a,\mu,\phi}$, imply
\begin{equation}
\label{eq:nonlinear-fixed-path-perturbation}
    \sup_{\gamma\in\mathfrak D_T}|G_a(\gamma)-G_0(\gamma)|
       \le C_{T,k,\beta}|a|,
    \quad \text{and} \quad
    \sup_\gamma|G_a(\gamma)|\le C_{T,k,\beta}.
\end{equation}
For initial perturbations
$\mu^\varepsilon
=\mathcal L(\xi+\varepsilon\phi(\xi))$, by 
the Lipschitz continuity of the drift and
$\|\phi\|_{L^2_\mu}=1$ we know
\[
    \sup_\gamma\sup_{t\le T}
    \|X_t^{\gamma,a,\mu^\varepsilon}
       -X_t^{\gamma,a,\mu}\|_{L^2(\mathbb P^1)}
       \le C_{T,\beta}|\varepsilon|,
       \quad |a|\le1.
\]
By Lemma~\ref{lem:initial-variation}, the
corresponding difference quotients converge
in $L^2(\mathbb P^1)$ for each fixed $\gamma$.
The displayed uniform bound and boundedness
of $f_{k,\theta}'$ justify differentiation
under expectation with respect to the switching process. Hence
\[
    D_\phi^L(P_T^{\lambda,a}f_{k,\theta})(\mu)
    =\mathbb E_0^\lambda[G_a(r^\lambda)].
\]
Since $G_a$ is independent of $\lambda$, by virtue of 
the likelihood-ratio formula we get
\[
    \partial_\lambda
    D_\phi^L(P_T^{\lambda,a}f_{k,\theta})(\mu)
    =\mathbb E_0^\lambda
       [G_a(r^\lambda)S_\lambda(T)].
\]
Using \eqref{eq:nonlinear-fixed-path-perturbation}
and $\mathbb E_0^\lambda[S_\lambda(T)^2]=T/\lambda$,
we obtain
\begin{equation}
\label{eq:nonlinear-mixed-perturbation}
    \left|\partial_\lambda
      D_\phi^L(P_T^{\lambda,a}f_{k,\theta})(\mu)
      -\partial_\lambda
      D_\phi^L(P_T^{\lambda,0}f_{k,\theta})(\mu)
     \right| 
    \le C_{T,k,\beta}|a|
                   \sqrt{\frac{T}{\lambda}}.
\end{equation}
By \eqref{eq:nonlinear-linear-positive},
there exists $a_*>0$ such that
\begin{equation}
\label{eq:nonlinear-smooth-positive}
    \left.\partial_\lambda
      D_\phi^L(P_T^{\lambda,a}f_{k,\theta})(\mu)
    \right|_{\lambda=\lambda_0}>0
    \quad\text{for every }|a|<a_*.
\end{equation}

We now use the Bismut formula
to obtain a discontinuous test function.
For each fixed $|a|<a_*$, set
\[
    \nu_a(B):=\mathbb E^{\lambda_0}\left[
      \mathbf1_{\{X_T^{\lambda_0,a,\mu}\in B\}}
      \mathcal W^{r^{\lambda_0},a,\mu,\phi}
      S_{\lambda_0}(T)\right],
      \quad B\in\mathscr B(\mathbb R).
\]
By \eqref{eq:nonlinear-weight-centering} and
$\mathbb E_0^{\lambda_0}[S_{\lambda_0}(T)^2]
=T/\lambda_0$, the random weight is integrable.
Thus $\nu_a$ is a finite signed Borel measure.
Theorem~\ref{thm:mixed-lions-sensitivity}(ii),
together with \eqref{eq:nonlinear-bismut-control},
implies, for every $f\in\mathscr B_b(\mathbb R)$,
\begin{equation}
\label{eq:nonlinear-signed-measure-representation}
    \left.\partial_\lambda
       D_\phi^L(P_T^{\lambda,a}f)(\mu)
    \right|_{\lambda=\lambda_0}
    =\int_{\mathbb R}f(y) \nu_a(\mathrm dy).
\end{equation}
Applying this to $f=f_{k,\theta}$ and using
\eqref{eq:nonlinear-smooth-positive}, we find
$\int f_{k,\theta} \mathrm d\nu_a>0$;
therefore $\nu_a\ne0$.
If $\nu_a((z,\infty))=0$ for every $z\in\mathbb R$,
then $\nu_a(\mathbb R)=0$ by continuity from
below and
$\nu_a((x,y])=\nu_a((x,\infty))
-\nu_a((y,\infty))=0$ for all $x<y$.
Uniqueness of finite signed measures would give
$\nu_a=0$, a contradiction.
Hence there exists $z_a\in\mathbb R$ with
$\nu_a((z_a,\infty))\ne0$.
Taking $f=f_{z_a}$ in
\eqref{eq:nonlinear-signed-measure-representation}
proves \eqref{eq:nonlinear-threshold-nonzero}.
The nonzero pairing also implies the stated nonvanishing of the $L^2(\mathbb R\to\mathbb R,\mu)$-valued derivative.
The conclusion for an arbitrary nonzero direction
follows by homogeneity.
\end{proof}

\subsection{Numerical illustration}

\begin{figure}[ht]
    \centering
    \includegraphics[width= 0.8\linewidth]
    {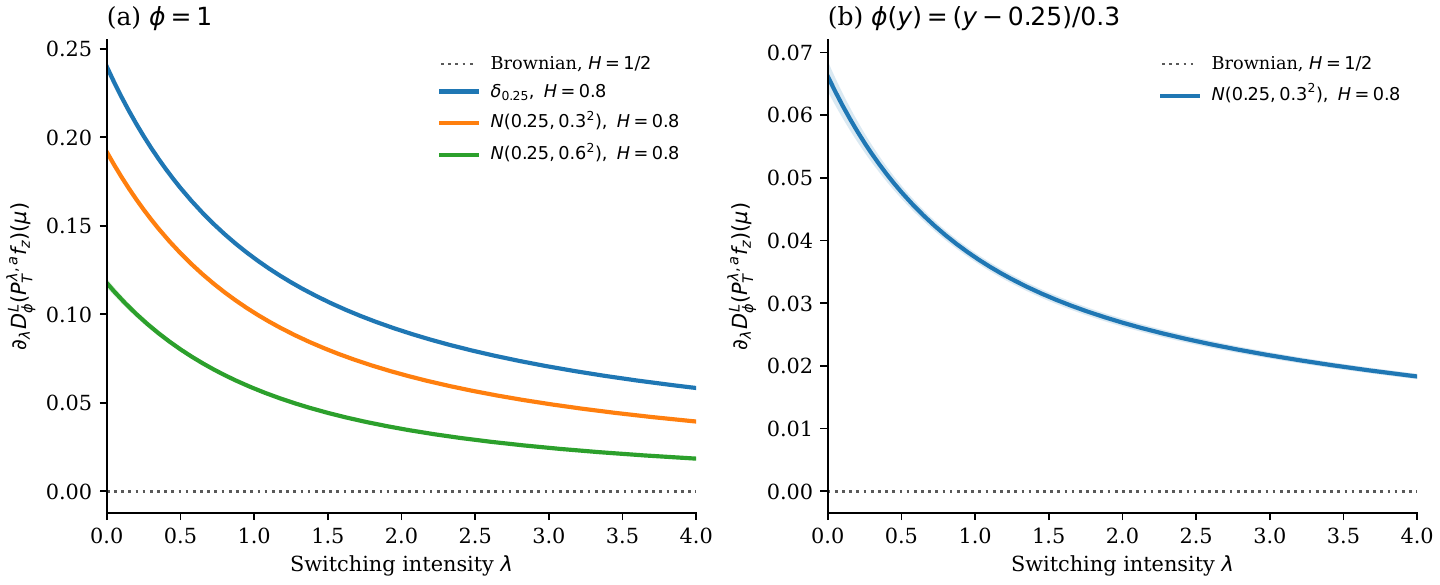}
    \caption{(a) location perturbations under different
    initial distributions;
    (b) scale perturbations of a Gaussian initial law. Shaded regions indicate pointwise 95\% Monte Carlo confidence intervals.
    }
    \label{fig:numerical-mixed-sensitivity}
\end{figure}

We illustrate the mixed transition-rate sensitivity
in Theorem~\ref{thm:nonlinear-memory-persistence}
using the nonlinear model
\eqref{eq:nonlinear-two-state-model}.
We fix $T=1$, $\beta=0.3$, $a=0.08$,
and compare fractional Brownian motion with $H=0.8$
against the Brownian case $H=1/2$.
For $f_z=\mathbf1_{(z,\infty)}$, we consider $\partial_\lambda
    D_\phi^L(P_T^{\lambda,a}f_z)(\mu)$. 
Figure~\ref{fig:numerical-mixed-sensitivity}(a)
considers $\phi\equiv1$ and $z=0.25e^{0.3}$
with
\[
    \mu=\delta_{0.25},\quad
    \mu= N(0.25,0.3^2),\quad \text{and} \quad \mu=N(0.25,0.6^2).
\]
Figure~\ref{fig:numerical-mixed-sensitivity}(b)
considers $\phi(y)=\frac{y-0.25}{0.3}$ and $z=0.25e^{0.3}+0.5$ with $\mu=N(0.25,0.3^2)$.

As shown in Figure~\ref{fig:numerical-mixed-sensitivity},
the fractional model exhibits nonzero mixed
sensitivities in both directions, in contrast
to the Brownian case.
For location perturbations, the sensitivity
decreases as the variance of the initial
distribution increases.
The nonzero scale sensitivity further shows
that the effect of fractional memory is not
restricted to translations of the initial distribution.
These numerical findings illustrate the
nonvanishing phenomenon established in
Theorem~\ref{thm:nonlinear-memory-persistence}(ii).

\appendix

\section*{Appendix}

\section{Well-posedness for the frozen equation}\label{Appendix_frozen equation}

\begin{proof}[Proof of Lemma \ref{lem:frozen_equation}]
Fix $\nu\in\mathbb{H}_{\mathscr{F}^0}^{T,p}$. By Assumption \ref{ass:wellposedness}, for $t\in[0,T]$ and $i\in\mathbb S$, we have
\[
    |b(t,x,\nu_t,i)|
    \le |b(t,0,\delta_0,i)|
    +L_T\big(|x|+W_p(\nu_t,\delta_0)\big) \le K_T+L_T|x|
    +L_T\left(\int_{\mathbb R^d}|y|^p\nu_t(\mathrm dy)\right)^{1/p},
\]
and
\begin{equation*}
    \|\sigma(t,\nu_t,i)\|
    \le K_T
        +L_T\left(\int_{\mathbb R^d}|y|^p\nu_t(\mathrm dy)\right)^{1/p}.
\end{equation*}
Consequently,
\begin{align*}
    \mathbb E^0\int_0^T
    \|\sigma(s,\nu_s,r_s)\|^p \mathrm ds
    &\le
    2^{p-1}K_T^pT
    +
    2^{p-1}L_T^p
    \int_0^T
    \mathbb E^0
    \left[
        \int_{\mathbb R^d}|y|^p\nu_s(\mathrm dy)
    \right]\mathrm ds \\
    &\le
    2^{p-1}K_T^pT
    +
    2^{p-1}L_T^pT
    \|\nu\|_{\mathbb H_{\mathscr F^0}^{T,p}}^p
    <\infty.
\end{align*}
Since $\nu_t$ and $r_t$ are $\mathscr F_t^0$-measurable, the process $\sigma(t,\nu_t,r_t)$ depends only on the common-environment variable $\omega^0$. Hence, for $\mathbb P^0$-a.e$.$ fixed $\omega^0$, the map
\(
    t\mapsto
    \sigma\left(t,\nu_t(\omega^0),r_t(\omega^0)\right)
\)
is deterministic on the spatial probability space carrying $B^H$. Therefore, the corresponding divergence integral coincides with the fractional Wiener integral. For such an $\omega^0$, set
\begin{equation*}
    Z_t^{(\nu)}(\omega^0,\cdot)
    :=
    \int_0^t
        \sigma(s,\nu_s(\omega^0),r_s(\omega^0)) \delta B_s^H.
\end{equation*}
It follows from Lemma \ref{lem:deterministic_skorohod_estimate} that $Z^{(\nu)}$ admits a continuous version and satisfies
\begin{equation*}
    \mathbb E\sup_{t\in[0,T]}|Z_t^{(\nu)}|^p
    \le
    C_{p,H}T^{pH-1}
    \mathbb E^0\int_0^T
    \|\sigma(s,\nu_s,r_s)\|^p \mathrm ds.
\end{equation*}

For fixed $\omega^0$, the frozen equation \eqref{eq:frozen_equation} can be written as the integral equation
\begin{equation*}
    X_t^{(\nu)}(\omega^0,\cdot)
    =
    \xi
    +
    \int_0^t
        b(s,X_s^{(\nu)}(\omega^0,\cdot),\nu_s(\omega^0),r_s(\omega^0)) \mathrm ds
    +
    Z_t^{(\nu)}(\omega^0,\cdot).
\end{equation*}
To construct the solution, define
\(
    X_t^{(\nu,0)}(\omega^0,\cdot):=\xi+Z_t^{(\nu)}(\omega^0,\cdot)
\),
and for $k\ge0$,
\begin{equation*}
    X_t^{(\nu,k+1)}(\omega^0,\cdot)
    :=
    \xi
    +
    \int_0^t
        b(s,X_s^{(\nu,k)}(\omega^0,\cdot),\nu_s(\omega^0),r_s(\omega^0)) \mathrm ds
    +
    Z_t^{(\nu)}(\omega^0,\cdot).
\end{equation*}
For fixed \(\omega^0\), set
\(
    a_k(t)
    :=
    \mathbb E^1\sup_{u\in[0,t]}
    |X_u^{(\nu,k)}-X_u^{(\nu,k-1)}|^p
\), for $k\ge1$.
By the preceding estimates, for $\mathbb P^0$-a.e$.$ fixed $\omega^0$, we get
\begin{equation*}
    \mathbb E^1\sup_{t\in[0,T]}
    |X_t^{(\nu,0)}|^p
    \le
    2^{p-1}
    \left(
        \mathbb E^1|\xi|^p
        +
        \mathbb E^1\sup_{t\in[0,T]}|Z_t^{(\nu)}|^p
    \right)
    <\infty.
\end{equation*}
Moreover, by H\"older's inequality and the growth estimate for $b$, we obtain
\begin{align*}
    a_1(T)
    &=
    \mathbb E^1\sup_{t\in[0,T]}
    |X_t^{(\nu,1)}-X_t^{(\nu,0)}|^p  =
    \mathbb E^1\sup_{t\in[0,T]}
    \left|
        \int_0^t
        b(s,X_s^{(\nu,0)},\nu_s,r_s) \mathrm ds
    \right|^p \\
    &\le
    T^{p-1}\int_0^T
    \mathbb E^1
    |b(s,X_s^{(\nu,0)},\nu_s,r_s)|^p \mathrm ds
    <\infty.
\end{align*}
Hence, by the Lipschitz continuity of $b$ and H\"older's inequality, we have 
\begin{equation*}
    a_{k+1}(t)
    \le
    L_T^pT^{p-1}\int_0^t a_k(s) \mathrm ds.
\end{equation*}
Iterating this inequality yields
    \( a_{k+1}(T)
    \le
    \frac{(L_T^pT^p)^k}{k!} a_1(T) \) for $k\ge1$.
Consequently, for $\mathbb P^0$-a.e$.$ fixed $\omega^0$,
\begin{align*}
    & \quad \sum_{k=1}^{\infty}
    \left\|
        X^{(\nu,k+1)}(\omega^0,\cdot)
        -
        X^{(\nu,k)}(\omega^0,\cdot)
    \right\|_{L^p(\Omega^1;C([0,T];\mathbb R^d))}\\
    & =
    \sum_{k=1}^{\infty}a_{k+1}(T)^{1/p} \le
    a_1(T)^{1/p}
    \sum_{k=1}^{\infty}
    \left(
        \frac{(L_T^pT^p)^k}{k!}
    \right)^{1/p}
    <\infty.
\end{align*}
Therefore, $\{X^{(\nu,k)}(\omega^0,\cdot)\}_{k\ge0}$ is Cauchy in
$L^p(\Omega^1;C([0,T];\mathbb R^d))$ for $\mathbb P^0$-a.e$.$
$\omega^0$. 
Let $X^{(\nu)}(\omega^0,\cdot)$ denote the limit. Then
\begin{equation*}
    X^{(\nu,k)}(\omega^0,\cdot)
    \longrightarrow
    X^{(\nu)}(\omega^0,\cdot)
    \quad\text{in }
    L^p(\Omega^1;C([0,T];\mathbb R^d))
\end{equation*}
for $\mathbb P^0$-a.e$.$ $\omega^0$. 
Taking $\mathbb E^0$ in the preceding factorial estimate
and using $\mathbb E^0a_1(T)<\infty$, we obtain
\[
    \sum_{k=1}^{\infty}
    \|X^{(\nu,k+1)}-X^{(\nu,k)}\|_{
    L^p(\Omega;C([0,T];\mathbb R^d))}
    <\infty.
\]
Consequently, the limit admits a continuous adapted
version on the product probability space.

By the Lipschitz continuity of $b$,
\begin{align*}
    & \quad \mathbb E^1\sup_{t\in[0,T]}
    \left|
        \int_0^t
        \Big[
            b(s,X_s^{(\nu,k)},\nu_s,r_s)
            -
            b(s,X_s^{(\nu)},\nu_s,r_s)
        \Big] \mathrm ds
    \right|^p \\
    & \le
    L_T^pT^{p-1}
    \int_0^T
    \mathbb E^1
    |X_s^{(\nu,k)}-X_s^{(\nu)}|^p \mathrm ds \\
    & \le
    L_T^pT^p
    \mathbb E^1\sup_{s\in[0,T]}
    |X_s^{(\nu,k)}-X_s^{(\nu)}|^p
    \longrightarrow0.
\end{align*}
Therefore, letting $k\to\infty$ in the Picard iteration gives
\begin{equation*}
    X_t^{(\nu)}
    =
    \xi
    +
    \int_0^t
        b(s,X_s^{(\nu)},\nu_s,r_s) \mathrm ds
    +
    Z_t^{(\nu)},
    \quad  t\in[0,T].
\end{equation*}
Hence $X^{(\nu)}$ is a solution of the frozen equation
\eqref{eq:frozen_equation}.

If \(X^{(\nu)}\) and \(Y^{(\nu)}\) are two solutions of the frozen equation, then the same estimate gives
\[
    \mathbb E^1\sup_{u\in[0,t]}|X_u^{(\nu)}-Y_u^{(\nu)}|^p
    \le
    C_T\int_0^t
    \mathbb E^1\sup_{v\in[0,s]}|X_v^{(\nu)}-Y_v^{(\nu)}|^p \mathrm ds.
\]
By Gronwall's inequality, we obtain \(X^{(\nu)}=Y^{(\nu)}\).

From the integral form of \eqref{eq:frozen_equation}, H\"older's inequality, and Lemma \ref{lem:deterministic_skorohod_estimate}, for $t\in[0,T]$ we obtain
\begin{align*}
    \mathbb{E}\sup_{s\in[0,t]}|X_s^{(\nu)}|^p
    &\le
    3^{p-1}\left(
        \mathbb{E}^1|\xi|^p
        +\mathbb{E}\sup_{s\in[0,t]}
        \left|
            \int_0^s b(u,X_u^{(\nu)},\nu_u,r_u) \mathrm du
        \right|^p
        +\mathbb{E}\sup_{s\in[0,t]}|Z_s^{(\nu)}|^p
    \right)\\
    &\le
    3^{p-1}\left(
        \mathbb{E}^1|\xi|^p
        +t^{p-1}\int_0^t
        \mathbb{E}|b(u,X_u^{(\nu)},\nu_u,r_u)|^p \mathrm du
        +\mathbb{E}\sup_{s\in[0,t]}|Z_s^{(\nu)}|^p
    \right)\\
    &\le
    3^{p-1}\Bigg(
        \mathbb{E}^1|\xi|^p
        +3^{p-1}t^{p-1}\int_0^t
        \Bigg[
            K_T^p
            +L_T^p\mathbb{E}|X_u^{(\nu)}|^p
            +L_T^p\mathbb{E}^0
            \left(
                \int_{\mathbb R^d}|y|^p\nu_u(\mathrm dy)
            \right)
        \Bigg]\mathrm du\\
    &\quad \quad \quad \quad
        +C_{p,H}t^{pH-1}
        \mathbb{E}^0\int_0^t
        \|\sigma(u,\nu_u,r_u)\|^p \mathrm du
    \Bigg)\\
    &\le
    3^{p-1}\Bigg(
        \mathbb{E}^1|\xi|^p
        +3^{p-1}t^{p-1}\int_0^t
        \Bigg[
            K_T^p
            +L_T^p\mathbb{E}|X_u^{(\nu)}|^p
            +L_T^p\mathbb{E}^0
            \left(
                \int_{\mathbb R^d}|y|^p\nu_u(\mathrm dy)
            \right)
        \Bigg]\mathrm du\\
    &\quad \quad \quad \quad
        +2^{p-1}C_{p,H}t^{pH-1}\int_0^t
        \Bigg[
            K_T^p
            +L_T^p\mathbb{E}^0
            \left(
                \int_{\mathbb R^d}|y|^p\nu_u(\mathrm dy)
            \right)
        \Bigg]\mathrm du
    \Bigg)\\
    &\le
    C_T\mathbb{E}^1|\xi|^p
    +C_T\int_0^t
        \mathbb{E}\sup_{v\in[0,u]}|X_v^{(\nu)}|^p \mathrm du 
    +C_T\int_0^t
        \mathbb{E}^0
        \left[
            \int_{\mathbb R^d}|y|^p\nu_u(\mathrm dy)
        \right]\mathrm du
    +C_T.
\end{align*}
Since $\nu\in\mathbb{H}_{\mathscr{F}^0}^{T,p}$, it follows from Gronwall's inequality that
\begin{equation*}
    \mathbb E\sup_{t\in[0,T]}|X_t^{(\nu)}|^p
    \le C_T\left(1+\mathbb{E}^1|\xi|^p+\|\nu\|_{\mathbb{H}_{\mathscr{F}^0}^{T,p}}^p\right).
\end{equation*}
The proof is complete.
\end{proof}

\section{Auxiliary results for the Bismut formula}

\subsection{The partial derivative and the Malliavin derivative}
\label{subsec:variational-equations}

For $\eta\in L^2(\Omega^1,\mathscr F_0^1,\mathbb P^1;\mathbb R^d)$ and
$\varepsilon\in\mathbb R$, let $X^{\varepsilon,\eta}$ denote the solution
of \eqref{eq:regularity-model} with initial value
$\xi+\varepsilon\eta$, and we set 
\(
    \mu_t^{\varepsilon,\eta}
    :=
    \mathcal L
    (X_t^{\varepsilon,\eta}\mid\mathscr F_t^0)
\).
For $\theta\in[0,1]$, define
\[
    X_t^{\varepsilon,\theta}
    :=
    X_t
    +
    \theta
    \left(
        X_t^{\varepsilon,\eta}-X_t
    \right),
    \quad \text{and} \quad
    \mu_t^{\varepsilon,\theta}
    :=
    \mathcal L
    (X_t^{\varepsilon,\theta}\mid\mathscr F_t^0).
\]

\begin{lemma}
\label{lem:initial-variation}

If Assumptions~\ref{ass:wellposedness}
and~\ref{ass:regularity-b}(i) hold, then, for
$\eta\in L^2(\Omega^1,\mathscr F_0^1,\mathbb P^1;\mathbb R^d)$,
the limit
\[
    \nabla_\eta X_t
    :=
    \lim_{\varepsilon\to0}
    \frac{X_t^{\varepsilon,\eta}-X_t}{\varepsilon}
\]
exists in
$L^2(\Omega;C([0,T];\mathbb R^d))$.
Moreover, $\nabla_\eta X$ is the unique solution of
\begin{equation}
\label{eq:initial-variation}
    \nabla_\eta X_t
    =
    \eta
    +
    \int_0^t
    \nabla_x b(s,X_s,\mu_s,r_s)
    \nabla_\eta X_s \mathrm ds
    +
    \int_0^t
    \widetilde{\mathbb E}^1
    \left[
        D^L b(s,X_s,\mu_s,r_s)
        (\widetilde X_s)
        \widetilde{\nabla_\eta X_s}
    \right]
    \mathrm ds .
\end{equation}
Furthermore, there exists a constant $C_T>0$ such that
\begin{equation}
\label{eq:initial-variation-estimate}
    \mathbb E
    \sup_{t\in[0,T]}
    |\nabla_\eta X_t|^2
    \le
    C_T\mathbb E^1|\eta|^2.
\end{equation}

\end{lemma}

\begin{proof}

For $\varepsilon\ne0$, set
\(
    \Upsilon_t^\varepsilon
    :=
    \frac{X_t^{\varepsilon,\eta}-X_t}{\varepsilon}
\).
Since $\sigma(t,r_t)$ is independent of both the state and the
distribution,
\[
    X_t^{\varepsilon,\eta}-X_t
    =
    \varepsilon\eta
    +
    \int_0^t
    \Big[
        b(s,X_s^{\varepsilon,\eta},
          \mu_s^{\varepsilon,\eta},r_s)
        -
        b(s,X_s,\mu_s,r_s)
    \Big]
    \mathrm ds .
\]
For $\mathbb P^0$-a.e.\ $\omega^0$, the pair
\(
    \left(
        X_t^{\varepsilon,\eta}(\omega^0,\cdot),
        X_t(\omega^0,\cdot)
    \right)
\)
is a coupling of
$\mu_t^{\varepsilon,\eta}(\omega^0)$ and $\mu_t(\omega^0)$.
Hence
\begin{equation}
\label{eq:conditional-W2-variation}
    W_2^2
    (\mu_t^{\varepsilon,\eta},\mu_t)
    \le
    \mathbb E^1
    |X_t^{\varepsilon,\eta}-X_t|^2,
    \quad
    \mathbb P^0\text{-a.s.}
\end{equation}
Using Assumption~\ref{ass:wellposedness}(i),
\eqref{eq:conditional-W2-variation}, and the Gronwall
argument as in \cite[Lemma~4.1]{MR4441502}, we obtain
\begin{equation}
\label{eq:upsilon-uniform-estimate}
    \sup_{0<|\varepsilon|\le1}
    \mathbb E
    \sup_{t\in[0,T]}
    |\Upsilon_t^\varepsilon|^2
    \le
    C_T\mathbb E^1|\eta|^2.
\end{equation}

We next identify the limit. 
For $\mathbb P^0$-a.e.\ fixed $\omega^0$, the conditional laws are
realized as laws under $\mathbb P^1$. Applying the ordinary
mean-value formula in the state variable and the lifted mean-value
formula for the Lions derivative
(cf.~\cite[(2.7)]{MR3630288}), we obtain
\begin{equation}
    \begin{aligned}
        \Upsilon_t^\varepsilon
    = &
    \eta
    +
    \int_0^t\int_0^1
    \nabla_x b
    \left(
        s,X_s^{\varepsilon,\theta},
        \mu_s^{\varepsilon,\eta},r_s
    \right)
    \Upsilon_s^\varepsilon
    \mathrm d\theta \mathrm ds
    \\
    & +
    \int_0^t\int_0^1
    \widetilde{\mathbb E}^1
    \left[
        D^L b
        \left(
            s,X_s,
            \mu_s^{\varepsilon,\theta},r_s
        \right)
        (\widetilde X_s^{\varepsilon,\theta})
        \widetilde{\Upsilon_s^\varepsilon}
    \right]
    \mathrm d\theta \mathrm ds .
    \end{aligned}
\label{eq:upsilon-equation}
\end{equation}
Let $\Gamma^\eta$ be the unique solution of
\[
    \Gamma_t^\eta
    =
    \eta
    +
    \int_0^t
    \nabla_x b(s,X_s,\mu_s,r_s)
    \Gamma_s^\eta
    \mathrm ds
    +
    \int_0^t
    \widetilde{\mathbb E}^1
    \left[
        D^L b(s,X_s,\mu_s,r_s)
        (\widetilde X_s)
        \widetilde{\Gamma_s^\eta}
    \right]
    \mathrm ds .
\]
By the boundedness of $\nabla_xb$ and $D^Lb$ given by Assumptions~\ref{ass:regularity-b}(i),
\begin{equation}
\label{eq:gamma-estimate}
    \mathbb E
    \sup_{t\in[0,T]}
    |\Gamma_t^\eta|^2
    \le
    C_T\mathbb E^1|\eta|^2.
\end{equation}
Moreover, \eqref{eq:upsilon-uniform-estimate} and
\eqref{eq:conditional-W2-variation} imply
\[
    \mathbb E
    \sup_{t\in[0,T]}
    |X_t^{\varepsilon,\eta}-X_t|^2
    \le
    C_T\varepsilon^2\mathbb E^1|\eta|^2,
\]
and consequently
\[
    \sup_{\theta\in[0,1]}
    \sup_{t\in[0,T]}
    \mathbb E^0
    W_2^2
    (\mu_t^{\varepsilon,\theta},\mu_t)
    \longrightarrow0,
    \quad
    \varepsilon\to0.
\]
Therefore, by the continuity and boundedness of $\nabla_xb$ and
$D^Lb$, the coefficient terms in
\eqref{eq:upsilon-equation} converge to those in the equation for
$\Gamma^\eta$. 
For $\mathbb P^0$-a.e.\ fixed $\omega^0$, the conditional
coupling estimate and the two Gronwall arguments of
\cite[Lemma~4.1]{MR4441502} also yield the conditional
counterpart of \cite[(4.6)]{MR4441502},
with $\mathbb E$ replaced by $\mathbb E^1$.
This provides the integrable domination required for the
argument of \cite[Proposition~4.2]{MR4441502}
on $(\Omega^1,\mathscr F^1,\mathbb P^1)$.
Hence, we conclude
\[
    \lim_{\varepsilon\to0}
    \mathbb E
    \sup_{t\in[0,T]}
    |\Upsilon_t^\varepsilon-\Gamma_t^\eta|^2
    =0.
\]
Here the constants in the above estimates are independent of
$\omega^0$, since the assumptions are uniform in $i\in\mathbb S$.
Hence
\(
    \nabla_\eta X=\Gamma^\eta
\)
in $L^2(\Omega;C([0,T];\mathbb R^d))$,
which proves \eqref{eq:initial-variation}.
Finally,
\eqref{eq:initial-variation-estimate} follows from
\eqref{eq:gamma-estimate}.
\end{proof}

For $\phi\in  L^2(\mathbb R^d \to \mathbb R^d , \mu)$, we shall use
$\eta=\phi(\xi)$. Since
\(
    \mu^{0,\phi}=\mu,
\) and \(
    X^{0,\phi}=X
\),
equation \eqref{eq:initial-variation} can be written as
\begin{equation*}
    \mathrm d\nabla_{\phi(\xi)}X_t
    =
    \left[
        \nabla_x b(t,X_t,\mu_t,r_t)
        \nabla_{\phi(\xi)}X_t
        +
        F_t^{0,\phi}
    \right]\mathrm dt,
    \quad
    \nabla_{\phi(\xi)}X_0
    =
    \phi(\xi),
\end{equation*}
where $F^{0,\phi}$ is defined by
\eqref{eq:epsilon-law-variation} with $\varepsilon=0$. Let $h\in L^2(\Omega;\mathcal H)$ be such that
$R_Hh$ admits an absolutely continuous version. If $X_t\in\mathbb D^{1,2}(\mathbb R^d)$ for every $t\in[0,T]$, then we set 
\(
    \mathbb D_{R_Hh}X_t
    :=
    \langle \mathbb DX_t,h\rangle_{\mathcal H}
\).

\begin{lemma}
\label{lem:malliavin-derivative}

Suppose that Assumptions~\ref{ass:wellposedness}
and~\ref{ass:regularity-b}(i) hold. Then
$X_t\in\mathbb D^{1,2}(\mathbb R^d)$ for every $t\in[0,T]$.
Moreover, 
\begin{equation}
\label{eq:malliavin-variation}
    \mathbb D_{R_Hh}X_t
    =
    \int_0^t
    \nabla_x b(s,X_s,\mu_s,r_s)
    \mathbb D_{R_Hh}X_s \mathrm ds
    +
    \int_0^t
    \sigma(s,r_s)
    \mathrm d(R_Hh)(s).
\end{equation}

\end{lemma}

\begin{proof}

Fix $\omega^0$ outside a $\mathbb P^0$-null set.
Since $\mu_s$ and $r_s$ are $\mathscr F_s^0$-measurable, for fixed
$\omega^0$ they do not depend on the $\Omega^1$-coordinate and are
therefore frozen under the Malliavin derivative with respect to $B^H$.
Thus, for each fixed $\omega^0$, we have
\[
    X_t
    =
    \xi
    +
    \int_0^t
    b
    \left(
        s,X_s,\mu_s(\omega^0),r_s(\omega^0)
    \right)
    \mathrm ds
    +
    \int_0^t
    \sigma(s,r_s(\omega^0))
    \delta B_s^H.
\]
In particular, the Malliavin derivative with respect to $B^H$
does not act on either $\mu_s$ or $r_s$. Moreover,
$\mathbb D\xi=0$, since $\xi$ belongs to the initial sigma-field
independent of $B^H$.

By Assumption~\ref{ass:regularity-b}(i),
$\nabla_x b$ is uniformly bounded. Therefore, for
$\mathbb P^0$-a.e.\ fixed $\omega^0$, the Picard argument
for additive fractional Brownian equations
(cf.~\cite[Lemma~3.1]{MR3646611}) applies on
$(\Omega^1,\mathscr F^1,\mathbb P^1)$ and yields
\(
    X_t(\omega^0,\cdot)
    \in
    \mathbb D^{1,2}(\mathbb R^d)
\) for $t\in[0,T]$.
Moreover,
\begin{equation}
\label{eq:full-malliavin-equation}
    \mathbb D X_t
    =
    \int_0^t
    \nabla_x b(s,X_s,\mu_s,r_s)
    \mathbb D X_s \mathrm ds
    +
    \sigma(\cdot,r_\cdot)\mathbf 1_{[0,t]} .
\end{equation}

Since
\(
    L^{1/H}([0,T];\mathbb R^{d\times d})
    \subset
    \mathcal H
\)
and $\sigma$ is uniformly bounded on $[0,T]\times\mathbb S$,
\[
    \sup_{\omega^0}
    \left\|
        \sigma(\cdot,r_\cdot(\omega^0))
        \mathbf 1_{[0,t]}
    \right\|_{\mathcal H}
    \le
    C_T.
\]
Using \eqref{eq:full-malliavin-equation}, the boundedness of
$\nabla_xb$, and Gronwall's inequality, we obtain
\[
    \sup_{\omega^0}
    \mathbb E^1
    \|\mathbb D X_t(\omega^0,\cdot)\|_{\mathcal H}^2
    \le
    C_T.
\]
Integration with respect to $\mathbb P^0$ therefore yields
\(
    X_t
    \in
    \mathbb D^{1,2}(\mathbb R^d)
\) for $t\in[0,T]$.
It follows from the
$\mathcal H$-inner product of
\eqref{eq:full-malliavin-equation} with $h$ that
\[
    \mathbb D_{R_Hh}X_t
    =
    \int_0^t
    \nabla_x b(s,X_s,\mu_s,r_s)
    \mathbb D_{R_Hh}X_s \mathrm ds
    +
    \left\langle
        \sigma(\cdot,r_\cdot)\mathbf 1_{[0,t]},
        h
    \right\rangle_{\mathcal H}.
\]
For each fixed common-environment path,
$\sigma(\cdot,r_\cdot)$ is deterministic with respect to the
fractional Brownian coordinate. Since $R_Hh$ is absolutely
continuous, by the reproducing-kernel identity
(cf.~\cite[Proposition~3.1]{MR3646611}) we have
\[
    \left\langle
        \sigma(\cdot,r_\cdot)\mathbf 1_{[0,t]},
        h
    \right\rangle_{\mathcal H}
    =
    \int_0^t
    \sigma(s,r_s)
    \mathrm d(R_Hh)(s).
\]
Consequently, we conclude \eqref{eq:malliavin-variation}.
\end{proof}

\subsection{Bismut formula: a general result}
\label{subsec:abstract-bismut}

We next give a general criterion for the $L$-differentiability of
$P_Tf$ and the corresponding Bismut-type representation.
Recall that $\mathscr F_{s_0}$ denotes the full filtration at time
$s_0$, namely the usual augmentation of
$\mathscr F_{s_0}^0\otimes\mathscr F_{s_0}^1$.

\begin{lemma}
\label{lem:abstract-bismut}

Suppose that the assumptions of Lemma~\ref{lem:initial-variation}
hold. Let
$\phi\in  L^2(\mathbb R^d \to \mathbb R^d , \mu)$.
Assume that, for every $\varepsilon\in[0,1]$ and
$s_0\in[0,T)$, there exists
$h_{s_0}^{\varepsilon,\phi}\in\operatorname{Dom}\delta$
such that $R_Hh_{s_0}^{\varepsilon,\phi}$ admits an absolutely
continuous version and, $\mathbb P$-a.s.,
\begin{equation}
\label{eq:abstract-matching-s0}
    \mathbb D_{R_Hh_{s_0}^{\varepsilon,\phi}}
    X_T^{\varepsilon,\phi}
    =
    \nabla_{\phi(\xi)}
    X_T^{\varepsilon,\phi}, \quad \text{and} \quad
    (R_Hh_{s_0}^{\varepsilon,\phi})(t)
    =
    0,
    \quad
    t\in[0,s_0].
\end{equation}
Assume further that, for every $s_0\in[0,T)$,
the map
\(
    \tau
    \mapsto
    \delta(h_{s_0}^{\tau,\phi})
\)
is measurable and
\begin{equation}
\label{eq:abstract-integrability}
    \int_0^1
    \left(
        \mathbb E
        \left|
            \delta(h_{s_0}^{\tau,\phi})
        \right|^2
    \right)^{1/2}
    \mathrm d\tau
    <\infty.
\end{equation}
Then, for every $\varepsilon\in[0,1]$, $s_0\in[0,T)$ and
$f\in\mathscr B_b(\mathbb R^d)$,
\begin{equation}
\label{eq:conditional-bismut-identity}
    \mathbb E
    \left[
        f(X_T^{\varepsilon,\phi})
        -
        f(X_T)
        \mid
        \mathscr F_{s_0}
    \right]
    =
    \int_0^\varepsilon
    \mathbb E
    \left[
        f(X_T^{\tau,\phi})
        \delta(h_{s_0}^{\tau,\phi})
        \mid
        \mathscr F_{s_0}
    \right]
    \mathrm d\tau.
\end{equation}
In particular,
\begin{equation}
\label{eq:unconditional-bismut-identity}
    \mathbb E
    \left[
        f(X_T^{\varepsilon,\phi})
        -
        f(X_T)
    \right]
    =
    \int_0^\varepsilon
    \mathbb E
    \left[
        f(X_T^{\tau,\phi})
        \delta(h_0^{\tau,\phi})
    \right]
    \mathrm d\tau.
\end{equation}

\end{lemma}

\begin{proof}
The proof is identical to that of
\cite[Lemma~4.5]{MR4441502} after fixing the common switching environment.
\end{proof}

We now state the general Bismut formula.

\begin{theorem}[General Bismut formula]
\label{thm:abstract-bismut}

Suppose that the assumptions of
Lemma~\ref{lem:initial-variation} hold.
Assume that, for every
$\phi\in  L^2(\mathbb R^d \to \mathbb R^d , \mu)$,
$\varepsilon\in[0,1]$ and $s_0\in[0,T)$, there exists
$h_{s_0}^{\varepsilon,\phi}\in\operatorname{Dom}\delta$
satisfying
\eqref{eq:abstract-matching-s0} and \eqref{eq:abstract-integrability}.
Assume further that
$\delta(h_0^{0,\phi})$ is $\mathscr F_T$-measurable and
\begin{equation}
\label{eq:abstract-divergence-convergence}
    \lim_{\varepsilon\downarrow0}
    \mathbb E
    \left|
        \delta(h_0^{\varepsilon,\phi})
        -
        \delta(h_0^{0,\phi})
    \right|
    =0,
    \quad
    \phi\in  L^2(\mathbb R^d \to \mathbb R^d , \mu).
\end{equation}
Then, for every $f\in\mathscr B_b(\mathbb R^d)$ and
$\phi\in  L^2(\mathbb R^d \to \mathbb R^d , \mu)$, the directional
$L$-derivative $D_\phi^L(P_Tf)(\mu)$ exists and
\begin{equation}
\label{eq:abstract-bismut-formula}
    D_\phi^L(P_Tf)(\mu)
    =
    \mathbb E
    \left[
        f(X_T)
        \delta(h_0^{0,\phi})
    \right].
\end{equation}

If, in addition, the map
\(
    \phi
    \mapsto
    \delta(h_0^{0,\phi})
\)
is linear and there exists a constant
$\widetilde C_{T,\mu}>0$ such that
\begin{equation}
\label{eq:abstract-divergence-L2}
    \mathbb E
    \left|
        \delta(h_0^{0,\phi})
    \right|^2
    \le
    \widetilde C_{T,\mu}
    \|\phi\|_{L^2_\mu} ^2,
    \quad
    \phi\in  L^2(\mathbb R^d \to \mathbb R^d , \mu),
\end{equation}
then $P_Tf$ is weakly $L$-differentiable at $\mu$.

If, moreover,
\begin{equation}
\label{eq:abstract-uniform-convergence}
    \lim_{\|\phi\|_{L^2_\mu} \to0}
    \sup_{\varepsilon\in(0,1]}
    \frac{
        \mathbb E
        \left|
            \delta(h_0^{\varepsilon,\phi})
            -
            \delta(h_0^{0,\phi})
        \right|
    }{
        \|\phi\|_{L^2_\mu} 
    }
    =0,
\end{equation}
then $P_Tf$ is $L$-differentiable at $\mu$.

\end{theorem}

\begin{proof}

The existence of the directional derivative follows from the
Bismut argument; see, for example,
\cite[Theorem~4.4]{MR4441502}. Here, we only give the details that are
needed in the present setting.

Set $Z:=\delta(h_0^{0,\phi})$.
For every fixed
$\phi\in  L^2(\mathbb R^d \to \mathbb R^d , \mu)$, we apply
Lemma~\ref{lem:abstract-bismut} as in the Bismut argument.
The only point that requires additional care is the
martingale-localization step. As $s_0\uparrow T$, we have
\(
    \mathscr F_{s_0}
    \uparrow
    \mathscr F_{T-}
    :=
    \bigvee_{s<T}\mathscr F_s
\).
Since $B^H$ has continuous sample paths, $\xi$ is
$\mathscr F_0$-measurable, and the finite-state continuous-time
Markov chain has no jump at the deterministic time $T$ almost surely,
every $\mathscr F_T$-measurable random variable considered here is
$\mathscr F_{T-}$-measurable after completion.
Hence, by the martingale convergence theorem, we have
\[
    \lim_{s_0\uparrow T}
    \mathbb E
    \left|
        \mathbb E[Z\mid\mathscr F_{s_0}]
        -Z
    \right|
    =0.
\]
By
\eqref{eq:unconditional-bismut-identity}, for
$\varepsilon\in(0,1]$, we have
\begin{align*} 
    & \quad \left|
        \frac{(P_Tf)(\mu^{\varepsilon,\phi})-(P_Tf)(\mu)}
             {\varepsilon}
        -\mathbb E[f(X_T)Z]
    \right| \\
    & = \left| \frac1\varepsilon
    \int_0^\varepsilon
    \mathbb E\left[
    f(X_T^{\tau,\phi})
    \delta(h_0^{\tau,\phi})
    -f(X_T)Z
    \right]d\tau  \right|\\
    & \le
    \frac{\|f\|_\infty}{\varepsilon}
    \int_0^\varepsilon
    \mathbb E\left|
        \delta(h_0^{\tau,\phi})-Z
    \right|\mathrm d\tau +
    \frac1\varepsilon\int_0^\varepsilon
    \left|
        \mathbb E\left[
            \left(f(X_T^{\tau,\phi})-f(X_T)\right)Z
        \right]
    \right|\mathrm d\tau.
\end{align*}
The first term tends to zero by
\eqref{eq:abstract-divergence-convergence}.
For the second term, fix $s_0\in(0,T)$.
Using \eqref{eq:conditional-bismut-identity} and
Cauchy--Schwarz inequality, we obtain
\begin{align*}
    \left|
        \mathbb E\left[
            \left(f(X_T^{\tau,\phi})-f(X_T)\right)Z
        \right]
    \right|
    & = \mathbb E\left[
    (f(X_T^{\tau,\phi})-f(X_T))
    (Z-Z_{s_0})
    \right] \\
    & \quad +
    \mathbb E\left[
    (f(X_T^{\tau,\phi})-f(X_T))
    Z_{s_0}
    \right] \\
    & \le
    2\|f\|_\infty
    \mathbb E\left|
        Z-\mathbb E[Z\mid\mathscr F_{s_0}]
    \right| \\
    & \quad +
    \|f\|_\infty\|Z\|_{L^2(\Omega)}
    \int_0^\tau
    \left(
        \mathbb E|\delta(h_{s_0}^{\theta,\phi})|^2
    \right)^{1/2}\mathrm d\theta.
\end{align*}
Consequently, we obtain by \eqref{eq:abstract-integrability} that
\begin{align*}
    & \quad \limsup_{\varepsilon\downarrow0}
    \left|
        \frac{(P_Tf)(\mu^{\varepsilon,\phi})-(P_Tf)(\mu)}
             {\varepsilon}
        -\mathbb E[f(X_T)Z]
    \right| \\
    &\le 2\|f\|_\infty
    \mathbb E\left|
        Z-\mathbb E[Z\mid\mathscr F_{s_0}]
    \right| + \limsup_{\varepsilon\downarrow0}
    \frac{\|f\|_\infty\|Z\|_{L^2(\Omega)}}{\varepsilon}
    \int_0^\varepsilon
    \int_0^\tau
    \left(
    \mathbb E|\delta(h_{s_0}^{\theta,\phi})|^2
    \right)^{1/2}
    d\theta d\tau \\
    & \le
    2\|f\|_\infty
    \mathbb E\left|
        Z-\mathbb E[Z\mid\mathscr F_{s_0}]
    \right|.
\end{align*}
Letting $s_0\uparrow T$ and using the martingale
convergence established above, we conclude that
\[
    D_\phi^L(P_Tf)(\mu)
    =\mathbb E\left[
        f(X_T)\delta(h_0^{0,\phi})
    \right],
\]
which proves \eqref{eq:abstract-bismut-formula}.

Suppose now that
$\phi\mapsto\delta(h_0^{0,\phi})$ is linear and that
\eqref{eq:abstract-divergence-L2} holds. Then
\[
    \left|
        D_\phi^L(P_Tf)(\mu)
    \right|
    \le
    \|f\|_\infty
    \left(
        \mathbb E
        \left|
            \delta(h_0^{0,\phi})
        \right|^2
    \right)^{1/2}
    \le
    \sqrt{\widetilde C_{T,\mu}} 
    \|f\|_\infty
    \|\phi\|_{L^2_\mu} .
\]
Therefore,
$\phi\mapsto D_\phi^L(P_Tf)(\mu)$ is a bounded linear functional on
$ L^2(\mathbb R^d \to \mathbb R^d , \mu)$, and hence $P_Tf$ is weakly
$L$-differentiable at $\mu$.

It remains to prove the full $L$-differentiability.
By \eqref{eq:abstract-divergence-L2} and
\eqref{eq:abstract-uniform-convergence}, we know
\[
    \sup_{\tau\in[0,1]}
    \mathbb E
    \left|
        \delta(h_0^{\tau,\phi})
    \right|
     \le  \sup_{\tau\in[0,1]}
    \mathbb E\left|
    \delta(h_0^{\tau,\phi})
    -\delta(h_0^{0,\phi})
    \right| + 
    \mathbb E\left|
    \delta(h_0^{0,\phi})
    \right|  \le 
    \left(
        o(1) + 
        \sqrt{\widetilde C_{T,\mu}}
    \right)
    \|\phi\|_{L^2_\mu},
\]
as $\|\phi\|_{L^2_\mu} \to0$.
Applying \eqref{eq:unconditional-bismut-identity} to indicator
functions, we get
\begin{equation}
\label{eq:abstract-tv-continuity}
    \sup_{\tau\in[0,1]}
    \left\|
        \mathcal L
        \left(
            X_T^{\tau,\phi}
        \right)
        -
        \mathcal L(X_T)
    \right\|_{\mathrm{TV}}
    \longrightarrow0
    \quad
    \text{as }
    \|\phi\|_{L^2_\mu} \to0.
\end{equation}
Moreover, by the stability estimate for
\eqref{eq:regularity-model}  we have
\begin{equation}
\label{eq:abstract-solution-stability}
    \sup_{\tau\in[0,1]}
    \mathbb E
    \left|
        X_T^{\tau,\phi}
        -
        X_T
    \right|^2
    \le
    C_T
    \|\phi\|_{L^2_\mu} ^2.
\end{equation}
Fix $\tau\in[0,1]$ and we define $\nu_{\tau,\phi}
    :=
    \mathcal L
    \left(
        X_T^{\tau,\phi}
    \right)$ and $\nu:=\mathcal L(X_T)$.
Let $\eta>0$.
By Lusin's theorem
(see, e.g., \cite[Theorem~7.4.4]{MR3098996}),
there exist a compact set $K\subset\mathbb R^d$ and
$g\in C_b(\mathbb R^d)$ such that
\[
    \nu(K^c)<\eta,
    \quad
    g=f \ \text{on }K,
    \quad \text{and} \quad
    \|g\|_\infty\le\|f\|_\infty.
\]
By \eqref{eq:abstract-tv-continuity}, we have
\[
    \nu_{\tau,\phi}(K^c)
    \le
    \nu(K^c)
    +
    \|\nu_{\tau,\phi}-\nu\|_{\mathrm{TV}}
    \longrightarrow
    \nu(K^c)
    \quad
    \text{as }
    \|\phi\|_{L^2_\mu} \to0.
\]
Moreover, \eqref{eq:abstract-solution-stability} implies
\[
    X_T^{\tau,\phi}
    \longrightarrow
    X_T
    \quad\text{in probability},
\]
and hence, since $g$ is bounded and continuous,
\[
    \mathbb E
    \left|
        g(X_T^{\tau,\phi})
        -
        g(X_T)
    \right|^2
    \longrightarrow0.
\]
Therefore,
\begin{align*}
    \mathbb E
    \left|
        f(X_T^{\tau,\phi})
        -
        f(X_T)
    \right|^2
    & \le
    3\mathbb E
    \left|
        f(X_T^{\tau,\phi})
        -
        g(X_T^{\tau,\phi})
    \right|^2
    +
    3\mathbb E
    \left|
        g(X_T^{\tau,\phi})
        -
        g(X_T)
    \right|^2
    \\
    & \quad +
    3\mathbb E
    \left|
        g(X_T)
        -
        f(X_T)
    \right|^2.
\end{align*}
According to Lusin's theorem, the first and third terms on the right-hand side satisfy
\[
    3\mathbb E
    \left|
        f(X_T^{\tau,\phi})
        -
        g(X_T^{\tau,\phi})
    \right|^2
    +
    3\mathbb E
    \left|
        g(X_T)
        -
        f(X_T)
    \right|^2
    \le
    12\|f\|_\infty^2
    \left(
        \nu_{\tau,\phi}(K^c)
        +
        \nu(K^c)
    \right).
\]
Letting first
$\|\phi\|_{L^2_\mu} \to0$ and then $\eta\downarrow0$, we obtain
\begin{equation}
\label{eq:abstract-bounded-measurable-continuity}
    \lim_{\|\phi\|_{L^2_\mu} \to0}
    \mathbb E
    \left|
        f(X_T^{\tau,\phi})
        -
        f(X_T)
    \right|^2
    =0.
\end{equation}

Finally, by \eqref{eq:unconditional-bismut-identity}, we obtain
\begin{align*}
    &\quad
    \frac{
        \left|
            (P_Tf)
            \left(
                \mu\circ(\operatorname{Id}+\phi)^{-1}
            \right)
            -
            (P_Tf)(\mu)
            -
            D_\phi^L(P_Tf)(\mu)
        \right|
    }{
        \|\phi\|_{L^2_\mu} 
    }
    \\
    &\le
    \|f\|_\infty
    \sup_{\tau\in(0,1]}
    \frac{
        \mathbb E
        \left|
            \delta(h_0^{\tau,\phi})
            -
            \delta(h_0^{0,\phi})
        \right|
    }{
        \|\phi\|_{L^2_\mu} 
    }
    +
    \sqrt{\widetilde C_{T,\mu}}
    \int_0^1
    \left(
        \mathbb E
        \left|
            f(X_T^{\tau,\phi})
            -
            f(X_T)
        \right|^2
    \right)^{1/2}
    \mathrm d\tau.
\end{align*}
The first term tends to zero by
\eqref{eq:abstract-uniform-convergence}, while the second tends to
zero by
\eqref{eq:abstract-bounded-measurable-continuity} and dominated
convergence. Therefore,
\[
    \lim_{\|\phi\|_{L^2_\mu} \to0}
    \frac{
        \left|
            (P_Tf)
            \left(
                \mu\circ(\operatorname{Id}+\phi)^{-1}
            \right)
            -
            (P_Tf)(\mu)
            -
            D_\phi^L(P_Tf)(\mu)
        \right|
    }{
        \|\phi\|_{L^2_\mu} 
    }
    =0.
\]
This proves the $L$-differentiability.
\end{proof}

\footnotesize
\bibliographystyle{plain}
\bibliography{reference}
\end{document}